\documentclass[10pt,reqno]{amsart}

\usepackage{amssymb}
\usepackage{amsfonts}
\usepackage{amsthm}
\usepackage{amsmath}
\usepackage{bbm}
\usepackage{MnSymbol}
\usepackage{enumitem}
\usepackage{xcolor}

\numberwithin{equation}{section}
\allowdisplaybreaks
\usepackage[style=alphabetic, sorting=nty, minnames=99, url=false,isbn=false, doi = false, maxnames=99]{biblatex}

\newtheorem{theorem}{Theorem}[section]
\newtheorem{proposition}[theorem]{Proposition}

\newtheorem{lemma}[theorem]{Lemma}
\newtheorem{corollary}[theorem]{Corollary}

\theoremstyle{definition}
\newtheorem{definition}[theorem]{Definition}

\theoremstyle{remark}
\newtheorem{remark}[theorem]{Remark}

\newcommand{\R}{\mathbb{R}}

\newcommand{\N}{\mathbb{N}}

\newcommand{\C}{\mathbb{C}}

\renewcommand{\hat}{\widehat}
\newcommand{\eps}{\varepsilon}
\newcommand{\grad}{\nabla}

\newcommand{\scriptE}{\mathcal{E}}

\newcommand{\scriptR}{\mathcal{R}}

\newcommand{\scriptT}{\mathcal{T}}

\newcommand{\jp}[1]{\langle{#1}\rangle}
\newcommand{\qtq}[1]{\quad\text{#1}\quad}
\newcommand{\ctc}[1]{\,\,\text{#1}\,\,}

\DeclareMathOperator*{\diam}{diam}
\DeclareMathOperator*{\dist}{dist}

\begin{document}
\title[Fourier restriction to hyperbolic surfaces in $\R^3$]{End-line Fourier restriction for hyperbolic surfaces in three dimensions}
\author{Sam Craig}
\address{Department of Mathematics, University of Wisconsin--Madison}
\email{secraig2@wisc.edu}
\author{Betsy Stovall}
\address{Department of Mathematics, University of Wisconsin--Madison}
\email{stovall@math.wisc.edu}

\begin{abstract}
We prove new bilinear adjoint Fourier restriction estimates for pairs of separated rectangles associated to a hyperbolic surface in $\R^3$.  We then employ these bilinear estimates to prove Fourier extension inequalities along the (sharp) parabolic scaling line for any smooth hyperbolic 2-surface that is either the graph of a polynomial or that obeys a mild finite type condition.  
\end{abstract}

\maketitle

\section{Introduction}

Define the Fourier extension operator associated with a surface $\Sigma \subset \R^3$ equipped with surface measure $d\sigma$ by $\mathcal{E}_{\Sigma}:f \mapsto \widehat{f ~d\sigma}$. We consider the problem of bounding the Fourier extension operator for a collection of hyperbolic surfaces $\Sigma$ satisfying a mild finite-type assumption. Going forward, we refer to these surface as \emph{admissible} and will give a precise definition in Section 2. They include all analytic hyperbolic surfaces in $\R^3$. 

The extension operator $\mathcal{E}_{\Sigma}$ is conjectured to obey $L^q \to L^p$ bounds when $p > 3, 2q' \le p$. In this paper, we prove these bounds for $p > \frac{22}{7}$ on the end-line $q = (p/2)'$.

\begin{theorem}\label{T:1}
    For any admissible surface $\Sigma$ and any $p > \frac{22}{7}$, there is a constant $C_{\Sigma, p}$ such that for any Schwartz function $f:\R^2 \to \C$,
    \begin{equation}\label{E:T1}
        ||\mathcal{E}_{\Sigma}f||_{L^{p}(\R^3)} \le C_{\Sigma, p} ||f||_{L^{(p/2)'}(\Sigma)}.
    \end{equation}
    
\end{theorem}

Outside of the hyperbolic paraboloid \cite{S19} and the hyperbolic hyperboloid \cite{BOS21}, Theorem \ref{T:1} is the first end-line extension estimate beyond the classical $L^2 \to L^4$ estimate for the extension operator. Our lower bound on $p$ arises from a restriction estimate for $p > \frac{22}{7}, 2q' < p$ by Guo, Liu, and Xi earlier this year \cite{GLX26}. Indeed, our approach implies that any restriction estimate for admissible surfaces up to the end-line can be upgraded to a restriction estimate on the end-line.

\begin{theorem}\label{T:2}
    Suppose $\Sigma$ is an admissible surface. Suppose for some $p_0 > 3$ we have \begin{equation}\label{E:T2 assumption}
        ||\mathcal{E}_{\Sigma}f||_{L^p(\R^3)} \le C_{\Sigma, p,q}||f||_{L^q(\Sigma)} \text{ for all }f \in \mathcal{S}(\R^3), p < p_0, \text{ and } 2q' < p.
    \end{equation}
    Then for any $p > p_0$, we also have \begin{equation}\label{E:T2 goal}
        ||\mathcal{E}_{\Sigma}f||_{L^{p}(\R^3)} \le C_{\Sigma, p} ||f||_{L^{(p/2)'}(\Sigma)}.
    \end{equation}
\end{theorem}

Our proof of these linear extension estimates works in two steps. The first step is a novel bilinear extension estimate which may be of independent interest. This extension estimate applies under weaker separation conditions than classical bilinear extension estimates. We give an informal description of the separation condition and defer the set-up required for the formal definition to the next section. Namely, a pair of points $(\xi,\eta) \in \Sigma$ will be said to be separated when $\xi$ and $\eta$ are not equal and do not lie on a common line in $\Sigma$.  In particular, in a surface containing no lines, all pairs of distinct points are separated.  We will quantify this separation with two parameters, one equal to the distance between the points, and the other equal to the maximum distance between the segment connecting $(\psi(\xi),\xi)$ and $(\psi(\eta),\eta)$ and $\Sigma$, normalized by the distance.  

This separation may be placed in a more geometric framework, namely, for a pair $\xi, \eta \in \Sigma$, let $v_{\xi, \eta} \in T_{\xi} \Sigma \cap T_{\eta} \Sigma$ be a unit vector. If $\xi \neq \eta$, then $v$ is unique up to signs. Let $\gamma:[0,1] \to \Sigma$ parametrize the curve $\{\eta' : v_{\xi, \eta} = v_{\xi, \eta'}\}$ between $\xi$ and $\eta$ and take $\mathrm{I\!I}_{\gamma(t)}:T_{\gamma(t)} \Sigma \times T_{\gamma(t)}\Sigma \to \R$ be the second fundamental form of $\Sigma$ at $\gamma(t)$. We say a pair of points $(\xi, \eta)$ is $(1,1)$-separated if $|\xi - \eta| \sim 1$ and \begin{equation}\label{E:weak condition}
    \sup_{t \in [0,1]} |\mathrm{I\!I}_{\gamma(t)}(v_{\xi, \gamma(t)}, v_{\xi, \gamma(t)})| \sim 1.
\end{equation}Our bilinear extension theorem follows. \begin{theorem}\label{T:3}
    For any admissible $\Sigma$ and any $s \in (\frac{11}{7}, 2]$, there exists $\varepsilon > 0$ such that the following holds for any $r > s' - \varepsilon$. Suppose $R_1, R_2\subset \R^3$ are two balls of radius $\sim 1$ such that for each pair $(\xi, \eta) \in (R_1 \cap \Sigma) \times (R_2 \cap \Sigma)$, $(\xi, \eta)$ are $(1,1)$-separated. Then 

    \begin{equation}\label{E:bil ext first statement}
        ||(\mathcal{E}_{\Sigma\cap R_1}f)(\mathcal{E}_{\Sigma \cap R_2}g)||_{L^s} \le C_{\Sigma, r}||f||_{L^r(\Sigma \cap R_1)}||g||_{L^r(\Sigma \cap R_2)}.
    \end{equation}Under the same assumptions as Theorem \ref{T:2}, we may improve our range of $s$ to $s > \frac{p_0}{2}$.
\end{theorem}

For comparison, Lee proved (\ref{E:bil ext first statement}) for a wider range of $r$ when each pair $(\xi, \eta)$ in the product of the domains satisfy \begin{equation}\label{E:uniform separation}
    |\xi - \eta| \sim 1, |\mathrm{I\!I}_{\xi}(v_{\xi, \eta}, v_{\xi, \eta})| \sim |\mathrm{I\!I}_{\eta}(v_{\xi, \eta}, v_{\xi, \eta})| \sim 1,
\end{equation}which is a stronger assumption than (\ref{E:weak condition}) \cite{L05}. 

Our weaker separation condition has the advantage that it produces rectangular pairs of separated sets.  By contrast, according to the condition determined in \cite{L05}, the set of points separated from a singleton $\xi \in \Sigma$ might (roughly) be a curved rectangle, and the sets of points separated from distinct elements of the curved rectangle might also be curved rectangles that are quite different from one another.  (This issue is explored in detail in \cite{BMV22}.)  Thus, in the second step, when we prove that Theorem \ref{T:3} implies Theorem \ref{T:1}, we may build upon the sharp bilinear-to-linear argument of \cite{V05, S19}, though the differences in geometry introduce significant new challenges.

\subsection{Background and motivation}

New $L^p$ bounds for the extension operator $\mathcal{E}_{\Sigma}$ imply improvements in the restriction problem for a compact surface $\Sigma$. The restriction problem asks for the range of $p, q$ for which the Fourier transform of $f \in L^p(\R^d)$ is well-defined as a function in $L^q(\Sigma)$. This problem originated from work of Stein in the 1960s and received sustained interest in the intervening years. Curvature plays an important role in the range of $(p,q)$ for which restriction estimates are conjectured; we expect the widest range of restriction estimates to hold when the principal curvatures of $\Sigma$ are uniformly non-vanishing. While the sign of the principal curvatures do not impact the conjectured range for the restriction problem, they do play an important role in proving bounds in the restriction problem. The case of hyperbolic surfaces (those with principal curvatures of different signs) tends to be more difficult than elliptic surfaces (those with principal curvatures of the same sign). In this article, we prove new restriction estimates for hyperbolic surfaces in $\R^3$. 

The conjectured range of exponents $p,q$ in the restriction problem would follow from a uniform bound $L^p(\R^d) \to L^q(\Sigma)$ for the operator $\mathcal{R}_{\Sigma}: f \mapsto 1_{\Sigma}\hat{f}$ for the same exponents. It is more common to prove equivalent bounds for the formal adjoint $\mathcal{E}_{\Sigma}$. We call $L^p$ bounds for $\mathcal{E}_{\Sigma}$ extension estimates. In $\R^3$, the restriction conjecture is equivalent to proving that $\mathcal{E}_{\Sigma}:L^q \to L^p$ for $p > 3, 2q' \le p$. Exponents pairs on the end-line $2q' = p$ are special for two reasons: standard techniques for the restriction problem cannot reach the end-line, and end-line restriction estimates are stable under rescaling, thereby enabling some uniformity. 
As mentioned previously, other than the hyperbolic paraboloid \cite{S19} and hyperbolic hyperboloid \cite{BOS21}, no end-line extension estimates are known beyond the classical $L^2 \to L^4$ bound \cite{T75}.  In fact, our separation condition sheds some light on why scale-invariant estimates had only been established for doubly-ruled surfaces.  Namely, it is precisely in this case that our separation condition matches that of Lee, and, thus, Lee's condition naturally results in separated pairs of rectangles.  

Research into extension estimates off of the scaling line has seen steady improvement over the last twenty years. For most of that time, researchers focused on the the hyperbolic paraboloid. Lee and Vargas independently proved (\ref{E:T2 assumption}) for $p > \frac{10}{3}$, $2q' > p$ around 2005 using a bilinear approach \cite{V05, L05}. Vargas's bilinear-to-linear reduction will play an important role in proving end-line restriction estimates in this paper. In 2017, Cho and Lee proved (\ref{E:T2 assumption}) for $p > \frac{13}{4}$ and $q = p$ using polynomial partitioning, which was later improved by Kim to $2q' > p$ \cite{CL17, K17}. In the following years, Buschenhenke, M\"uller, and Vargas and Guo and Oh independently worked on proving linear restriction estimates for more general hyperbolic surfaces, eventually extending Kim's result to apply for all hyperbolic surfaces in $\R^3$ \cite{BMV23, GO24}. In 2024, Demeter and Wu used refined decoupling and two-ends Furstenberg estimates to prove (\ref{E:T2 assumption}) for $p > \frac{22}{7}$ and $q=p$ for the hyperbolic paraboloid \cite{DW26}. In 2026, Guo, Liu, and Xi extended this argument to apply for all hyperbolic surfaces in $\R^3$ and for $2q' > p$ \cite{GLX26}.

Bilinear extensions estimates and their relatives play an important part in all of these linear results. In chronological order, Lee and Vargas used a variant Theorem \ref{T:3} under the stronger separation conditions (\ref{E:uniform separation}), Cho and Lee and Guo and Oh used a bilinear reverse square function estimate for ``tangential'' wave packets, and Demeter and Wu and Guo, Liu, and Xi used bilinear refined decoupling estimates. All these bilinear estimates apply under the separation condition (\ref{E:uniform separation}). While the first element $|\xi - \eta|\sim 1$ is easy enough to account for, the lower bounds on the second fundamental form are more difficult to manage. In the un-ruled case, it not clear how to completely avoid loss in the bilinear estimates due to the complicated nature of the sets where (\ref{E:uniform separation}) fails. Past literature prove a restriction estimate localized to $R$-balls, which can accept the loss from (\ref{E:uniform separation}) in an $R^{\varepsilon}$ coefficient. These losses make it impossible to prove an end-line extension estimate. While we also use a bilinear extension estimate in this paper, our argument is of a substantially different nature than recent research on restriction to hyperbolic surfaces: our extension estimate holds under a worse range of exponents but requires a weaker separation condition that does not have the same losses that would appear with a uniform separation assumption like (\ref{E:uniform separation}).

\subsection{Outline of paper}

\begin{itemize}
    \item In Section 2, we give formal definitions of admissible surfaces, define some important vector fields associated with our surface, a formal statement of Theorem \ref{T:3}, and establish some preliminary results on the geometry of admissible surfaces and on rescaling for such surfaces.
    \item In Section 3, we prove Theorem \ref{T:3}. We focus on proving an $L^s \times L^s \to L^2$ bound for some $s < 2$, since the general case follows by interpolation. The $L^2$ bound in the output space allows us to use Plancherel on $\mathcal{E}_{\Sigma \cap R_1}f\mathcal{E}_{\Sigma \cap R_2} g$ and then dualize to reduce the proof of Theorem \ref{T:3} to proving a bound for a generalized Radon transform. We use Christ's method of refinements and a careful analysis of those sets which extremize the generalized Radon transform to prove a bound of the form $L^{2 - \varepsilon}\times L^{2-\varepsilon} \to L^2$. Interpolating with linear restriction estimates yields Theorem \ref{T:3}.
    \item In Section 4, we prove a restricted strong-type $((p/2)', p)$ variant of Theorem \ref{T:1}, with the added restriction that we only consider sets $\Omega$ with essentially constant length intersection with integral curves for a null vector field for $\Sigma$. We do so by generalizing an argument of Vargas for the hyperbolic paraboloid \cite{V05}. The novel part in this part of the paper relative to Vargas's result is finding a Whitney decomposition for $\Sigma \times \Sigma$ adapted to the geometry of $\Sigma$, which is rather more difficult for un-ruled $\Sigma$ compared to ruled $\Sigma$. We do most of the technical work for the Whitney decomposition in Section 2 and finish the construction in Section 4.
    \item In Section 5, we follow an argument of the second named author to upgrade the restricted strong-type inequality for sets of constant fiber length to a restricted strong-type inequality for general sets \cite{S19}. We prove that those constant fiber-length sets $\Omega$ which extremize the restricted strong type inequality from Section 4 must be a union of a small number of ``rectangles'' on the surface. This special structure and Theorem \ref{T:3} allows us to remove the fiber length restriction in our restricted strong-type inequality. We prove Theorem \ref{T:1} by interpolating the restricted strong-type $((p/2)', p)$ inequality with the classical $L^2 \to L^4$ bound for the extension operator.
\end{itemize}

\subsection{Notation} 
We use three notations for multiplicative constants arising in different constants. In all notations, our multiplicative constant may depend on a collection of other parameters. In this section, we use $\alpha$ as a short-hand for those parameters. 
\begin{itemize}
    \item We write $A \lesssim_{\alpha} B$ if $A \le C B$ for some constant $C$ depending only on the parameters $\alpha$. If $A \lesssim_{\alpha} B$ and $B \lesssim_{\alpha} A$, we write $A \sim_{\alpha} B$. 
    \item For exponents $j, k$, we write $k \ge j + C_{\alpha}$ if there exists a constant $C$ depending only on $\alpha$ for which $k \ge j +C$. In particular, the value of $C_{\alpha}$ may change from line to line. This notation will only be used when $j, k$ are exponents, so $C_{\alpha}$ will in fact correspond to a multiplicative error.
    \item For $f:\N \to \R^+$, we write $\Theta_{\alpha}(f)$ to denote a function $g: \N \to \R^+$ satisfying $C^{-1}f(n) \le g(n) \le Cf(n)$ for all $n \in \N$.
\end{itemize}

We use $\delta, \varepsilon, \rho$ to denote small constants. We can always require these constants to be dyadic. We use $j, k$ to denote scale parameters and we use $n, m, \ell$ to denote parameters related to the type condition.

\subsection{Acknowledgements}  The first author was supported by NSF grants DMS-2037851 and DMS-2246906.  The second author was supported by NSF DMS-2246906 and and Simons Foundation grant  SFI-MPS-SFM-00011865.

\subsection*{AI Statement} We didn't use any.  Formulation of the main and auxillary theorems, method of proof, and exposition are the authors', as is each word of the text, as well as any errors grammatical, mathematical, bibliographical, or otherwise, appearing therein.

\section{Separated rectangles}

\subsection{Admissible surfaces}
We let $Q$ denote the unit square, $Q:=[-1,1]^2$, and we fix a function $\psi$ that is smooth with 
\begin{equation} \label{E:psi near psi0}
\psi(0) = 0,\: \nabla \psi(0) = 0,\: D^2\psi(0) = D^2 \psi_0(0), \qtq{and}\|D^2\psi - D^2\psi_0\|_{C^0(2Q)} < c_0,
\end{equation}
for $c_0$ a sufficiently small constant.  Here we set the notation $\psi_0(\xi):=\xi_1\xi_2$.  We define 
\begin{equation}\label{E:Sigma}
\Sigma:=\{(\psi(\xi),\xi) : \xi \in 2Q\}. 
\end{equation}
Extension estimates for these surfaces $\Sigma$ are equivalent to $L^p$ bounds for the following operator: \[\tilde{\mathcal{E}}_{\Sigma}(x, t) = \int_Q e^{-i (x \cdot \xi + t\psi(\xi))}f(\xi)~d\xi.\]As we only consider $\tilde{\mathcal{E}}_{\Sigma}$ going forward, we define $\mathcal{E} = \tilde{\mathcal{E}}_{\Sigma}$. If we need to emphasize the role of $\psi$ in the operator, we instead write $\mathcal{E}_{\psi}$. 

We note that, via a partition of unity and composition with an affine transformation, restriction estimates for any smooth, precompact hyperbolic 2-surface in $\R^3$ may be obtained from restriction estimates for such a surface. 
We recall the null directions $N_j$ defined in the introduction.  By \eqref{E:psi near psi0} we may choose $N_1$ and $N_2$ to be near the coordinate directions,
\begin{equation} \label{E:N_j near e_j}
    \|N_j-e_j\|_{C^0(2Q)} < c_0, \qquad j=1,2.
\end{equation}

It will be convenient later to observe that under the parabolic rescaling in which $\psi(\xi)$ is replaced by $\lambda^{-2}\psi(\lambda \xi)$, for $0 < \lambda < 1$, $N_j(\xi)$ is replaced by $N_j(\lambda\xi)$, and the conditions \eqref{E:psi near psi0} and \eqref{E:N_j near e_j} are not only preserved, but $c_0$ can be replaced by a smaller value.  

We denote by $\tilde N_j(\psi(\xi),\xi):=\bigl( \nabla \psi(\xi) \cdot N_j(\xi),N_j(\xi)\bigr)$ the pushfoward via the parametrization $\xi \mapsto (\psi(\xi),\xi)$ of $N_j$ to a vector field on $\Sigma$.  We will often abuse notation by writing $\tilde N_j(\psi(\xi),\xi) = \tilde N_j(\xi)$.  

\begin{definition}\label{D:admissible}
We say that $\psi:2Q \to \R$ (or the surface $\Sigma$ given by \eqref{E:Sigma}) is admissible if $\psi$ is a smooth function obeying \eqref{E:psi near psi0}, and each $N_j$, j=1,2, satisfies one of the following properties: 
\begin{itemize}
\item $\Sigma$ is ruled in the direction $\tilde N_j$, in the sense that for each $\xi,t$ such that $\xi, \xi+tN_j(\xi) \in 2Q$, $\psi(\xi+tN_j(\xi)) = \psi(\xi) + tN_j(\xi) \cdot \nabla \psi(\xi)$; or 
\item $\Sigma$ is of \textit{finite line type in the direction $\tilde N_j$}, in the sense that there exists some $m_j > 2$ such that for each $t,\xi$ such that $\xi,\xi+tN_j(\xi) \in 2Q$, 
\begin{equation} \label{E:finite type}
|\partial_t^{m_j}\psi(\xi+tN_j(\xi))| \sim c_\psi \neq 0; \qtq{or}
\end{equation}
\item $\Sigma$ is \textit{nearly of finite line type in the direction $\tilde N_j$}, in the sense that it contains the $\xi_j$-axis and no other line in the direction $\tilde N_j$, and there exists $\ell_j > 0$ such that 
\begin{equation} \label{E:nearly finite type}
\partial_{\xi_{j'}}^{\ell_j}|_{\xi_{j'}=0}\partial_{t}^2 \psi(\xi+tN_j(\xi)) \not\equiv 0.
\end{equation}
Moreover, with $\ell_j$ the minimum such exponent, there exists $m_j > 2$ such that for $\xi,t$ with $\xi,\xi+tN_j(\xi) \in 2Q$,
\begin{equation} \label{E:admissible}
|\partial_{\xi_{j'}}^{\ell_j}\partial_t^{m_j}\psi(\xi+tN_j(\xi))| \sim c_\psi \neq 0. 
\end{equation}
Here $j'$ is the complementary index to $j$, so that $\{j,j'\} = \{1,2\}$.
\end{itemize}
\end{definition}

We will say, more generally, that $\psi:2Q \to \R$ is \textit{admissible} if there exist a finite decomposition of $\Sigma$ and corresponding affine transformations of the pieces into admissible surfaces.

By minimality of $\ell_j$ in \eqref{E:nearly finite type}, and consequently \eqref{E:admissible}, we may assume in addition that if $\Sigma$ is nearly of finite type,
\begin{equation} \label{E:not ruled}
|\partial_t^{m_j}\psi(\xi+tN_j(\xi))| \sim c_\psi |\xi_{j'}|^{\ell_j},
\end{equation}
with $\ell_j,m_j$ as above.  

We observe that in the finite type case, \eqref{E:nearly finite type}, \eqref{E:admissible}, and \eqref{E:not ruled} hold with $\ell_j = 0$, and we will henceforth use these three inequalities in both the finite type and nearly finite type cases.  By a partition of unity and parabolic rescaling, we may assume that if $\Sigma$ is either finite type or nearly so and $N$ is some fixed integer, for all $\ell_j+m_j < |\alpha| \leq N$, 
\begin{equation} \label{E:polynomial like}
|\partial^\alpha \psi| \leq c_N c_\psi,
\end{equation}
with $c_N$ as small as we like.  (Here $c_\psi$ is the constant in \eqref{E:admissible} and \eqref{E:not ruled}.)  Henceforth we fix some large-enough $N$, and $N > (m_1+\ell_1+m_2+\ell_2+4)^2$ (say) will do.

We will say that an admissible surface $\Sigma$ with these additional properties is \textit{sufficiently localized}.  

\begin{definition} \label{D:type}
We definite a notion of essential type for $\psi$ (or, when convenient, $\Sigma$), both overall and in the directions $N_j$ (resp.\ $\tilde N_j$).  We say that $\psi$ is of type $m_j=1$ in the direction $N_j$ if $\Sigma$ is ruled in the direction $\tilde N_j$.  If $\Sigma$ is either nearly finite type or finite type in the direction $\tilde N_j$, then we say that $\psi$ is of essential type $n_j$ in the direction $N_j$, where $n_j=m_j+\ell_j$ is the minimal such sum for which \eqref{E:finite type} or \eqref{E:admissible} (resp.) holds.  The essential type $n$ of $\psi$ is defined to be the maximum of $n_1,n_2$, $n:=\max\{n_1,n_2\}$.  
\end{definition}

We remark here that it is a bit of an oddity that ruled surfaces have low (rather than infinite) type, and the reason is that the type may be thought of as the type of the surface along lines joining separated pairs of points (definition forthcoming).  In the ruled case, such lines must be transverse to the $N_j$, but in the finite and nearly finite type cases, the lines may be tangent to integral curves of the $N_j$ to some finite order.

We will prove a bilinear adjoint restriction estimate for sufficiently localized $\Sigma$, based on the type, and then we will use this estimate to prove linear extension estimates.  By the triangle inequality, this will imply linear extension estimates for all admissible $\Sigma$.  We note that the class in question includes (locally) all analytic hyperbolic surfaces.  

\begin{lemma} \label{L:analytic admissible}
Every precompact, hyperbolic, analytic 2-surface in $\R^3$ is admissible.
\end{lemma}
\begin{proof}
We may assume $\Sigma$ takes the form \eqref{E:Sigma}, with $\psi$ obeying \eqref{E:psi near psi0} and that $\Sigma$ is not ruled in the direction $\tilde N_j$.  From the definition of the $N_j$, $\Sigma$ contains a line through $(\psi(\xi),\xi)$ if and only if it contains a line in the direction of one of the $\tilde N_j$, or equivalently, if and only if 
$$
    \partial_t^2 \psi(\xi+tN_j(\xi)) = 0,
$$
for all $t$ such that $\xi,\xi+tN_j(\xi) \in 2Q$.  Thus if $\Sigma$ contains no lines in the direction $\tilde N_j$, by analyticity, for each $\xi$, $\Sigma$ makes at most finite order of contact with the line through $(\psi(\xi),\xi)$ in the direction $\tilde N_j(\xi)$, whence \eqref{E:finite type} holds.  

Alternatively, if $\Sigma$ contains a line but isn't ruled, any lines in the direction $\tilde N_j$ are isolated, so we may assume that there is only one.  Making a linear transformation of $\Sigma$, we may assume that this unique line passes through the origin, whence it must be the $\xi_j$ axis and, after suitable localization, we must be in the nearly finite type case.  
\end{proof}

\subsection{Separated points}

We assume henceforth that $\Sigma$ is sufficiently localized.  

We say that the pair of points $(\xi,\eta) \in Q^2$ is \textit{separated} (or that $\xi$ and $\eta$ are separated) if $\xi\neq \eta$ and the line segment containing $(\psi(\xi),\xi)$ and $(\psi(\eta),\eta)$ is not contained in $\Sigma$. We let $Z \subseteq Q^2$ denote the set of non-separated pairs.  For admissible $\psi$, it is easy to check that $Z$ is a closed set of measure zero.  

We now introduce a more quantitative notion of separation.  Before doing so, we observe that $\Sigma$ contains the line segment through $\zeta:=(\psi(\xi),\xi)$ and $\zeta':=(\psi(\eta),\eta)$ if and only if $\partial_t^2\psi(\xi+t(\eta-\xi))\equiv 0$.  We may expand the left hand side to express the latter condition as
\begin{equation} \label{E:IP condition linear}
\langle D^2\psi(\xi+t(\eta-\xi))(\eta-\xi),\eta-\xi \rangle \equiv 0,
\end{equation}
which in turn is equivalent to 
\begin{equation} \label{E:wedge condition linear}
N_j(\xi+t(\eta-\xi)) \wedge (\eta-\xi) \equiv 0.
\end{equation}

\begin{definition} \label{D:separated}
    Fix some large dimensional constant $C$, with $2^{-C} < c_0$.  Let $(\xi,\eta)$ be separated with $|\xi-\eta| \sim 2^{-j}$.  Let $k \geq j+C$.  We say that $(\xi,\eta)$ is:  
    \begin{itemize}
        \item $(2^{-j},2^{-k})$-separated if 
            \begin{equation} \label{E:jk separation}
                \max_{t \in [-1,2]} \bigl|N_1(\xi+t(\eta-\xi))\wedge (\xi-\eta)| \sim 2^{-k};
            \end{equation}
        \item $(\xi,\eta)$ is $(2^{-k},2^{-j})$ separated if 
            \begin{equation} \label{E:kj separation}
                \max_{t \in [-1,2]} \bigl|N_2(\xi+t(\eta-\xi))\wedge (\xi-\eta)| \sim 2^{-k}.
            \end{equation}
    \end{itemize}
    Finally, we say that $(\xi,\eta)$ is:
    \begin{itemize}
        \item $(2^{-j},2^{-k})$-separated with $|k-j|\leq C$ if $2^{j}(\xi-\eta)$ makes angle at least $2c_0$ with $N_1(\xi)$ and $N_2(\xi)$;
        \item $(2^{-j},0)$-separated if $|\xi-\eta| \sim 2^{-j}$ and the line joining $(\psi(\xi),\xi)$ and $(\psi(\eta),\eta)$ lies in $\Sigma$ and has direction $\tilde N_1(\xi)$;
        \item $(0,2^{-j})$-separated if $|\xi-\eta| \sim 2^{-j}$ and the line joining $(\psi(\xi),\xi)$ and $(\psi(\eta),\eta)$ lies in $\Sigma$ and has direction $\tilde N_2(\xi)$;
        \item $(0,0)$-separated if $\xi=\eta$.
    \end{itemize}
\end{definition}

This separation condition agrees up to constants with the separation condition stated for Theorem \ref{T:3}. For a heuristic argument as to why these are equivalent, not that (\ref{E:IP condition linear}) agrees with $|\mathrm{I\!I}_{\xi}(v_{\xi, \eta}, v_{\xi, \eta})|$ if $\langle\xi - \eta, \grad \psi(\xi) - \grad\psi(\eta) \rangle= 0$. One can check that if either (\ref{E:weak condition}) fails or if (\ref{E:IP condition linear}) fails, then the line segment joining $(\psi(\xi), \xi)$ and $(\psi(\eta), \eta)$ must stay close to the curve $\gamma(t)$ and along both we must have $\langle\xi - \eta, \grad \psi(\xi) - \grad\psi(\eta) \rangle \ll 1$, which implies that the other of (\ref{E:weak condition}) fails or if (\ref{E:IP condition linear}) must fail as well. This argument can be made formal using similar techniques as used in the rest of this section.

Recalling that $N_1$ and $N_2$ are nearly transverse to one another and that the conditions \eqref{E:IP condition linear} and \eqref{E:wedge condition linear} are equivalent to $\Sigma$ containing the line through $\xi$ with direction $\tilde N_j(\xi)$, we see that the preceding conditions are mutually exclusive for a given pair $\xi,\eta$, provided $C$ is chosen sufficiently large.  

\subsection{Rectangles}

With this quantitative definition of separation in place, we will define a a few different equivalent families of ``rectangles,'' which we will prove are rightly called rectangles under certain additional conditions on their lengths and widths. Before doing so, we note that by the Implicit Function Theorem, there exist smooth submersions $\pi_j:2Q \to \R$, $j=1,2$, such that $d\pi_j N_j \equiv 0$ on $2Q$.  By transversality of $N_1,N_2$, $\pi_j(e^{tN_{j'}}(\xi))$ is a diffeomorphism in $t$ for all $\xi \in 2Q$.

We are now ready to define our rectangles.  Fix $\xi \in Q$, $j \geq 0$, and $k \geq j$.  We define
\begin{align*}
R_1(\xi,j,k)&:=\{\eta \in 2Q : (\xi,\eta) \,\text{is $(2^{-j'},2^{-k'})$-separated,}\\
&\qquad\qquad\text{for some $\infty \geq j' \geq j$ and $\infty \geq k' \geq k$}\}\\
R_1(\xi,k,j)&:=\{\eta \in 2Q : (\xi,\eta) \, \text{is $(2^{-k'},2^{-j'})$-separated,} \\
&\qquad\qquad \text{for some $\infty \geq k' \geq k$ and $\infty \geq j' \geq j$}\}\\
R_2(\xi,j,k)&:=\{\eta \in 2Q : \eta = \xi+tN_1(\xi) + \zeta, \ctc{some} |t| \leq 2^{-j}, |\zeta| \leq 2^{-k}\}\\
R_2(\xi,k,j)&:=\{\eta \in 2Q : \eta = \xi+tN_2(\xi) + \zeta, \ctc{some} |t| \leq 2^{-j}, |\zeta| \leq 2^{-k}\}\\
R_3(\xi,j,k)&:=\{\eta \in 2Q : |\pi_2(\xi)-\pi_2(\eta)| \leq 2^{-j}, |\pi_1(\xi)-\pi_1(\eta)| \leq 2^{-k}\}\\
R_3(\xi,k,j)&:=\{\eta \in 2Q : |\pi_1(\xi)-\pi_1(\eta)| \leq 2^{-j}, |\pi_2(\xi)-\pi_2(\eta)| \leq 2^{-k}\}\\
R_4(\xi,j,k)&:=\{e^{t_2N_2}e^{t_1N_1}(\xi) : |t_2| \leq 2^{-k}, |t_1| \leq 2^{-j}\}\\
R_4(\xi,k,j)&:=\{e^{t_2N_1}e^{t_1N_2}(\xi) : |t_2| \leq 2^{-j}, |t_1| \leq 2^{-k}\}\\
R_5(\xi,j,k)&:=\{e^{t_2N_1}e^{t_1N_2}(\xi) : |t_1| \leq 2^{-k}, |t_2| \leq2^{-j}\}\\
R_5(\xi,k,j)&:=\{e^{t_2N_2}e^{t_1N_1}(\xi) : |t_1| \leq 2^{-j}, |t_2| \leq2^{-k}\}.
\end{align*}
We are using the convention $2^{-\infty}=0$ in the definitions of $R_1$, and we note the relative order of $j,k$ alternates between subsequent lines.)

We observe that $R_2(\xi,\cdot,\cdot)$ is always a (near-)rectangle (with half-disks glued on the ends).  The other rectangles are each useful in their own way, and it turns out that all of these rectangles are comparable under certain constraints on $j,k$ for which the following lemma will lay the groundwork.  

\begin{lemma}\label{L:minimum separation}
If $\xi \in 2Q$, then for all points $\xi \neq \eta \in 2Q$ and $i=1,2$,
\begin{equation} \label{E:minimum sep on line}
\max_{t \in [-1,2]}\bigl|N_i(\xi+t(\eta-\xi)) \wedge (\eta-\xi)\bigr| 
\gtrsim 
|\eta-\xi| \max_{t \in [-1,2]}\bigl|N_i\bigl(\xi+t|\eta-\xi|N_i(\xi)\bigr) \wedge N_i(\xi)\bigr|
\end{equation}
Furthermore, if $\Sigma$ does not contain the line through $(\psi(\xi),\xi)$ in the direction $\tilde N_i(\xi)$, then for all points $\eta \neq \xi$,
    \begin{equation}\label{E:minimum max}
        \max_{t \in [-1,2]}\bigl|N_i(\xi+t(\eta-\xi)) \wedge (\eta-\xi)\bigr|
        \gtrsim
        c_\psi |\xi-\eta|^{m_i-1}|\xi_{i'}|^{\ell_i},
    \end{equation}
    with $c_\psi, m_i,\ell_i$ as in \eqref{E:not ruled}.
\end{lemma}
\begin{proof}
 The roles of the indices being symmetric, we will prove this with $i=1$.  
 
 We take the claims in order.  If $|N_i(\xi) \wedge (\eta-\xi)| \gtrsim \rm{RHS}\eqref{E:minimum sep on line}$, then we are done, by taking $t=0$ on the left hand side of \eqref{E:minimum sep on line}.   Otherwise, $\eta = \xi+sN_i(\xi)+\zeta$, with $s = |\xi-\eta|$ and $|\zeta| \ll \rm{RHS}\eqref{E:minimum sep on line}$, and \eqref{E:minimum sep on line} follows from the triangle inequality.  
 
We turn now to \eqref{E:minimum max}, which, in view of \eqref{E:minimum sep on line} and $i=1$, it suffices to prove for $\eta := \xi+sN_1(\xi)$  with $|s| = |\xi-\eta|$.  By the definition of $N_1,N_2$ and $|N_i-e_i| \ll 1$ throughout $Q$,
\begin{equation} \label{E:wedge to IP}
|N_1(\zeta) \wedge N_1(\xi)| \sim |\langle D^2\psi(\zeta) N_1(\xi),N_1(\xi) \rangle|, \qquad \zeta \in Q.
\end{equation}
The right hand side will be easier to work with.  We have assumed that 
$$
|\partial_t^{m_1}\psi(\xi+tN_1(\xi))| \sim c_\psi |\xi_2|^{\ell_1},
$$
and either $\ell_1=0$ or $\xi_2 \neq 0$ (since we have assumed $\eta$ is separated from $\xi$).  Equivalently, 
\begin{equation} \label{E:m1-2 deriv}
    \bigl|\partial_t^{m_1-2} \langle D^2 \psi(\xi+tN_1(\xi))N_1(\xi),N_1(\xi)\rangle\bigr| \sim c_\psi |\xi_2|^{\ell_1}.
\end{equation}
By the assumption that $\Sigma$ is sufficiently localized, replacing the $m_1-2$ derivative on the left of \eqref{E:m1-2 deriv} with a higher-order $t$ derivative gives a quantity much smaller than the right hand side of \eqref{E:m1-2 deriv}, so we can use the $m_1-2$ order Taylor expansion to approximate:
\begin{equation} \label{E:polynomial approx}
    \begin{aligned}
        &\max_{t \in [-1,2]}|N_1(\xi+t(\eta-\xi)) \wedge (\eta-\xi)|
        \sim
        \max_{t \in [-1,2]} |\xi-\eta| |\langle D^2\psi(\xi+stN_1(\xi))N_1(\xi),N_1(\xi)\rangle|\\
        &\qquad \gtrsim
        \max_{t \in [-1,2]} |\xi-\eta| \bigl(|\sum_{m=0}^{m_1-2} \tfrac1{m!}t^m\tfrac{d^m}{dt^m}|_{t=0}\langle D^2\psi(\xi+stN_1(\xi))N_1(\xi),N_1(\xi)\rangle| \\
        &\qquad\qquad\qquad\qquad\qquad\qquad- cc_\psi|\xi-\eta|^{m_1-1}|\xi_2|^{\ell_1}\bigr)
        \gtrsim c_\psi |\xi-\eta|^{m_1-1}|\xi_2|^{\ell_1}.
    \end{aligned}
    \end{equation}
This verifies \eqref{E:minimum max} for $\eta = \xi+sN_1(\xi)$, and hence (as  noted above) in general.  
\end{proof}

In view of the preceding, for $\xi \in Q$, $i \in \{1,2\}$, and $j \in \N$, we define an integer $k^i(\xi,j)$ by the relation
\begin{equation} \label{E:k max}
2^{-k^i(\xi,j)} \sim \max_{t \in [-1,2]}|N_i(\xi+2^{-j}tN_i(\xi)) \wedge 2^{-j}N_i(\xi)|.
\end{equation}
We observe that 
$$
|N_i(\xi+2^{-j}tN_i(\xi)) \wedge 2^{-j}N_i(\xi)| \leq 2^{-j}|N_i(\xi+2^{-j}tN_i(\xi))  - N_i(\xi)| \leq  2^{-2j+1},
$$
so $k^i$ will always be quite large, \begin{equation}\label{E:ki lower bound}
    k^i(\xi,j) \geq 2j.
\end{equation}

\begin{lemma} \label{L:comparable rectangles}
For $1 \leq m,n \leq 5$, there exists a constant $C$ such that for every $\xi \in Q$, $j \geq C$, and $j \leq k$, 
\begin{equation} \label{E:equiv rectangles}
\begin{aligned}
R_m(\xi,j+C,k+C) \subseteq &R_n(\xi,j,k) \subseteq R_m(\xi,j-C,k-C), \qquad k \leq k^1(\xi,j),\\
R_m(\xi,k+C,j+C) \subseteq &R_n(\xi,k,j) \subseteq R_m(\xi,k-C,j-C), \qquad k \leq k^2(\xi,j).
\end{aligned}
\end{equation}
\end{lemma}

\begin{proof}
It suffices by symmetry to prove the first line of \eqref{E:equiv rectangles}.  More precisely, we will prove that for every $\xi \in Q$, $j \geq 0$, $k \leq j \leq k^1(\xi,j)$, $1 \leq m \leq 5$, and $\eta \in R_m(\xi,j,k)$, 
\begin{equation} \label{E:eta in Rm+1}
\eta \in R_{m+1}(\xi,j-C,k-C),
\end{equation}
where $R_{5+1}$ is defined to be $R_1$.

\textbf{Case $m=1$.}  Let $\eta \in R_1(\xi,j,k)$.  Then there exists some $j' \geq j, k' \geq k$ such that $\eta$ is $(2^{-j'},2^{-k'})$-separated from $\xi$.  Let $j'':=\min\{j',k'\}$, $k'':=\max\{j',k'\}$.  Then $|\xi-\eta| \sim 2^{-j''}$.  If $j'' = k'$, then in fact, $|\xi-\eta| \lesssim 2^{-k}$, so \eqref{E:eta in Rm+1} is immediate.  Alternatively, if $j''=j'$, \eqref{E:jk separation} implies that $\dist(\eta,\xi+\jp{N_1(\xi)}) \lesssim 2^{-k}$, so, once again, \eqref{E:eta in Rm+1} holds.  

\textbf{Case $m=2$.}  For $\eta \in R_2(\xi,j,k)$, $|\pi_2(\xi)-\pi_2(\eta)| \lesssim |\xi-\eta| \lesssim 2^{-j}$.  By the definition of $R_2(\xi,j,k)$, $\eta = \xi+tN_1(\xi)+\zeta$, some $|\zeta| \leq 2^{-k}$.  Thus by continuity of $N_1$ and the definition of $k^1(\xi,j)$,
\begin{align*}
|N_1(\eta)-N_1(\xi)| &\sim |N_1(\eta)\wedge N_1(\xi)| \leq \max_{|t| \sim 2^{-j}}|N_1(\xi+ tN_1(\xi)) \wedge N_1(\xi)| + C2^{-k} \\
&\qquad \lesssim 2^{j-k^i(\xi,j)} + 2^{-k} \lesssim 2^{j-k}.  
\end{align*}
From this, we see that the line through $\xi$ in the direction $N_1(\xi)$ stays within $2^{-k+C}$ of the integral curve of $N_1$ through $\xi$, so $\eta = e^{tN_1}(\xi)+\zeta'$, $|\zeta'| \lesssim 2^{-k}$.  As $\pi_1$ is constant along the integral curves of $N_1$, we see that $|\pi_1(\eta)-\pi_1(\xi)| \lesssim 2^{-k}$, so \eqref{E:eta in Rm+1} holds.  

\textbf{Case $m=3$.}  Let $\eta \in R_3(\xi,j,k)$.  Since $N_2$ is transversal to $N_1$, and $|N_2|\equiv|N_1|\equiv1$, the integral curve of $N_2$ through $\eta$ must intersect the integral curve of $N_1$ through $\xi$.  That is, we may write $\eta = e^{t_2 N_2}e^{t_1N_1}(\xi)$, some small $t_1,t_2$.  Moreover, since $N_2$ is transversal to the kernel of $d\pi_1$ (which is spanned by $N_1$), $|t_2| \sim |\pi_1(\eta)-\pi_1(e^{t_1N_1}(\xi)| = |\pi_1(\eta)-\pi_1(\xi)| \leq 2^{-j}$.  Therefore \eqref{E:eta in Rm+1} holds.  

\textbf{Case $m=4$.}  Let $\eta \in R_4(\xi,j,k)$, that is, $\eta = e^{t_2N_2}e^{t_1N_1}(\xi)$, $|t_2| \leq 2^{-k}$, $|t_1| \leq 2^{-k}$.  Arguing as in the preceding case, the integral curve of $N_1$ through $\eta$ intersects the integral curve of $N_2$ through $\xi$.  Thus $\eta = e^{s_2N_1}e^{s_1N_2}(\xi)$, with $s_1,s_2$ small.  Applying the Baker--Campbell--Hausdorff formula (e.g., \cite{CNSW99}) to this second expression, $\eta = e^{s_1N_2}e^{s_2N_1}(\xi) + O(|s_1s_2|)$.  On the other hand, direct computation shows that $t \mapsto e^{t_2N_2}e^{t_1N_1}(\xi)$ is a local diffeomorphism, whence the equation
$$
e^{s_1N_2}e^{s_2N_1}(\xi) + O(|s_1s_2|) = e^{t_2N_2}e^{t_1N_1}(\xi)
$$
implies that $(s_1,s_2)+O(|s_1s_2|) = (t_2,t_1)$.  Since $s_1,s_2$ are both small, $|s_1| \lesssim |t_2| \leq 2^{-k}$ and $|s_2| \lesssim |t_1| \leq 2^{-j}$.  Therefore \eqref{E:eta in Rm+1} holds.  

\textbf{Case $m=5$.}  Let $\eta \in R_5(\xi,k,j)$.  Writing $\eta = e^{t_2N_2}e^{t_1N_1}(\xi)$ with $|t_1| \leq 2^{-j}$, $|t_2| \leq 2{-k} \leq 2^{-j}$, $|\eta-\xi| \lesssim 2^{-j}$, and $|\eta-e^{t_1N_1}(\xi)| \lesssim 2^{-k}$.  Thus, it suffices to prove that
\begin{equation} \label{E:etN1xi-xi near tN1}
\max_{s \in [-1,2]}|N_1(\xi+s(e^{t_1N_1}(\xi)-\xi)) \wedge (e^{t_1N_1}(\xi)-\xi)| \lesssim 2^{-k}.  
\end{equation}
We have assumed that $k \leq k^1(\xi,j)$, so 
\begin{equation} \label{E:N1 on line near N1}
\max_{s \in [-1,2]}|N_1(\xi+2^{-j}sN_1(\xi)) \wedge 2^{-j}N_1(\xi)| \lesssim 2^{-k}.  
\end{equation}
Since the wedge of two unit vectors approximates the distance when the distance is small,
$$
\max_{s \in [-1,2]}|N_1(\xi+2^{-j}sN_1(\xi) + O(2^{-k})) \wedge N_1(\xi)| \lesssim 2^{j-k}.
$$
By Picard iteration, then, $|(e^{t_1N_1}(\xi)-\xi) - t_1N_1(\xi)| \lesssim 2^{-k}$, for $|t_1| \lesssim 1$.  By the triangle inequality, \eqref{E:etN1xi-xi near tN1} then follows directly from \eqref{E:N1 on line near N1}.  We finally see that \eqref{E:eta in Rm+1} holds in the final case.  
\end{proof}

The equivalence of the different families of rectangles established in Lemma~\ref{L:comparable rectangles} will allow us to establish several important properties.  Henceforth, we will use the first definition in the series as our default rectangles and denote $R(\xi,j,k):=R_1(\xi,j,k)$.  Using comparability with the $R_2$-rectangles, we see that the doubling property holds.  Namely, for $0 \leq j \leq k$, and $C$ a fixed constant, 
\begin{equation} \label{E:doubling}
\begin{gathered}
|R(\xi,j,k)| \sim |R(\xi,j+C,k+C)|, \ctc{if} k \leq k^1(\xi,j), \ctc{and}\\
 |R(\xi,k,j)| \sim |R(\xi,k+C,j+C)|, \ctc{if} k \leq k^2(\xi,j),
\end{gathered}
\end{equation}
where the implicit constants depend on $C$.  Because our separation conditions \eqref{E:jk separation} and \eqref{E:kj separation} are symmetric, our balls obey the symmetry property
\begin{equation} \label{E:symmetry}
\eta \in R(\xi,j,k) \implies \xi \in R(\eta,j,k).
\end{equation}

We have also have an engulfing property, whose proof is a little more involved.    

\begin{lemma} \label{L:engulfing1}
There exists a constant $C$ such that if $C \leq j \leq k \leq k^1(\xi,j)$
\begin{itemize}
\item[{\rm(a)}] If $(\xi,\eta)$ is $(2^{-j},2^{-k})$-separated, then any $(\xi',\eta') \in R(\xi,j+C,k+C) \times R(\eta,j+C,k+C)$ is $(2^{-j'},2^{-k'})$-separated, for some $|j-j'| \leq 1$, $|k-k'| \leq 1$; and
\item[{\rm(b)}] If $R(\xi,j,k) \cap R(\eta,j,k) \neq \emptyset$, then 
$$
R(\xi,j,k) \cup R(\eta,j,k) \subseteq R(\xi,j-C,k-C) \cap R(\eta,j-C,k-C).
$$
\end{itemize}
\end{lemma}

Lemma~\ref{L:engulfing1} has an immediate corollary by symmetry in the two indices.  

\begin{lemma} \label{L:engulfing2}
There exists a constant $C$ such that if $C \leq j \leq k \leq k^2(\xi,j)$, then items {\rm(a)} and {\rm(b)} from Lemma~\ref{L:engulfing1} hold with the places of $j,k$ interchanged each time they appear.  
\end{lemma}

It suffices to prove the first lemma.

\begin{proof}[Proof of Lemma~\ref{L:engulfing1}]
Let $C \leq j \leq k \leq k^1(\xi,j)$, with $C$ sufficiently large.  Both (a) and (b) involve variations on the triangle inequality.  We begin with (a).  Let $\xi,\eta,\xi',\eta'$ satisfy the hypotheses in (a). 

By assumption, $|\xi-\xi'| \lesssim 2^{-j-C}$, $|\eta-\eta'| \lesssim 2^{-j-C}$, $|\xi-\eta| \sim 2^{-j}$, whence by direct application of the Euclidean triangle inequality, $|\xi'-\eta'| \sim 2^{-j'}$ (the change from $j$ to $j'$ is so `$\sim$' is used with the same implicit constants as in the definition).  

Controlling the second separation term is a bit more involved.  We need to prove that
\begin{equation} \label{E:second dist xi' eta'}
\max_{t \in [-1,2]}|N_1(\xi'+t(\eta'-\xi')) \wedge (\xi'-\eta')| \sim 2^{-k'}.
\end{equation}
\begin{equation} \label{E:xi' to xi}
\xi'+t(\eta'-\xi') = \xi+s(\eta-\xi) + O(2^{-k-C}).
\end{equation}
To this end, we rearrange
$$
\xi'+t(\eta'-\xi') = \xi+t(\eta-\xi) + (\xi'-\xi) + t(\eta'-\eta-\xi'+\xi).
$$
By Lemma~\ref{L:comparable rectangles}, our assumption implies
$$
\xi' = \xi-u N_1(\xi) + O(2^{-k-C}), \qquad \eta' = \eta - vN_1(\eta) + O(2^{-k-C}),
$$
with $|u|,|v| \lesssim 2^{-j-C}$.  We can thus rewrite
$$
(\xi'-\xi) + t(\eta'-\eta-\xi'+\xi) = uN_1(\xi) + t(vN_1(\eta)-uN_1(\xi)) + O(2^{-k-C}). 
$$
As the roles of $\xi,\eta$ are symmetric in our separation condition, we have
$$
\eta = \xi + w N_1(\xi) + O(2^{-k}), \qquad \xi = \eta + sN_1(\eta) + O(2^{-k}),
$$
with $|w|,|s| \sim 2^{-j}$.  Thus 
$$
uN_1(\xi) + t(vN_1(\eta)-uN_1(\xi)) = (-\tfrac{u}{w} + \tfrac{tv}{s} - \tfrac{tu}{w})(\xi-\eta) + O(2^{-k-C}),
$$
and we observe that $|-\tfrac{u}{w} + \tfrac{tv}{s} - \tfrac{tu}{w}| \lesssim 2^{-(k-j)-C} \lesssim 2^{-C}$, whence the claim \eqref{E:xi' to xi} is proved.  Conversely to the claim, given $s \in [-1+C2^{-C},2-C2^{-C}]$, by the same argument, there exists $t_s \in [-1,2]$, with $|t_s-s| = O(2^{-C})$, such that \eqref{E:xi' to xi} holds.  Finally, arguing as above,
$$
\xi'-\eta' = (\xi-\eta)(1+O(2^{-C})) + O(2^{-k-C}).
$$

Since $N_1$ is smooth with bounded derivatives, we then obtain
\begin{align*}
&\max_{s \in [-1+C2^{-C},2-C2^{-C}]} |N_1(\xi+s(\eta-\xi)) \wedge (\xi-\eta)|(1+O(2^{-C})) + O(2^{-k-C})\\
& \qquad \leq \max_{t \in [-1,2]} |N_1(\xi'+t(\eta'-\xi')) \wedge (\xi'-\eta')| \\
&\qquad \leq \max_{s \in [-1-C2^{-C},2+C2^{-C}]} |N_1(\xi+s(\eta-\xi)) \wedge (\xi-\eta)|(1+O(2^{-C})) + O(2^{-k-C}).  
\end{align*}
Since we have assumed that 
$$
\max_{t \in [-1,2]} |N_1(\xi+t(\eta-\xi)) \wedge (\xi-\eta)| \sim 2^{-k},
$$
to prove \eqref{E:second dist xi' eta'}, it suffices to prove that
\begin{equation} \label{E:outer to inner}
\max_{s \in I_C^+} |N_1(\xi+s(\eta-\xi)) \wedge (\xi-\eta)| \leq \max_{t \in I_C^-} |N_1(\xi+t(\eta-\xi)) \wedge (\xi-\eta)| + O(2^{-k-C}),
\end{equation}
where $I_C^\pm:=[-1-\pm C2^{-C},2+\pm C2^{-C}]$.  

Fix $s \in I_C^+$, and let $t:=\tfrac{1-C2^{-C}}{1+C2^{-C}} s$.  Then $t \in I_C^-$, and $|t-s| = O(2^{-C})$.  Let us decompose
$$
\xi-\eta = (\xi-\eta)^\parallel + (\xi-\eta)^\perp,
$$ 
the parallel and perpendicular parts with respect to $N_1(\xi)$.  By assumption, $|(\xi-\eta)^\perp| \lesssim 2^{-k}$.  

Thus, in the ruled case, when the integral curves of $N_1$ are lines (along which $N_1$ is constant),
\begin{align*}
&N_1(\xi+s(\eta-\xi))\wedge (\xi-\eta) = N_1(\xi+s(\eta-\xi)^\parallel)\wedge(\xi-\eta) + O(2^{-k-C})\\
&\quad =N_1(\xi+t(\eta-\xi)^\parallel)\wedge(\xi-\eta) + O(2^{-k-C}) = N_1(\xi+t(\eta-\xi))\wedge(\xi-\eta) + O(2^{-k-C}),
\end{align*}
and we are done.  (We may argue similarly in the nearly finite type case when $\xi_2=0$.)

In the finite type and nearly finite type cases, we decompose further.  By smallness of $(\xi-\eta)^\perp$,
\begin{equation} \label{E:N1 wedge perp}
N_1(\xi+s(\eta-\xi)) \wedge (\xi-\eta)^\perp = N_1(\xi+t(\eta-\xi)) \wedge (\xi-\eta)^\perp + O(2^{-k-C}).
\end{equation}
Writing $(\xi-\eta)^\parallel = uN_1(\xi)$, with $|u| \sim 2^{-j}$, we have also that
\begin{equation} \label{E:N1 wedge parallel}
N_1(\xi+s(\eta-\xi))(\xi-\eta)^\parallel = N_1(\xi-suN_1(\xi))\wedge uN_1(\xi) + O(2^{-k-C}).
\end{equation}
We would like to replace $s$ with $t$ on the right hand side, and are willing to pay the price of an $(O2^{-k-C})$ error term.  
We rewrite the main term from the right hand side of \eqref{E:N1 wedge parallel}:
$$
N_1(\xi-suN_1(\xi))\wedge uN_1(\xi) = (\tfrac12 + O(2^{-C}))uN_1(\xi)^T D^2\psi(\xi-suN_1(\xi))N_1(\xi),
$$
and we want to show that 
$$
|uN_1(\xi)^T D^2\psi(\xi-suN_1(\xi))N_1(\xi) - uN_1(\xi)^T D^2\psi(\xi-tuN_1(\xi))N_1(\xi)| \lesssim 2^{-k-C}.
$$
By the mean value theorem, the left hand side above is bounded by
\begin{equation} \label{E:s-t u^2 error}
|s-t| u^2 \sup_{w \in J_C}|\tfrac{d^3}{dw^3} \psi(\xi-wN_1(\xi))|.
\end{equation}
Here $J_C$ is the set of all $s'u$ with $s'$ lying between $s$ and $t$.  

Looking back to \eqref{E:not ruled}, we recall that $\ell_j$ was the minimal index (possibly zero) with $\partial_{\xi_{j'}}^{\ell_j} \partial_t^2 \psi(\xi+tN_j(\xi))\not\equiv 0$.  Recalling also the conditions for $\Sigma$ to be sufficiently localized, we consider $P$, the $(m_1+3)$-th order Taylor polynomial of $\psi$, centered at $\xi$. We have
$$
\sup_{w \in J_C}|\tfrac{d^3}{dw^3} \psi(\xi-wN_1(\xi))|
\leq \sup_{w \in J_C}|\tfrac{d^3}{dw^3} P(\xi-wN_1(\xi))| + O(c_\psi 2^{-(m_1+1)j}|\xi_{j'}|^{\ell_j}).
$$
We denote $I_0:=[-1,2]$.  We will invoke the equivalence of norms on finite dimensional vector spaces (such as single-degree polynomials of degree at most $m_1+3$).  It will be helpful to observe in advance that $\tfrac12 I_0 \subset I_C^-$, $I_C^+ \subset 2I_0$.  We have:  
$$
\sup_{w \in J_C}|\tfrac{d^3}{dw^3} P(\xi-wN_1(\xi))| 
\lesssim 2^j \sup_{t' \in I_0}|\tfrac{d^2}{dw^2}|_{w=t'u} P(\xi-wN_1(\xi))|.  
$$
Moreover, since
$$
|\tfrac{d^{m_1}}{dw^{m_1}} P(\xi-w-N_1(\xi))| \sim c_\psi|\xi_{j'}|^{\ell_j},
$$
by hypothesis, 
$$
\sup_{t' \in I_0} |\tfrac{d^2}{dw^2}|_{w=t'u} P(\xi-wN_1(\xi))| \gtrsim c_\psi |\xi_{j'}|^{\ell_j}2^{-(m_1-2)j}.
$$
We may thus absorb Taylor remainders, and, by the triangle inequality, write
\begin{align*}
&2^j \sup_{t' \in I_0}|\tfrac{d^2}{dw^2}|_{w=t'u} P(\xi-wN_1(\xi))| + O(c_\psi 2^{-(m_1+1)j}|\xi_{j'}|^{\ell_j})
\\
&\qquad \lesssim 2^j \sup_{t' \in I_0} |\tfrac{d^2}{dw^2}|_{w=t'u} \psi(\xi-wN_1(\xi))|
\sim 2^j \sup_{t' \in I_0}|N_1(\xi-t'uN_1(\xi)) \wedge N_1(\xi)|.
\end{align*}
So finally,
$$
|s-t|u^2 \sup_{w \in J_C}|\tfrac{d^3}{dw^3} \psi(\xi-wN_1(\xi))|
\lesssim 2^{-C} \sup_{t' \in I_0} |N_1(\xi-t'u N_1(\xi)) \wedge (t'u N_1(\xi))|
\lesssim 2^{-k-C},
$$
which is what we wanted to prove.  We note that in the last inequality above, we have used Lemma~\ref{L:minimum separation} and our hypothesis that $(\xi,\eta)$ is $(2^{-j},2^{-k})$-separated.  This completes the proof of Part (a) of Lemma~\ref{L:engulfing1}.  

Part (b) is much simpler, and we omit its proof.  
\end{proof}

\subsection{Rescaled rectangles}

In this section, we prove a few rescaling lemmas that will simplify several later arguments.  

\begin{lemma} \label{L:rescaling 1 rectangle}
Let $\psi:2Q \to \R$ be an admissible, sufficiently localized function.  Let $\xi_0 \in Q$ and $C \leq j \leq k \leq k^1(\xi,j)$.  Define a matrix $A:=(2^{-j}N_1(\xi_0),2^{-k}N_2(\xi_0))$, and define 
$$
\tilde\psi(\eta):=2^{j+k}\bigl(\psi(A\eta+\xi_0)-\psi(\xi_0) - \nabla \psi(\xi_0) \cdot(A\eta) \bigr).
$$
Then $\tilde \psi$ is admissible.  Furthermore, if $\xi_1,\xi_2 \in AQ+\xi_0$ are $(2^{-j},2^{-k})$-separated with respect to $\psi$ and $\eta_i:=A^{-1}(\xi_i-\xi_0)$, $i=1,2$, then $\eta_1,\eta_2$ are $(1,1)$-separated with respect to $\tilde\psi$.  
\end{lemma}

We note that by Lemma~\ref{L:comparable rectangles}, $R(\xi_0,j+C,k+C) \subseteq AQ+\xi_0 \subseteq R(\xi_0,j-C,k-C)$.

\begin{proof}
Before we begin, it will be helpful to make a small reduction.  Namely, we may assume that for $i=1,2$, either $\xi_{0,i'} = 0$, or $\Sigma$ is nearly finite type in the $\tilde N_i$ direction and $|\xi_{0,i'}| \gg 2^{-\ell_{i'}}$.  Here $i':=2$ and $\ell_{i'}:=k$ for $i=1$, $i':=1$ and $\ell_{i'}:=j$ for $i=2$, and we recall that $\tilde N_i$ denotes the pushforward of $N_i$ to the graph $\Sigma$ of $\psi$.  Our reduction is easiest (but least helpful) when $\Sigma$ is either ruled or finite type in the $\tilde N_i$ direction.  Indeed, in these cases, there are no distinguished values of $\xi_{i'}$, and so the reduction may be obtained via a simple affine map that depends on $\xi_0,\psi$, but not on $j,k$.  If the surface is nearly finite type in the $\tilde N_i$ direction, then the line $\{\xi_{i'} = 0\}$ is special, being the only line in the surface in the $\tilde N_i$ direction.  However, if $|\xi_{0,i'}| \lesssim 1$, then by Lemma~\ref{L:engulfing1}, $R(\xi_0,j,k)$ is contained in a slightly larger rectangle $R(\xi_0',j',k')$, with $\eta_{0,i'} = 0$ and $2^{-j'} \lesssim 2^{-j}$, $2^{-k'}\lesssim 2^{-k}$.

Now we are ready to begin.  First we verify the $C^0$ condition \eqref{E:psi near psi0} for $\tilde \psi$.  We immediately see that $\tilde\psi(0) =0$, $\nabla \tilde\psi(0) = 0$, $\tfrac{\partial^2}{\partial \eta_i^2}\tilde\psi(0) = 0$, $i=1,2$.  Furthermore, for $\eta \in 2Q$, $A\eta+\xi_0 \in 2Q$, and 
$$
\tfrac{\partial^2}{\partial \eta_1 \partial \eta_2 }\tilde\psi(\eta) = N_1^T(\xi_0) D^2\psi(A\eta+\xi_0) N_2(\xi_0).
$$
By assumption, $\|N_i-e_i\|_{C^0(2Q)} < c_0$, $i=1,2$, and $\|D^2 \psi - D^2\psi_0\|_{C^0(2Q)} < c_0$, so we obtain
$$
\|\tfrac{\partial^2}{\partial \eta_1 \partial \eta_2 }\tilde\psi(\eta)-1\|_{C^0(2Q)} \ll 1,
$$
and the implicit constant may be taken to be independent of $j,k,\xi_0,\psi$.  We now need to control the growth of the diagonal entries of $D^2\tilde\psi$.  We have
\begin{align*}
|\tfrac{\partial^2}{\partial \eta_2^2} \tilde\psi(\eta)| &= 2^{j-k}|N_2(\xi_0)^T D^2 \psi(A\eta+\xi_0)N_2(\xi_0)|
\lesssim 2^{j-k}|N_2(A\eta+\xi_0) - N_2(\xi_0)| \\
&\qquad \lesssim 2^{j-k}2^{-j} \leq 2^{-k} \leq c_0,
\end{align*}
since $|A\eta| \lesssim 2^{-j}$.  For the other diagonal term, we start similarly and see
$$
|\tfrac{\partial^2}{\partial \eta_1^2} \tilde\psi(\eta)| \lesssim 2^{k-j}|N_1(A\eta+\xi_0)\wedge N_1(\xi_0)|.
$$
Since $|2^{-k}N_2(\xi_0)\eta_2| \lesssim 2^{-k}$, 
$$
|N_1(A\eta+\xi_0)\wedge N_1(\xi_0)| \leq |N_1(2^{-j}\eta_1N_1(\xi_0)+\xi_0) \wedge N_1(\xi_0)| + O(2^{-k}) \lesssim 2^{-k_1(\xi_0,j)}+2^{-k} \lesssim 2^{-k},
$$
resulting in 
$$
|\tfrac{\partial^2}{\partial \eta_1^2} \tilde\psi(\eta)| \lesssim 2^{k-j}2^{-k} = 2^{-j} \leq c_0.
$$
Putting these together, $\tilde\psi$ obeys \eqref{E:psi near psi0}.  

Next, we bound derivatives of $\tilde\psi$ of order at least 3.  By construction (and $k \geq j$), derivatives $\frac{\partial^{m+n}}{\partial \eta_1^m\partial \eta_2^n} \tilde \psi$, with $m+n \geq 3$ and $n \geq 1$ obey the same bounds as the corresponding derivatives of $\psi$ (indeed, these derivatives also obey stronger bounds), and so it remains to bound derivatives $\frac{\partial^m}{\partial \eta_1^m}$ with $m \geq 3$.  Our proof of these bounds naturally breaks into consideration of four cases, which depend on the type of $\Sigma$ in the direction $\tilde N_1$.  

If $\Sigma$ is ruled in the direction $\tilde N_1$, then the graph of $\tilde \psi$ is also ruled in the corresponding direction, and  $\frac{\partial^m}{\partial \xi_1^m}\tilde\psi \equiv 0$ for $m \geq 2$.  

For the remaining cases, it will be helpful to note that, by the Mean Value Theorem (and definition of $\tilde\psi$), 
\begin{equation} \label{E:simplify mth deriv}
\frac{\partial^m}{\partial \eta_1^m} \tilde \psi(\eta) = \frac{\partial^m}{\partial \eta_1^m} \tilde \psi(\eta_1,0) + O(2^{-(m-1)j}\|\psi\|_{C^{m+1}}), \qquad \eta \in 2Q.
\end{equation}
Since second term on the right hand side is always much smaller than 1, we need only bound the first term on the right.

When $\Sigma$ is nearly finite type in the direction $\tilde N_1$ and $\xi_{0,2} = 0$, 
$$
\frac{\partial^m}{\partial \eta_1^m} \tilde \psi(\eta_1,0) \equiv 0,
$$
so the first term on the right of \eqref{E:simplify mth deriv} is identically zero.  Hence, by the observation just below \eqref{E:simplify mth deriv}, derivatives in $\eta_1$ of order at least 3 are suitably bounded.  
Furthermore, the graph $\tilde\Sigma$ of $\tilde\psi$ contains exactly one line in the direction corresponding to $\tilde N_1$, namely, the $x$-axis.  We see that \eqref{E:nearly finite type} holds with the same value of $\ell_1$, and the minimality condition on $\ell_1$  and Mean Value Theorem give \eqref{E:admissible} and \eqref{E:not ruled} with $c_{\tilde\psi} = 2^{j+k-k\ell_1}c_\psi$.  Inequality \eqref{E:polynomial like} likewise follows (indeed, high-order derivatives become even smaller relative to low-order ones after rescaling).  

We now turn to the case that either $\Sigma$ is finite type, or it is nearly so and $\xi_{0,2} \neq 0$.  We note that the main term from \eqref{E:simplify mth deriv} is
$$
|2^{k-(m-1)j} \partial_t^m|_{t=2^{-j}\eta_1} \psi(\xi_0+tN_1(\xi_0)) |,
$$
and we recall that $k \leq k_1(\xi_0,j)$.  By he definition of $k_1(\xi_0,j)$,
\begin{align*}
2^{-k_1(\xi_0,j)} &\sim \max_{t \in [-1,2]} |N_1(\xi_0+2^{-j}tN_1(\xi_0)) \wedge 2^{-j}N_1(\xi_0)|\\
&\sim \max_{t \in [-1,2]}|2^j \tfrac{d^2}{dt^2} \psi(\xi_0+2^{-j}t N_1(\xi_0))|.
\end{align*}
By the argument in \eqref{E:polynomial approx}, 
\begin{align*}
\max_{t \in [-1,2]}|2^j \tfrac{d^2}{dt^2} \psi(\xi_0+2^{-j}N_1(\xi_0))| \sim \sum_{m=2}^{m_1} 2^{-(m-1)j}\max_{t \in [-2^{-j},2^{-j+1}]}|\partial_t^m \psi(\xi_0+tN_1(\xi_0))|,
\end{align*}
allowing us to bound those derivatives with $m \leq m_1$.  Larger order derivatives are smaller, thanks to \eqref{E:polynomial like}, minimality of $\ell_1$, and \eqref{E:minimum max}.  In either case, the graph $\tilde\Sigma$ of $\tilde \psi$ is of finite type, with $c_{\tilde \psi} = 2^{j+k-k\ell_1}c_\psi$.  

The type conditions in the complementary directions may be obtained in a similar way in each of the following cases:  $\Sigma$ ruled in the direction $\tilde N_2$, $\Sigma$ nearly finite type in the direction $\tilde N_2$ and $\xi_{0,1}=0$, and $\Sigma$ either nearly finite type with $\xi_{0,1} \neq 0$ or finite type.  

Finally, let $\xi_1,\xi_2 \in AQ+\xi_0$ be a pair of $(2^{-j},2^{-k})$-separated points, and let $\eta_i:=A^{-1}(\xi_i-\xi_0)$.  We wish to show that $\eta_1,\eta_2$ are $(1,1)$-separated.  The initial separation condition is elementary, since $|\xi_1-\xi_2| \leq \|A\||\eta_1-\eta_2|$, and $\|A\| \sim 2^{-j}$.  For the second separation condition, we can use the definition of the vector fields $N_i$ and the chain rule to see that $N_1^{\tilde\psi}(\eta)$ equals $2^jA^{-1}N_1^\psi(A\eta+\xi_0)$, times a constant comparable to 1.  (Here the superscript indicates the function with respect to which the vector field is computed.)   Therefore,
\begin{align*}
&\max_{t \in [-1,2]} |N_1^{\tilde\psi}(\eta_1+t(\eta_2-\eta_1)) \wedge (\eta_1-\eta_2)| \\
&\qquad \sim 2^{-j}\max_{t \in [-1,2]} |A^{-1} N_1(\xi_1+t(\xi_2-\xi_1)) \wedge A^{-1}(\xi_1-\xi_2)|\\
&\qquad = 2^{-j} \det A^{-1}\max_{t \in [-1,2]} |N_1(\xi_1+t(\xi_2-\xi_1)) \wedge (\xi_1-\xi_2)|
\sim 2^{-j}2^{j+k}2^{-k} \sim 1,
\end{align*}
and the proof that $\eta_1$ and $\eta_2$ are $(1,1)$-separated is complete.    
\end{proof}

Combined with Lemma~\ref{L:engulfing1}, Part (a), Lemma~\ref{L:rescaling 1 rectangle} implies that two $(2^{-j},2^{-k})$-separated $2^{-j-C} \times 2^{-k-C}$ rectangles may be rescaled (via a single transformation) to a $(1,1)$-separated pair of $2^{-C} \times 2^{-C}$ rectangles, which is the statement of the following lemma.

\begin{lemma} \label{L:rescale rectangles to 1x1 pair}
Let $n \geq 0$, $\psi:2Q \to \R$ an admissible, sufficiently localized smooth function of essential type $n$, $\xi_0 \in Q$, $C \leq j \leq k \leq k^1(\xi_0,j)$, and $\eta_0$ be $(2^{-j},2^{-k})$-separated from $\xi_0$.  Set $R_1:=R(\xi_0,2^{-j-C},2^{-k-C})$ and $R_2:=R(\eta_0,2^{-j-C},2^{-k-C})$.  Then after rescaling by the map $A$ defined in Lemma~\ref{L:rescaling 1 rectangle}, the rectangles $R_1,R_2$ are $(1,1)$-separated $(2^{-C},2^{-C})$ rectangles for the admissible function $\tilde\psi$.  

\end{lemma}

As noted above, this lemma is a direct application of results already established.

\section{Bilinear extension estimates for admissible $\psi$}

Our main goal in this section is to prove the following.  

\begin{theorem}\label{T:bilinear rest}
Given $n \geq 0$, if $\psi:2Q \to \R$ is an admissible, sufficiently localized smooth function of essential type $n$, then for all $(1,1)$-separated $(2^{-C},2^{-C})$-rectangles $R_1,R_2 \subseteq Q$ and smooth $f_1,f_2$ supported on $R_1,R_2$, resp.,
\begin{equation} \label{E:bilinear rest}
\|\scriptE f_1 \scriptE f_2\|_2 \lesssim \|f_1\|_s\|f_2\|_s, \qquad s=s_n:= 2-\tfrac1{8n}.
\end{equation}
The implicit constant may be taken to depend only on $n$.  
\end{theorem}

We we have made no effort to optimize the value of $s_n$ (indeed, through some minor modifications, our proof gives a smaller--though less aesthetically pleasing--value), but it seems an interesting question to identify the smallest possible exponent.  We note that the case $s=2$ is well-known; it just follows from the standard $L^2$ result for $\scriptE$ and Cauchy--Schwarz.  

\begin{corollary}\label{C:rescaled bilinear rest}
For any $s > \frac{11}{7}$, the following holds. Let $m \geq 0$, $\psi:2Q \to \R$ an admissible, sufficiently localized smooth function of essential type $n$, $\xi_0 \in Q$, $C \leq j \leq k \leq k^1(\xi_0,j)$, and $\eta_0$ be $(2^{-j},2^{-k})$-separated from $\xi_0$.  Set $R_1:=R(\xi_0,2^{-j-C},2^{-k-C})$ and $R_2:=R(\eta_0,2^{-j-C},2^{-k-C})$. There exists $\varepsilon_{n,s}$ depending only on $n$ and $s$ such that for smooth $f_1,f_2$ supported on $R_1,R_2$, respectively and any $r > s'-\varepsilon_{n,s}$
\begin{equation} \label{E:bilinear rest arb scale}
\|\scriptE f_1 \scriptE f_2\|_s \lesssim 2^{(j+k)\left(\frac2r + \frac2s - 2\right)}\|f_1\|_r\|f_2\|_r,.
\end{equation}
The implicit constant may be taken to depend only on $n,r$, and $s$. 

Furthermore, if (\ref{E:T2 assumption}) holds for up to some exponent $p_0$, then this corollary holds for any $s> \frac{p_0}{2}$.
\end{corollary}

The first claim of the corollary follows by Lemma~\ref{L:rescale rectangles to 1x1 pair}, the change of variables formula, and interpolating with \cite[Theorem 1.1]{GLX26}. The ``furthermore'' uses the same reasoning with the assumption (\ref{E:T2 assumption}) replacing \cite[Theorem 1.1]{GLX26}. As previously noted, the separation condition in Corollary \ref{C:rescaled bilinear rest} is essentially equivalent to the informal separation condition in Theorem \ref{T:3}, so Theorem \ref{T:3} follows from the $j = k = 0$ case of Corollary \ref{C:rescaled bilinear rest}.

\begin{proof}[Proof of Theorem~\ref{T:bilinear rest}]
Our proof takes as its initial inspiration the now-classical bilinear restriction inequality of Moyua--Vargas--Vega, \cite{MVV99}.  By Plancherel, \eqref{E:bilinear rest} will follow if we can show that
$$
\|\widehat{\scriptE f_1\,\overline{\scriptE f_2}}\|_2 \lesssim \|f_1\|_s\|f_2\|_s.  
$$
By duality, it suffices to prove that
\begin{equation} \label{E:grt ineq}
\scriptT(f_1,f_2,f_3) \lesssim \|f_1\|_s\|f_2\|_s\|f_3\|_2, \qquad f_1,f_2,f_3 \geq 0,
\end{equation}
where
$$
\scriptT(f_1,f_2,f_3):= \iint_{R_1 \times R_2} \prod_{j=1}^3 f_j \circ \pi_j(\xi,\eta)\, d\xi\, d\eta,
$$
with $\pi_1,\pi_2:\R^4 \to \R^2$, $\pi_3:\R^4 \to \R^3$, and
$$
\pi_1(\xi,\eta):=\xi, \qquad \pi_2(\xi,\eta):=\eta, \qquad \pi_3(\xi,\eta):=(\psi(\xi)-\psi(\eta),\xi-\eta).
$$ 

Our goal is now to bound a multilinear generalized Radon transform, and our next step will be to reduce our problem to proving certain geometric inequalities, relating volumes with average fiber measure.  This reduction has now become routine; see \cite{TW03}, and the references therein.  By multilinear real interpolation, it suffices to prove  
\begin{equation} \label{E:rwt grt}
\scriptT(\chi_{F_1},\chi_{F_2},\chi_{F_3}) \lesssim |F_1|^{1/p_1}|F_2|^{1/p_2}|F_3|^{1/p_3},
\end{equation}
for open sets $F_1 \subseteq R_1$, $F_2 \subseteq R_2$, $F_3 \subseteq \pi_3(R_1 \times R_2)$ and all $(p_1,p_2,p_3)$ in some neighborhood of $(s,s,2)$.  Setting $\Omega:=\pi_1^{-1}(F_1) \cap \pi_2^{-1}(F_2) \cap \pi_3^{-1}(F_3)$, \eqref{E:rwt grt} will follow from
$$
|\Omega| \lesssim |\pi_1(\Omega)|^{1/p_1}|\pi_2(\Omega)|^{1/p_2}|\pi_3(\Omega)|^{1/p_3},
$$
which is equivalent, after a bit of arithmetic, to 
$$
|\Omega| \gtrsim \alpha_1^{b_1}\alpha_2^{b_2}\alpha_3^{b_3},
$$
with $\alpha_j:=|\Omega|/|\pi_j(\Omega)|$ and $b_j:=p_j^{-1}/(p_1^{-1}+p_2^{-1}+p_3^{-1}-1)$, $j=1,2,3$.  Let $\eps=\eps_n>0$ be a sufficiently small constant, and cover $\Omega$ by balls of radius $(\alpha_1\alpha_2)^\eps$.  One of these balls then contains an open subset $\tilde\Omega \subseteq \Omega$ with $|\tilde\Omega| \gtrsim (\alpha_1\alpha_2)^{4\eps}|\Omega|$.  Letting $\tilde\alpha_j$ denote the parameters associated to $\tilde\Omega$, $j=1,2,3$, we then have that $\tilde\alpha_j \gtrsim (\alpha_1\alpha_2)^{4\eps}\alpha_j$, whence $\tilde\alpha_1\tilde\alpha_2 \gtrsim (\alpha_1\alpha_2)^{1+8\eps}$.  Thus, $\diam(\tilde\Omega) \lesssim (\tilde\alpha_1\tilde\alpha_2)^{\delta}$, where $\delta:=\tfrac{\eps}{1+8\eps}$.  

Removing the tildes, our main theorem is thus reduced to the following lemma.

\begin{lemma} \label{L:isoperim reduction}
There exists a neighborhood of $(2/(4-s),2/(4-s),s/(4-s))=:b(n)$ such that for each $\rho > 0$, the inequality
\begin{equation} \label{E:rwt rev}
|\Omega| \gtrsim_\rho \alpha_1^{b_1}\alpha_2^{b_2}\alpha_3^{b_3}
\end{equation}
holds for open subsets $\Omega \subseteq R_1 \times R_2$ of diameter at most $C(\alpha_1\alpha_2)^\rho$ and all $(b_1,b_2,b_3)$ in the neighborhood, where the $\alpha_j$ are defined as above.  
\end{lemma}
Here we recall that $s=s_n=2-\tfrac1{2(n+1)}$, and we emphasize that, while the implicit constant is allowed to depend on $\rho$, the neighborhood is not.  We note that in the case $s=2$ (long-known, as noted just after the theorem statement), $b_1=b_2=b_3=1$.  

Our task is now to prove Lemma~\ref{L:isoperim reduction}.  First, we prove a small lemma whose use motivates the diameter reduction just made.   

\begin{lemma}\label{L:poly approx}[Christ, \cite{C08}]
Let $W \subseteq B \subseteq Q \subseteq \R^n$, with $W$ an open set, $B$ a ball, and $Q$ the unit cube.  Let $\Phi:2Q \to \R^n$ be a $C^2$ function, and let $\Psi:\R^n \to \R^n$ be a polynomial of degree at most $N$.  Assume that 
\begin{gather}
\label{E:Phi C2}
\|\Phi\|_{C^2(Q)} \leq 4,
\\
\label{E:Jacobian lb}
\inf_{W}|\det D\Phi| =: J_\Phi > 0,
\\
\label{E:Phi-Psi C2}
\|\Phi-\Psi\|_{C^2(B)} \leq c_n J_\Phi^2.
\end{gather}
Then $|\Phi(W)| \gtrsim_N \int_W|\det D\Phi|$.  
\end{lemma}

\begin{proof}
For the lemma, we credit a 20-year-old preprint available on the webpage of Michael Christ \cite{C08}.  The statement in \cite{C08} (Lemma 7.1) is different from what we need, and there is no arXiv/published version of \cite{C08}, so we give a complete, self-contained proof.  This argument draws heavily from the ideas of Christ but differs in some respects.

We note that, for $c_n$ sufficiently small, our assumptions imply that 
\begin{equation} \label{E:det DPhi DPsi}
|\det D\Phi(t)| \sim |\det D\Psi(t)| \gtrsim J_\Phi, \quad t \in B.
\end{equation}  Let $t_0 \in W$, and consider the axis-parallel cube $\rho_{t_0}$ centered at $t_0$ of diameter $c_n J_\Phi$, with $c_n$ sufficiently small.  By our hypotheses and \eqref{E:det DPhi DPsi}, both $\Phi$ and $\Psi$ are one-to-one on $\rho_{t_0}$.  (Indeed, this is a straightforward adaptation of the standard contraction mapping proof of the Inverse Function Theorem, such as that in \cite{R76}.)  Let $\{\rho_j\}$ be a disjoint collection of such cubes covering $B$.  Then, by the change of variables formula, for each $j$, 
\begin{equation} \label{E:injective on rhoj}
|\Phi(W \cap \rho_j)| \geq \int_{W \cap \rho_j} |\det D\Phi|.
\end{equation}

Thus, our proof will be complete if we can prove that a point $x \in \Phi(W)$ can lie in $\Phi(W \cap \rho_j)$ for at most $C_{N,n}$ values of $j$.  We note that the analogous statement for $\Psi$ follows from Bezout's theorem.  Indeed, by \eqref{E:det DPhi DPsi} and Bezout's theorem, something stronger holds.  Namely, a point $x$ can lie in $\Phi(\rho_j^*)$ for a bounded number of $j$ with $\rho_j \cap W \neq \emptyset$.  Here $\rho_j^* = C_n\rho_j$.  Thus, by \eqref{E:Phi-Psi C2}, it suffices to prove that
$$x
\Psi(\rho_j) + b \subseteq \Psi(\rho_j^*),
$$
when $\rho_j \cap W \neq \emptyset$.  Here $b$ is the ball of radius $c_n J_\Phi^2$, centered at 0.  As argued above, if the constant $c_n$ used in defining the $\rho_t$ is sufficiently small, $\Phi$ is one-to-one on $\rho_j^*$.  Thus it suffices to prove that 
$$
\Psi^{-1}(\Psi(\rho_j) + b) \subseteq \rho_j^*, 
$$ 
which would follow from 
\begin{equation} \label{E:b to b'}
\Psi^{-1}(\Psi(\rho_j) + b) \subseteq \Psi^{-1}(\Psi(\rho_j)) + b',
\end{equation}
with $b'$ the cube of diameter $\tfrac12 (C_n-1)c_n J_\Phi$.  Finally, Cramer's rule, our hypotheses, and \eqref{E:det DPhi DPsi} imply that $\|D\Psi^{-1}\| \leq C_n J_\Phi^{-1}$, and hence, by the mean value theorem, $\Psi^{-1}$ is $C_n J_\Phi^{-1}$-Lipschitz.  Recalling the definitions of $b,b'$, \eqref{E:b to b'} follows.  
\end{proof}

Now we turn to the heart of our proof.  By way of motivation, we observe that, in the limit as $n \to \infty$, $s \to 2$ and $b(n) \to (1,1,1)$.  The case $b_1=b_2=b_3=1$ corresponds to \eqref{E:grt ineq} in the case $p_1=p_2=p_3=2$, a well-known consequence of the $L^2$ theory \cite{T75, S77} and Cauchy--Schwarz.  The next step in the proof of our theorem will be to reprove this classical estimate through less-efficient means.  We will then optimize our proof to obtain \eqref{E:rwt rev} for $b$ in a neighborhood of $b(n)$.   

We begin this stage of the argument with a familiar setup, laid out in \cite{C98, TW03}, and many subsequent works.  

\begin{lemma}\label{L:refinements setup}
There exists a sequence of open sets
$$
\Omega:=\Omega^{(0)} \supseteq \Omega^{(1)}\supseteq \Omega^{(2)} \supset \cdots \supset \Omega^{(j)} \supset \cdots,
$$
with $|\Omega^{(j)}| \gtrsim_j |\Omega|$, such that for each $i=1,2,3$ and $(\xi,\eta) \in \Omega^{(j)}$, 
\begin{gather*}
|\{u \in \R^2 : |u|<C(\alpha_1\alpha_2)^\rho, (\xi+u,\eta) \in \Omega^{(j-1)}\}| \gtrsim_j \alpha_2,\\
 |\{v \in \R^2 : |v| < C(\alpha_1\alpha_2)^\rho, (\xi,\eta+v) \in \Omega^{(j-1)}\}| \gtrsim_j \alpha_1\\
|\{t \in \R: |t| \leq C(\alpha_1\alpha_2)^\rho, e^{tZ_3}(\xi,\eta) \in \Omega^{(j-1)}\}| \gtrsim_j \alpha_3.
\end{gather*}
Here $Z_3$ is the vector field given by
$$
Z_3(\xi,\eta):=(\nabla^\perp\psi(\xi)-\nabla^\perp\psi(\eta))\cdot(\nabla_\xi+\nabla_\eta),
$$
and $\nabla^\perp\psi:=(-\partial_2\psi,\partial_1\psi)$.  
\end{lemma}

\begin{proof}[Proof of Lemma~\ref{L:refinements setup}]
The proof of statements involving $\alpha_1,\alpha_2$ is a routine adaptation of the argument in \cite{C98}.  In brief, we simply observe that by Fubini, the set of $(\xi,\eta) \in \Omega^{(j-1)}$ for which the inequality
$$
|\pi_2(\Omega^{(j-1)})| \cdot |\{u \in \R^2 : (\xi+u,\eta) \in \Omega^{(j-1)}\}| < \tfrac18|\Omega^{(j-1)}|
$$
holds has volume less than $\tfrac18|\Omega^{(j-1)}|$.  The complement of this ``bad'' set in $\Omega^{(j-1)}$ is $\Omega^{(j)}$.  The upper bound on $|u|$ comes from our assumption that $\diam(\Omega)| < C(\alpha_1\alpha_2)^\rho$.  

We turn now to the statement involving $\alpha_3$.  We note that the vector field $Z_3$ is tangent to the fibers of $\pi_3$, and $\pi_3$ has full rank when $Z_3$ is nonvanishing.  We also observe that $Z_3$ has size comparable to one on $R_1 \times R_2$.  The statement involving $\alpha_3$ may then be proved by noting that, by examining standard proofs of the implicit function theorem (e.g., \cite{R76}) and our hypotheses, there exists a global section $\sigma_3$ from some neighborhood of $\pi_3(R_1 \times R_2)$ to a neighborhood of $R_1 \times R_2$, and $(t,y) \mapsto e^{tZ_3}(\sigma_3(y))$ is a diffeomorphism of some open subset of $\R^4$ onto $R_1 \times R_2$ with Jacobian determinant comparable to one.  We then use Fubini and deletion of small bad sets as above.  
\end{proof}

There are many proofs of the inequality $|\Omega| \lesssim |\pi_1(\Omega)||\pi_2(\Omega)|$, or, equivalently, $|\Omega| \gtrsim \alpha_1\alpha_2$.  The most immediately relevant of these is that, given a point $(\xi_0,\eta_0) \in \Omega^{(2)}$, we may construct a parameter set consisting of those $(u,v)$ such that $(\xi_0+u,\eta_0) \in \Omega^{(1)}$ and $(\xi_0+u,\eta_0+v) \in \Omega$.  This set has measure $\gtrsim \alpha_1\alpha_2$, and hence so does its image under $(u,v) \mapsto (\xi_0+u,\eta_0+v)$.  A variant of this argument permits us to prove a nontrivial estimate.  

\begin{lemma} \label{L:alphas squared}
\begin{equation} \label{E:alphas squared}
|\Omega| \gtrsim \alpha_1^2\alpha_3^2 \qtq{and} |\Omega| \gtrsim \alpha_2^2\alpha_3^2.
\end{equation}
\end{lemma}


\begin{proof}[Proof of Lemma~\ref{L:alphas squared}]
We will prove the second of these inequalities, as it implies the first by interchanging the roles of $\xi$ and $\eta$.  

Fix $(\xi_0,\eta_0) \in \Omega^3$.  Dimensional considerations suggest the construction of a map
\begin{gather*}
\Phi_{(\xi_0,\eta_0)} : \R \times \R^2 \times \R \to \R^4\\
(t_1,u,t_2) \mapsto e^{t_2Z_3}(e^{t_1Z_3}(\xi_0,\eta_0)+(u,0)),
\end{gather*}
on the set $W$ of triples $(t_1,u,t_2)$ with $e^{t_1Z_3}(\xi_0,\eta_0) \in \Omega^2$, $e^{t_1Z_3}(\xi_0,\eta_0)+(u,0) \in \Omega^1$, and $\Phi_{(\xi_0,\eta_0)} \in \Omega$.  Then $W$ has measure at least $\alpha_3^2\alpha_2$ and diameter at most $C(\alpha_1\alpha_2)^\rho$.  The Jacobian determinant of $\Phi_{(\xi_0,\eta_0)}$ we readily compute as 
\begin{equation} \label{E:det DPhi}
|\det D\Phi_{(\xi_0,\eta_0)}(t_1,u,t_2)| = |\det(\nabla \psi(\xi_1+u)-\nabla\psi(\eta_1),\nabla \psi(\xi_1)-\nabla\psi(\eta_1)|,
\end{equation} 
where $(\xi_1,\eta_1):=e^{t_1Z_3}(\xi_0,\eta_0)$.  

By assumption, $\nabla \psi$ is a diffeomorphism with bounded derivatives and Jacobian determinant comparable to 1.  Thus, for each $t_1$ with $(\xi_1,\eta_1) \in \Omega^2$, the set of 
$$
\{\nabla \psi(\xi_1+u) : (\xi_1+u,\eta_1) \in \Omega^1\}
$$
is a set of measure $\gtrsim \alpha_2$ and diameter $\leq 1$.  On the other hand, given such $t_1$, the set of such $\nabla \psi(\xi_1+u)$ such that $|\det(\nabla \psi(\xi_1+u)-\nabla\psi(\eta_1),\nabla \psi(\xi_1)-\nabla\psi(\eta_1)| \ll \alpha_1$ is (since it has diameter $\leq 1$) contained in a set of measure $\ll \alpha_1$.  Thus, on most of $W$ (by measure),
\begin{equation} \label{E:jacobian lb holds}
|\det D\Phi_{(\xi_0,\eta_0)}(t_1,u,t_2)| \gtrsim \alpha_1.  
\end{equation}
We let $\tilde W$ denote the set of $(t_1,u,t_2) \in W$ such that \eqref{E:jacobian lb holds} holds.  We apply Lemma~\ref{L:poly approx} with $B$ the ball of radius $C(\alpha_1\alpha_2)^\rho$ centered at zero and $\Psi$ the degree $2/\rho$ (or higher) Taylor polynomial of $\Phi_{(\xi_0,\eta_0)}$, 
$$
|\Omega| \geq |\Phi_{(\xi_0,\eta_0)}(W)| \gtrsim \int_{\tilde W} |\det D\Phi_{(\xi_0,\eta_0)}(t_1,u,t_2)| \, dt_1\, du\, dt_2 \gtrsim \alpha_1|\tilde W| \gtrsim \alpha_1|W| \gtrsim \alpha_1^2\alpha_3^2.
$$
\end{proof}

Having proved 
$$
|\Omega| \gtrsim \max\{\alpha_1\alpha_2,\alpha_1^2\alpha_3^2,\alpha_2^2\alpha_3^2\},
$$
we have also 
\begin{equation} \label{E:rwt strichartz}
|\Omega| \gtrsim (\alpha_1\alpha_2)^{1/2}(\alpha_1\alpha_3)^{1/2}(\alpha_2\alpha_3)^{1/2} = \alpha_1\alpha_2\alpha_3,
\end{equation}
which, of course, we already knew.  

We wish to beat \eqref{E:rwt strichartz} by only a little bit and with a bit of arithmetic we see that Theorem \ref{T:bilinear rest} follows from proving
$$
|\Omega| \gtrsim \min\{\alpha_1^{1-\eps_n}\alpha_2,\alpha_1\alpha_2^{1-\eps_n}, \alpha_1^{2-\delta_n}\alpha_3^2,\alpha_1^2\alpha_3^{2-\delta_n}\},
$$
with 
\begin{equation} \label{E:def eps-n delta-n}
\eps_n> 2(2-s), \qquad \delta_n>4(2-s).  
\end{equation}.

Equivalently, it suffices to prove the conclusion of Lemma~\eqref{L:isoperim reduction} under the assumption that
\begin{equation} \label{E:Omega too small}
|\Omega| \ll \min\{\alpha_1^{1-\eps_n}\alpha_2,\alpha_1\alpha_2^{1-\eps_n}, \alpha_1^{2-\delta_n}\alpha_3^2,\alpha_2^{2-\delta_n}\alpha_3\}.  
\end{equation}
Assumption \eqref{E:Omega too small} will remain in force for the remainder of the proof of Lemma~\ref{L:isoperim reduction}.  

Let us suppose for the moment that there existed $\tilde\Omega^2 \subseteq \Omega^2$ such that $|\tilde\Omega^2| \gtrsim |\Omega^2|$ and for all $(\xi,\eta) \in \tilde\Omega^2$, 
$$
|\{u : (\xi+u,\eta) \in \Omega^1 \ctc{and} |\det(\nabla\psi(\xi+u)-\nabla\psi(\eta),\nabla\psi(\xi)-\nabla\psi(\eta))| \gtrsim \alpha_2^{1-\delta_n}\}| \gtrsim \alpha_2,
$$
then we have a slight improvement on \eqref{E:jacobian lb holds}, and, correspondingly adapting the proof of the lemma, we would have
$$
|\Omega| \gtrsim \alpha_2^{2-\delta_n}\alpha_3^2,
$$
which contradicts \eqref{E:Omega too small}.  

Thus, there exists $\tilde\Omega^2 \subseteq \Omega^2$ of size $|\tilde\Omega^2| \gtrsim |\Omega^2|$ such that for all $(\xi,\eta) \in \tilde\Omega^2$, 
$$
|\{u : (\xi+u,\eta) \in \Omega^1 \ctc{and} |\det(\nabla\psi(\xi+u)-\nabla\psi(\eta),\nabla\psi(\xi)-\nabla\psi(\eta))| \ll \alpha_2^{1-\delta_n}\}| \gtrsim \alpha_2,
$$
and likewise with the roles of $\xi,\eta$ reversed and $\alpha_2$ replaced by $\alpha_1$.  

We now introduce a new map.  Let $(\xi_0,\eta_0) \in \Omega^3$.  We know that for ``most'' $u$ such that $(\xi_0+u,\eta_0) \in \Omega^2$, $\nabla\psi(\xi_0+u)-\nabla\psi(\xi_0)$ lies within $\ll \alpha_2^{1-\delta_n}$ of the line in the direction $\nabla \psi(\xi_0)-\nabla\psi(\eta_0)$.  Since $\nabla \psi$ is a diffeomorphism on $2Q$, we may express such $u$ as 
$$
u=u(\xi_0,\eta_0,\sigma,\theta) := (\nabla \psi)^{-1}(\nabla \psi(\xi_0) + \sigma(\nabla \psi(\xi_0)-\nabla \psi(\eta_0))) - \xi_0 + \theta \omega_{(\xi_0,\eta_0)},
$$
where $|\sigma| \lesssim (\alpha_1\alpha_2)^\rho$, $|\theta| \ll \alpha_2^{1-\delta_n}$, and $\omega_{(\xi_0,\eta_0)}$ is a fixed unit vector transversal to the curve
$$
\{(\nabla \psi)^{-1}(\nabla \psi(\xi_0) + \sigma(\nabla \psi(\xi_0)-\nabla \psi(\eta_0))) : |\sigma| \ll 1\}
$$
at each point on the curve.  (Here it is useful to recall that $|\xi_0-\eta_0| \sim |\nabla \psi(\xi_0)-\nabla \psi(\eta_0)| \sim 1$ and that $D^2\psi$ is nondegenerate, so we know that the curve has a well-defined and continuous tangent vector at each point.)

With this definition in place, we define
\begin{align*}
&\Psi_{(\xi_0,\eta_0,\theta)}: \R \times \R^2 \times \R \to \R^4\\
&\Psi_{(\xi_0,\eta_0,\theta)}(\sigma,v,t) :=e^{tZ}(\xi_0+u(\xi_0,\eta_0,\sigma,\theta),\eta_0+v).
\end{align*}
Setting $(\xi_1,\eta_1):=(\xi_0+u(\xi_0,\eta_0,\sigma,\theta),\eta_0+v)$, we may compute the Jacobian determinant of $\Psi_{(\xi_0,\eta_0,\theta)}$:
\begin{equation} \label{E:det DPsi}
|\det D\Psi_{(\xi_0,\eta_0,\theta)}(\sigma,v,t)| \sim |(\nabla \psi(\xi_1)-\nabla\psi(\eta_1))^T D^2\psi(\xi_1)^{-1}(\nabla \psi(\xi_0)-\psi(\eta_0))|.
\end{equation}
Recalling that $|\theta| \ll \alpha_2^{1-\delta_n}$ and the definition of $\tilde\Omega^2$,
$$
\xi_1=\xi_0+u = (\nabla \psi)^{-1}(\nabla \psi(\xi_0) + \sigma(\nabla \psi(\xi_0) - \nabla \psi(\eta_0)))+O(\alpha_2^{1-\delta_n}),
$$
and 
\begin{align*}
\eta_1&=(\nabla \psi)^{-1}(\nabla \psi(\eta_0) + \tilde\sigma(\nabla \psi(\eta_0) - \nabla \psi(\xi_1)))+O(\alpha_2^{1-\delta_n}),\\
&=\nabla \psi^{-1}(\nabla \psi(\eta_0) + \tilde\sigma(\nabla\psi(\eta_0) - \nabla \psi(\xi_0) - \sigma(\nabla \psi(\xi_0)-\nabla\psi(\eta_0))+O(\alpha_2^{1-\delta_n}).
\end{align*}
Therefore, after some arithmetic, and in view of \eqref{E:det DPsi}, 
\begin{align}\notag
&|\det D\Psi_{(\xi_0,\eta_0,\theta)}(\sigma,v,t)|
\sim | (\nabla \psi(\xi_0)-\nabla \psi(\eta_0)^T\\\notag
&\qquad \times D(\nabla \psi)^{-1}(\nabla \psi(\xi_0)+\sigma(\nabla \psi(\xi_0)-\nabla \psi(\eta_0)))(\nabla \psi(\xi_0)-\nabla \psi(\eta_0))| + O(\alpha_2^{1-\delta_n}) \\ \label{E:no vt}
&\sim |\det D\Psi_{(\xi_0,\eta_0,\theta)}(\sigma,0,0)| + O(\alpha_2^{1-\delta_n}).
\end{align}

Let us now consider two quantities:  
\begin{align*}
A(t)&:=(\xi_0-\eta_0)^TD^2\psi(\xi_0+t(\xi_0-\eta_0))(\xi_0-\eta_0)\\
B(\sigma)&:=(\nabla \psi(\xi_0)-\nabla \psi(\eta_0)^TD(\nabla \psi)^{-1}(\nabla \psi(\xi_0)+\sigma(\nabla \psi(\xi_0)-\nabla \psi(\eta_0)))\\
&\qquad \times (\nabla \psi(\xi_0)-\nabla \psi(\eta_0)).  
\end{align*}

\begin{lemma} \label{L:sup B}
For each $(\xi_0,\eta_0) \in \Omega^3$,
\begin{equation} \label{E:sup B}
\sup_{\sigma \in [0,1]} |B(\sigma)| \sim 1.
\end{equation}
\end{lemma}

\begin{proof}[Proof of Lemma~\ref{L:sup B}]
We have assumed that $\xi_0 \in R_1,\eta_0 \in R_2$ are $(1,1)$-separated; thus, $|\xi_0-\eta_0| \sim 1$ and 
\begin{equation} \label{E:sup A}
\sup_{t \in [-1,2]}|A(t)| \sim 1.
\end{equation}
In view of \eqref{E:polynomial like} and 
$$
A(t) = \tfrac{d}{dt}[\nabla \psi(\xi_0+t(\eta_0-\xi_0)) - \nabla \psi(\xi_0)] \cdot (\xi_0-\eta_0),
$$
\eqref{E:sup A} implies that the maximum distance from $\nabla \psi([\xi_0,\eta_0])$ to the line through $\nabla \psi(\xi_0)$ in the direction $(\xi_0-\eta_0)^\perp$ is also comparable to 1.  That is,
\begin{equation} \label{E:line dist 1}
\dist(\nabla \psi([\xi_0,\eta_0]),\nabla \psi(\xi_0) + (\xi_0-\eta_0)^\perp \R) \sim 1.
\end{equation}
 (Here we use the standard notation $[a,b]$ to denote the line segment with endpoints $a,b$.)

If \eqref{E:sup B} were to fail, then by the mean value theorem and 
$$
B(\sigma) = \tfrac{d}{d\sigma}[(\nabla \psi)^{-1}(\nabla \psi(\xi_0)+\sigma(\nabla \psi(\eta_0)-\nabla \psi(\xi_0))) - \xi_0] \cdot (\nabla \psi(\xi_0)-\nabla \psi(\eta_0)),
$$
we have
\begin{equation} \label{E:line dist 2}
\dist((\nabla \psi)^{-1}([\nabla \psi(\xi_0),\nabla \psi(\eta_0)]), \xi_0 + (\nabla \psi(\xi_0)-\nabla \psi(\eta_0))^\perp) \ll 1.
\end{equation}
(Here we are using the asymmetric set distance, so that this is the maximum distance from a point in the first set to the second set.)    Noting that the endpoints of $(\nabla\psi)^{-1}([\nabla \psi(\xi_0),\nabla \psi(\eta_0)])$ are $\xi_0,\eta_0$, this implies that 
$$
\dist([\xi_0,\eta_0], \nabla \psi(\eta_0)]), \xi_0 + (\nabla \psi(\xi_0)-\nabla \psi(\eta_0))^\perp) \ll 1.
$$
This implies, in addition, $|(\xi_0,\eta_0) \cdot (\nabla \psi(\xi_0)-\nabla \psi(\eta_0))| \ll 1$.  Thus \eqref{E:line dist 2} implies that
$$
\dist([\xi_0,\eta_0],(\nabla \psi)^{-1}(\nabla \psi(\xi_0) + (\xi_0-\eta_0)^\perp)) \ll 1,
$$
which directly contradicts \eqref{E:line dist 1}, since $\nabla \psi$ is bi-Lipschitz.
\end{proof}

%

We are now in a position to use \eqref{E:det DPsi}.  With $(\xi_0,\eta_0) \in \Omega^3$ fixed, let us now consider, for $|\theta| \ll \alpha_2^{1-\delta_n}$,  
\begin{align*}
U_\theta&:=\{(\sigma,v,t) : |(\sigma,v,t)| \lesssim (\alpha_1\alpha_2)^\rho, \Psi_{(\xi_0,\eta_0,\theta)}(\sigma,0,0) \in \Omega^2,\\
&\qquad \Psi_{(\xi_0,\eta_0,\theta)}(\sigma,v,0) \in \Omega^1, \Psi_{(\xi_0,\eta_0,\theta)}(\sigma,v,t) \in \Omega\}.  
\end{align*}
Noting that $(\theta,\sigma) \mapsto u(\xi_0,\eta_0,\sigma,\theta)$ is a diffeomorphism, and $|\{u : (\xi_0+u,\eta_0) \in \Omega^2\}| \gtrsim \alpha_2$, there exists $|\theta| \lesssim \alpha_2^{1-\delta_n}$ such that $|S_\theta|:=|\{|\sigma| \lesssim (\alpha_1\alpha_2)^\rho : \Psi_{(\xi_0,\eta_0,\theta)} \in \Omega^2\}| \gtrsim \alpha_2^{\delta_n}$.  With this value of $\theta$, $|U_\theta| \gtrsim \alpha_2^{\delta_n}\alpha_1 \alpha_3$.  With $(\xi_0,\eta_0,\theta)$ now fixed, we henceforth omit them from the subscript on $\Psi$.

Let $\Psi_{vlo}$ denote the degree $n$ Taylor polynomial of $\Psi$ centered at 0, $\Psi_{lo}$ the degree $N$ polynomial, $\Psi_{med} := \Psi_{lo}-\Psi_{vlo}$, $\Psi_{hi}:=\Psi-\Psi_{vlo}$, and $\Psi_{vhi} := \Psi-\Psi_{lo}$.  By \eqref{E:polynomial approx}, $\|\Psi_{hi}(\cdot,0,0)\|_{C^1([0,1])} \ll 1$.  On the other hand, by \eqref{E:sup B} and \eqref{E:no vt}, we have seen that $\|\det D\Psi(\cdot,0,0)\|_{C^0([0,1])} \sim 1$, so $\|\det D\Psi_{vlo}(\cdot,0,0)\|_{C^0([0,1])} \gtrsim 1$.  As $\det D\Psi_{vlo}$ is a polynomial of degree $2(n-1)$, $\sum_{k \leq 2(n-1)} |\partial_\sigma^k \det D\Psi_{vlo}(\sigma,0,0)| \gtrsim 1$, $\sigma \in [0,1]$.  Therefore, $\sum_{k \leq 2(n-1)} |\partial_{\sigma}^k \det D\Psi_{lo}(\sigma,0,0)| \gtrsim 1$, $\sigma \in [0,1]$.  Therefore $|\det D\Psi_{lo}(\sigma,0,0)| \gtrsim \alpha_2^{2(n-1)\delta_n}$ on most of $S_\theta$.  By construction, 
\begin{equation} \label{E:Psi hi small}
\|\Psi-\Psi_{lo}\|_{C^2(\{|(\sigma,v,t)| \lesssim (\alpha_1\alpha_2)^\rho\})} \lesssim (\alpha_1\alpha_2)^{\rho N},
\end{equation}
and for $\rho$ sufficiently small.  Thus, 
$$
\|\det D\Psi_{lo}-\det D\Psi\|_{C^0(\{|(\sigma,v,t)| \lesssim (\alpha_1\alpha_2)^\rho\})} \lesssim (\alpha_1\alpha_2)^{\rho N}.
$$
For $\rho$ sufficiently small, $2(n+1)\delta_n < \rho N$, which implies that
$$|\det D\Psi| \gtrsim \alpha_2^{2(n-1)\delta_n}
$$
on most of $U_\theta$.  By \eqref{E:Psi hi small} again, we may thus apply Lemma~\ref{L:poly approx} to see that
$$
|\Omega| \geq |\Psi(U_\theta)| \gtrsim \int_{U_\theta}|\det D\Psi| \gtrsim |U_\theta|\alpha_2^{2(n-1)\delta_n} \gtrsim \alpha_1\alpha_2^{(2n-1)\delta_n}\alpha_3.
$$
By symmetry, we also obtain
$$
|\Omega| \gtrsim \alpha_1^{(2n-1)\delta_n}\alpha_2\alpha_3,
$$
whence
\begin{equation} \label{E:Omega too big}
|\Omega| \gtrsim \alpha_1^{\frac{(2n-1)\delta_n}2+\frac12} \alpha_2^{\frac{(2n-1)\delta_n}2+\frac12} \alpha_3.
\end{equation}
On the other hand, by \eqref{E:Omega too small},
$$
|\Omega| \ll \alpha_1^{1-\frac{\eps_n}4-\frac{\eps_n}4}\alpha_2^{1-\frac{\eps_n}4-\frac{\eps_n}4}\alpha_3,
$$
which contradicts \eqref{E:Omega too big} when 
\begin{equation} \label{E:def2 eps-n delta_n}
\frac{\eps_n}4 + (n-\frac14)\delta_n \leq \frac12.
\end{equation} 
Finally, there exist $\eps_n,\delta_n$ obeying \eqref{E:def eps-n delta-n} and \eqref{E:def2 eps-n delta_n} whenever $(2-s) \leq \tfrac1{8n}$.  
\end{proof}

\section{Linear restriction for sets of constant fiber length}
    
We begin our bilinear-to-linear reduction by proving a restricted weak-type bound under further assumptions on the domains for which the restricted weak-type bounds apply. Our assumption will be that the domain has essentially constant fiber length for the $N_1$-flows.

\begin{definition}
  For $\Omega \subset [-1,1]^2$, define and $K \in \N$, we define \[\Omega(K) = \{\xi \in \Omega: \mathcal{H}^1(\{t : e^{tN_1}(\xi) \in \Omega\}) \sim 2^{-K}\}.\]We say that $\Omega$ has constant $N_1$-fiber length if $\Omega = \Omega(K)$. We definite sets of constant $N_2$-fiber length analogously. If $\Omega$ has constant $N_1$-fiber length or constant $N_2$-fiber length, we say $\Omega$ has \emph{constant fiber length}.
\end{definition}

We have a scaling line linear restriction estimate in the same range as our bilinear restriction estimates for sets of constant fiber length. 

\begin{proposition}\label{P:RR}
  Suppose $s > \frac{3}{2}$ and for any $r$ in a neighborhood of $s'$, Corollary \ref{C:rescaled bilinear rest} holds.
  
  Then for any $\Omega \subset [-1,1]^2$ of constant fiber length, 
  \begin{equation}\label{E:restricted restricted linear}||E1_{\Omega'}||_{2s} \lesssim |\Omega(K)|^{1/s'}\end{equation}for any measurable $\Omega' \subseteq \Omega(K)$. 
\end{proposition}
Proposition \ref{P:RR} follows essentially the same reasoning as earlier work of Vargas for the hyperbolic paraboloid \cite{V05}. The novel part of the proof Proposition \ref{P:RR} is a Whitney decomposition adapted to the surface. Vargas chose axis-parallel rectangles for the Whitney decomposition for the hyperbolic paraboloid, that is rectangles of the form $I_j \times I_k$ for dyadic intervals $I_j$ of length $2^{-j}$ and $I_k$ of length $2^{-k}$. Since the projections $\pi_1, \pi_2$ associated with the hyperbolic paraboloid are simply coordinate projections, a natural generalization to other surfaces be tiles of the form $\tau_{j,k} = \pi_1^{-1}(I_j) \cap \pi_2^{-1}(I_k)$. When $j$ and $k$ are near to each other, these tiles form Euclidean rectangles, but they may gain non-trivial curvature as the dimensions separate. An important part of our argument is that the sum set of tiles paired in the Whitney decomposition must be boundedly overlapping to efficiently decompose the extension operator over the tiles and we could only prove that was the case when the tiles were essentially rectangles, so we must further refine our Whitney decomposition. 

Fortunately, the separation condition in the bilinear extension estimate in Corollary \ref{C:rescaled bilinear rest} improves when the tiles are far from being Euclidean rectangles. We can use this improvement to group together tiles that are too narrow to be Euclidean rectangles to form one larger tile which is a Euclidean rectangle and still satisfies a useful bilinear extension estimate. We give the details in the following subsection.

\subsection{A Whitney decomposition adapted to $\psi$}
  Denote by $\mathcal{I}_{j}$ a partition of $\R$ into dyadic $2^{-j}$-intervals. For $I_j, I'_j \in \mathcal{I}_j$, we write $I_j \sim I'_j$ if the nearest points of $I_j$, $I'_j$ differ by $2^{-j+1}$ (in other words, if the convex hull of $I_j \cup I'_j$ is an element of $\mathcal{I}_{j-2}$). Any pair $x, y \in \R$ satisfies $(x,y) \in I_j \times I'_j$ for precisely one choice of $j$ and one choice of $I_j, I'_j$, specifically $j = - \lceil\log_2 |x-y|\rceil-1$ and the choices of $I_j, I'_j$ respectively containing $x, y$. 

  Recall the projections $\pi_1, \pi_2:[-1,1]^2 \to \R$ defined in Subsection \textbf{Rectangles}. For any pair $j, k$, we partition \[Q = \bigcup_{\tau \in \mathcal{T}_{j,k}} \tau,\]where $ \mathcal{T}^0_{j,k} = \{\pi_1^{-1}(I_j) \cap \pi_2^{-1}(I_k) : I_j \in \mathcal{I}_{j}, I_k \in \mathcal{I}_{k}\}$. For choices of $\tau=\pi_1^{-1}(I_j) \cap \pi_2^{-1}(I_k), \tau' = \pi_1^{-1}(I_j') \cap \pi_2^{-1}(I_k')$, we write $\tau \sim_0 \tau'$ if $I_j' \sim_0 I_j$ and $I_k' \sim_0 I_k$. For any pair $(\xi, \eta)$ for which $\pi_1(\xi) \neq \pi_1(\eta)$ and $\pi_2(\xi) \neq \pi_2(\eta)$, we have exactly one pair of scales $(j,k)$ and one pair $\tau, \tau' \in \mathcal{T}^0_{j,k}$ with $\tau \sim_0 \tau'$ and $(\xi, \eta) \in \mathcal{T}^0_{j,k}$. The points $(\xi, \eta)$ with $\pi_1(\xi) = \pi_1(\eta)$ or $\pi_2(\xi) = \pi_2(\eta)$ form a null set in $\R^4$ and hence will play no role in any of our estimates. Going forward, we disregard them. The tiles $\mathcal{T}^0_{j,k}$ and the relation $\sim_0$ give the naive attempt at a Whitney decomposition described previously.

  Next, we expand the set of possible tiles to form the basis for our actual Whitney decomposition. For $k \ge j+2$, set $\mathcal{T}_{j,k} = \mathcal{T}^0_{j,k} \cup \mathcal{T}^0_{j, k-2}$. If $j \ge k + 2$, we similarly define $\mathcal{T}_{j,k} = \mathcal{T}^0_{j,k} \cup \mathcal{T}^0_{j-2, k}$. For the diagonal terms, we simply set $\mathcal{T}_{j,k} = \mathcal{T}^0_{j,k}$. While $\sim_0$ still gives a Whitney decomposition for $\bigcup_{j,k \in \N}\mathcal{T}_{j,k}$ since it does not see any of the additions we made, we need to further refine this Whitney decomposition to reflect the geometry of the surface. The diagonal terms $|k -j | \le 2$ are always essentially Euclidean rectangles, so we set $\tau \sim \tau'$ if $\tau \sim_0 \tau'$. For the off-diagonal case, we treat the $k \ge j + 2$ terms and $j \ge k+2$ symmetrically, so we only cover the $k \ge j+2$ case. 
  
  Consider two dyadic intervals $I_j$ and $I'_j$ with $I_j \sim I'_j$. Let $\mathcal{T}_{j,k}(I_j) = \{\tau \in \mathcal{T}_{j,k} : \pi_1(\tau) = I_j\}$. We will define a relation $\sim$ starting with by iterating over the scales $k \ge j + 2$. Assume we have defined the relation for all $k \in [j+2, k_0) \cap \N$. At scale $k_0$, we set $\tau \sim \tau'$ if the following conditions are met:

  \begin{enumerate}[label = (C\arabic*)]
    \item $\tau \in \mathcal{T}^0_{j,k_0}(I_j)$, $\tau' \in \mathcal{T}^0_{j,k_0}(I'_j)$, and $\tau \sim_0 \tau'$.
    \item There does not exist $\tau_*, \tau'_*$ with $\tau_* \sim \tau'_*$, $\tau \subset \tau_*$, and $\tau'\subset \tau'_*$. 
    \item $k_0 \le k^2(\xi, j-C)$ for $\xi$ the center point of $\tau$. 
  \end{enumerate}
  If (C1) is true, (C2) is true, but (C3) is false, we also add $\tilde{\tau} \sim \tilde{\tau}'$, where \begin{equation}\label{E:widened tiles}\tilde{\tau} = \pi_1^{-1}(I_j) \cap \pi_2^{-1}(\tilde{I}_k), \tilde{\tau}' = \pi_1^{-1}(I_j) \cap \pi_2^{-1}(\tilde{I}_k),\end{equation}and $\tilde{I}_k \in \mathcal{T}_{k-2}$ is the convex hull of $\pi_2(\tau)$ and $\pi_2(\tau')$. Note that $\tilde{\tau}\in \mathcal{T}_{j,k_0}(I_j), \tilde{\tau}' \in \mathcal{T}_{j,k_0}(I'_j)$, so adding a relation between the pair makes sense. 
  
  In other words, we start with a pair $\tau \sim_0 \tau'$ from the naive set of tiles (this is condition C1) and we ask in C2 if we have already defined a pair of tiles covering $\tau \times \tau'$. If we have not defined such a pair and C3 is true then $\tau, \tau'$ are essentially Euclidean so we take $\tau \sim \tau'$. If C3 is false, then we widen $\tau$ and $\tau'$ to have the same $\pi_2$ projection. We will verify that C2 and $\neg$C3 implies that these widened tiles are essentially Euclidean and satisfy the right separation condition to apply the bilinear extension estimate. Lemma \ref{L:comparable rectangles} tells us that the condition $k \le k^1(\xi, j-C)$ implies that a tile $\tau \in \mathcal{T}_{j,k}$ centered at $\xi$ is essentially Euclidean; this is the condition that will appear in our construction. 
  
  Here and henceforth $C$ is an admissible value allowed to change from line to line.

  \begin{lemma}\label{L:whitney props}
    For any pair $j, k \in \N^2$, the following hold
    \begin{enumerate}[label = (\roman*)]
      \item For any $\tau, \tau' \in \mathcal{T}_{j,k}$ with $\tau \sim \tau'$, if $|\pi_2(\tau)| = 2^{-k'}$, then $k' \le k^2(\xi, j-C)$ for $\xi$ the center point of $\tau$.
      \item For any $\tau_0, \tau'_0 \in \mathcal{T}^0_{j,k}$ with $\tau_0 \sim_0 \tau'_0$, for any $(\xi, \eta) \in \tau_0 \times \tau'_0$, there is exactly one choice of $(j', k')$ and (up to reordering) one pair $\tau_0, \tau'_0 \in \mathcal{T}_{j',k'}$ with $\tau_0 \sim \tau'_0$ and $(\xi, \eta) \in \tau_0 \times \tau'_0$. Moreover, for any $\tau$, there are $\lesssim 1$ $\tau'$ such that $\tau \sim \tau'$.
      \item For any $\tau, \tau' \in \mathcal{T}_{j,k}$ with $\tau \sim \tau'$, each point $(\xi, \eta) \in \tau \times \tau'$ is $(2^{-k}, 2^{-j})$-separated.
    \end{enumerate}
  \end{lemma}

  \begin{proof}
    We break this proof into the diagonal case $|j - k| < 2$ and the off-diagonal case $|j - k| \ge 2$. For the off-diagonal case, we assume $k \ge j+2$, since the opposite case is the same after swapping the roles of $N_1$ and $N_2$. 

    \vspace{4pt}

    \noindent\textbf{The diagonal case}. 

    \vspace{4pt}

    In the diagonal case, each $\tau \in \mathcal{T}_{j,k}$ has $\pi_2(\tau) = 2^{-k}$, so (i) follows from the fact that $k^2(\xi, j-C) \ge j + 2 \ge k$. For (ii), we simply take $\tau = \tau_0, \tau' = \tau'_0$. We already know that this pair is unique up to reordering. It is an immediate consequence of the definition $\sim_0$ that each $\tau$ has $\le 4$ $\tau'$ with $\tau \sim \tau'$. Hence, we have (ii) 
    
    For each $(\xi_0, \eta_0) \in \tau \times \tau'$, we certainly have $|\xi_0 - \eta_0| \sim 2^{-j}$. For the other separation element, we compute \begin{align}2^{-k} \sim |\pi_2(\eta_0) - \pi_2(\xi_0)| &= \left|\int_0^{1} \langle d\pi_2(\xi_0 + t(\eta_0 - \xi_0)), (\eta_0 - \xi_0)\rangle~dt\right|\\ &= \left|\int_0^1 N_2(\xi_0 + t(\eta_0 - \xi_0)) \wedge (\eta_0 - \xi_0) ~dt\right| \label{E:pi_2 separation?}.\end{align}Hence, $\sup_{t \in [-1, 2]} |N_2(\xi_0 + t(\eta_0 - \xi_0)) \wedge (\eta_0 - \xi_0)| \sim 2^{-k}$, so $(\xi_0, \eta_0)$ are $(2^{-k}, 2^{-j})$ separated. Hence, we have (iii).

    \vspace{4pt}

    \noindent\textbf{The off-diagonal case (i)}.

    \vspace{4pt}

    If $\tau \sim \tau'$ is a pair satisfying C1, C2, and C3, then we must have the desired conclusion. Otherwise, $\tau, \tau'$ is given by widening a pair $\tau_0, \tau'_0 \in \mathcal{T}_{j,k}^0$ satisfying C1 and C2 and $k' = k-2$. If $k-2 \le k^2(\xi_0, j-C)$, then we have the desired conclusion for the widened tiles $\tilde{\tau}, \tilde{\tau}'$ defined in (\ref{E:widened tiles}). Instead, suppose $k-2 > k^2(\xi_0, j-C)$ for $\xi_0$ the center of $\tau_0$. Let $\tau_1, \tau'_1 \in \mathcal{T}_{j,k-1}^0$ be the tiles adjacent to the outside of $\tau_1, \tau'_1$ satisfying $\tau_1 \sim \tau'_1$. For an example of what we mean, if $\pi_2(\tau_0)= [0, 2^{-k}]$ and $\pi_2(\tau'_0) = [2^{-k} + 2^{-k+1}, 2^{-k+2}]$, we would take $\pi_2(\tau_1) = [-2^{-k+1}, 0]$ and $\pi_2(\tau'_1) = [2^{-k+2}, 2^{-k+2} + 2^{-k+1}]$. Up to translation and swapping $\tau_0, \tau'_0$, this example characterizes the construction of $\tau_1, \tau'_1$ no matter our initial choice of $\tau_0, \tau'_0$. 

    It is not hard to see that $\xi \mapsto 2^{-k^2(\xi, j-C)}$ is Lipschitz with constant $O(||\psi||_{C^3(Q)})$. We may assume $||\psi||_{C^3(Q)} < c$, where $c$ is an absolute constant is chosen so that the Lipschitz constant of $\xi \mapsto 2^{-k^2(\xi, j-C)}$ is $\le \frac12$. If $\xi_1$ is the center of $\tau_1$, we know that $|\xi_0 - \xi_1| \le 2^{-k} + 2^{-k-1}$, and hence \[|2^{-k^2(\xi_0, j-C)} - 2^{-k^2(\xi_1, j-C)}| \le \frac{1}{2}|\xi_0 - \xi_1| = 2^{-k-1} + 2^{-k-2}.\]Since we assume $k_2(\xi_0, j-C) < k - 2$, we see that $k^2(\xi_1, j-C) < k-1$. We conclude that C3 must fail for $\tau_1, \tau'_1$. Since C2 succeeds for $\tau_0, \tau'_0$, C2 must also succeed for $\tau_1, \tau'_1$ and C1 since $\tau_1 \sim \tau'_1$. Hence, the widening $\tilde{\tau}_1, \tilde{\tau}'_1$ of $\tau_1, \tau'_1$ must have $\tilde{\tau}_1 \sim \tilde{\tau}'_1$, but $\tau_0 \times \tau'_0 \subset \tilde{\tau}_1 \times \tilde{\tau}'_1$, so we contradict C2 for $\tau_0, \tau'_0$.

    If $k = j+2$, this argument would fail because the pairs $\tau_1, \tau'_1$ would not come from the off-diagonal construction, but $k^2(\xi, j-C) \ge j+2$ by (\ref{E:ki lower bound}), giving the desired result.

    \vspace{4pt}
    
    \noindent\textbf{The off-diagonal case (ii)}

    \vspace{4pt}

    First, we prove that each point $(\xi, \eta)$ has a pair $\tau \sim \tau'$ with $(\xi, \eta) \in \tau \times \tau'$. Let $\tau_0 \sim_0 \tau'_0$ be the pair in $\mathcal{T}^0_{j,k}$ for which $(\xi, \eta) \in \tau_0 \times \tau'_0$. C1 always succeeds for the pair $\tau_0 \times \tau'_0$. If C2 and C3 succeed, then $\tau_0 \sim \tau'_0$ so we take $\tau = \tau_0, \tau' = \tau'_0$. If C2 fails, then exists a superset $\tau_*, \tau'_*$ with $\tau_* \sim \tau'_*$, so we take $\tau = \tau_*, \tau' = \tau'_*$. If C3 fails, then $(\xi, \eta) \in \tilde{\tau} \times \tilde{\tau}'$, where that pair is as defined in (\ref{E:widened tiles}). In any case, we have existence of a pair with the desired consequences. 

    Now, we prove uniqueness. Suppose $(\xi, \eta) \in (\tau_0 \times \tau'_0) \cap (\tau_1 \times \tau'_1)$, with $\tau_i, \tau'_i \in \mathcal{T}_{j,k_i}$ and $\tau_i \sim \tau'_i$ for both pairs. Assume $k_1 \ge k_2$. We know that $\pi_2(\tau_0) \cap \pi_2(\tau_1) \neq \emptyset$, so since dyadic intervals are nesting, we have $\pi_2(\tau_0) \subseteq \pi_2(\tau_1)$. Similarly, $\pi_2(\tau'_0) \subseteq \pi_2(\tau'_1)$. Since we have specified the other projection for both pairs, we know $\tau_0 \subseteq \tau_1, \tau'_0 \subseteq \tau'_1$. Then either C2 fails for $\tau_0 \times \tau'_0$ or $\tau_1 = \tau_0, \tau'_1 = \tau'_0$. Either case is a contradiction. 

    Finally, we prove that for any $\tau$, there is $\lesssim 1$ pairs $\tau'$ with $\tau \sim \tau'$. We know if $\tau \sim \tau'$, then $\pi_1(\tau')$ can be one of the two intervals $I'_j$ satisfying $I'_j \sim \pi_1(\tau)$. For each $I'_j$, if $\tau \sim_0 \tau'$ we must have either $\pi_1(\tau') = I'_k$ for one of two choices of $I'_k$ with $I'_k \sim I_k$. If $\tau \sim \tau'$, we must have $\pi_1(\tau') = \pi_1(\tau)$. In either case, we have at most two options for $I'_k$. In total, we have at most four options for $\tau'$. 
    
    \vspace{4pt}

    \noindent\textbf{The off-diagonal case (iii)}

    \vspace{4pt}

    Fix a pair $\tau, \tau' \in \mathcal{T}_{j,k}$ with $\tau \sim \tau'$ and a pair $(\xi_0, \eta_0) \in \tau \times \tau'$. Our definition of $\pi_1, \pi_2$ ensures that $\xi \mapsto (\pi_1(\xi), \pi_2(\xi))$ is bi-Lipschitz with constant $\sim 1$, so \[|\xi_0 - \eta_0| \sim |\pi_1(\xi_0) - \pi_1(\eta_0)| + |\pi_2(\xi_0) - \pi_2(\eta_0)| \sim 2^{-j} + 2^{-k} \sim 2^{-j}.\]It remains to prove that \begin{equation}\label{E:whitney props goal}\sup_{t \in [-1, 2]} |N_2(\xi_0 + t(\eta_0 - \xi_0)) \wedge (\eta_0 - \xi_0)| \gtrsim 2^{-k}.\end{equation}
    If \begin{equation}\label{E:N_2 close to eta_0 - xi_0}|2^j(\eta_0 - \xi_0) \wedge N_2(\xi_0)| \gtrsim 2^{j-k},\end{equation}then (\ref{E:whitney props goal}) is immediate, so suppose otherwise. It suffices then to prove that \begin{equation}\label{E:whitney props goal 1}\sup_{t \in [-1, 2]} |N_2(\xi_0 + t(\eta_0 - \xi_0)) \wedge 2^{-j}N_2(\xi_0)| \sim 2^{-k}.\end{equation}We also use (\ref{E:N_2 close to eta_0 - xi_0}) to conclude that the line segment $\{\xi_0 + t(\eta_0 - \xi_0): t \in [-1, 2]\}$ stays within $\sim 2^{-k}$ of the line segment $\{\xi_0 + 2^{-j+C}tN_2(\xi_0): t \in [-1, 2]\}$. Since we know $k \le k^2(\xi, j-C) + 2$ from (i) of this lemma, we see that \begin{equation}\label{E:whitney props goal 2}\sup_{t \in [-1, 2]} |N_2(\xi_0 + 2^{-j+C}tN_2(\xi_0)) \wedge N_2(\xi_0)| \sim 2^{-k}.\end{equation}The map $\xi \mapsto N_2(\xi) \wedge N_2(\xi_0)$ is Lipschitz with controlled constant, so (\ref{E:whitney props goal 2}) implies (\ref{E:whitney props goal 1}). 
  \end{proof}

  Now, we prove the spatial orthogonality we need for sum-sets of related pairs in $\mathcal{T}_{j,k}$.

  \begin{lemma}\label{L:almost orthog}
    For any $j, k \in \N$, $\Theta = \{\tau + \tau' \subset \R^2: \tau,\tau' \in \mathcal{T}_{j,k}, \tau \sim \tau'\}$ is boundedly overlapping, in the sense that for any $\theta \in \Theta$, \begin{equation}\label{E:almost orthog goal}|\{\theta' \in \Theta: \theta \cap \theta'\neq\emptyset\}| \lesssim 1.\end{equation}
  \end{lemma}

  The proof uses the fact that the projection-based rectangles which appear in $\mathcal{T}_{j,k}$ are essentially equivalent to the Euclidean rectangles over the form $R_2(\xi, k,j)$ and Euclidean rectangles work will with sum-sets. 

  \begin{proof}[Proof of Lemma \ref{L:almost orthog}]
    There are two cases to consider: if $|j - k| \le C$ or if $|j -k| > C$. First, we treat the $|j - k| \le C$ case. Without loss of generality, assume $j \le k$. For each $\tau \in \mathcal{T}_{j,k}$, we claim that there exists a $\sim 2^{-j}$-square $Q(\tau)$ containing $\tau$. Choose a point $\xi \in \tau$ and let $\xi'$ be the farthest point from $\xi$. The line segment $\ell$ joining $\xi$ and $\xi'$ must make angle $\gtrsim 1$ with one of the $N_i$ by the assumption (\ref{E:N_j near e_j}). Without loss of generality, assume $v = \xi - \xi'$ is quantitatively transverse to $N_1$. Then $\pi_1(\xi) - \pi_1(\xi') = d\pi_1(\xi^*)v$ for some choice of $\xi^* \in \ell$. Since $v$ is quantitatively transverse to $N_1(\xi^*) = \text{ker}(d\pi_1(\xi^*))$, we have $|d\pi_1(\xi^*)v| \sim |v|$. We know $|\pi_1(\xi) - \pi_1(\xi')| \lesssim 2^{-j}$, so we conclude that $|\xi - \xi^*| \lesssim 2^{-j}$ as well. Since $k \ge j$, the same reasoning applies to $\pi_2$.

    If $\tau \sim \tau'$, then $\pi_1(\tau), \pi_1(\tau')$ and $\pi_2(\tau), \pi_2(\tau')$ are each contained in a common $\sim 2^{-j}$. By the same reasoning we used to find $Q(\tau)$, we know there is a larger square $\tilde{Q}(\tau, \tau')$ of side length $\sim 2^{-j}$ containing $Q(\tau)$ and $Q(\tau')$. Then $\tau + \tau' \subset Q(\tau) + Q(\tau') \subset 2\tilde{Q}(\tau, \tau')$, so it suffices to show that the choices of $\{\tilde{Q}(\tau,\tau') : \tau\sim \tau'\}$ are $\sim 1$-overlapping. For a fixed $\tilde{Q}(\tau_0,\tau'_0)$, all intersecting $\tilde{Q}(\tau,\tau')$ are contained in a common $\sim 2^{-j}$ ball $B$. The $\pi_j$ are Lipschitz with constant $\sim 1$, so we know $\pi_1(B)$ are both contained in intervals $I_*, I'_*$ of length $\lesssim 2^{-j}$. For each $\tau \in \mathcal{T}_{j,j}$, $\tau = \pi^{-1}_1(I_j) \cap \pi^{-1}_1(I_k)$ for some unique pair $I_j, I_k$ of length $2^{-j}, 2^{-k}$ and respectively contained in $I_*, I'_*$. Since the intervals are unique and $|j-k| \le C$ there can only be $\sim 2^{C} \lesssim 1$ many different possible choices of $I_j, I_k$. Hence $\Theta$ is boundedly overlapping when $|j - k| \le C$.

    Now, we treat the case that $|j -k|>C$. Without loss of generality, we assume $j + C< k$ and in particular that we are always in the off-diagonal setting. Let $\mathcal{L}_m$ be the set of dyadic fractions in $\R$ with denominator $2^{m}$. Define \[\Xi = \{\pi_1^{-1}(\{t\}) \cap \pi_2^{-1}(\{s\}) : t \in \mathcal{L}_j, s \in \mathcal{L}_k\}.\]Since the flows of the $N_i$ are everywhere quantitatively transverse, $\Xi$ consists of a collection of isolated points. Recall that \[R_3(\xi, k, j) = \{\eta \in 2Q : |\pi_1(\eta) - \pi_1(\xi)| \le 2^{-j}, |\pi_2(\eta) - \pi_2(\xi)| \le 2^{-k}\}.\]Set \begin{align*}
    \mathcal{R}_{S} &= \{R_3(\xi, k-C/3, j-C/3): \xi \in \Xi, k < k^2(\xi, j-C)\},\\\mathcal{R}_{M} &= \{R_2(\xi, k-2C/3, j-2C/3): \xi \in \Xi, k < k^2(\xi, j-C)\},\\ \mathcal{R}_L &= \{R_3(\xi, k-C, j-C):\xi \in \Xi, k < k^2(\xi, j-C)\}.
    \end{align*}We choose $C$ large enough relative to the constants in Lemma \ref{L:comparable rectangles} so that \begin{equation}\label{E:subsets}R_3(\xi, k-C/3, j-C/3) \subset R_2(\xi, k-2C/3, j-2C/3) \subset R_3(\xi, k-C, j-C).\end{equation}We first claim that elements of $\mathcal{R}_L$ are $\sim 1$ overlapping, and hence the same is true for $\mathcal{R}_M$. Fix some $R = R_3(\xi, k-C,j-C) \in \mathcal{R}_L$. Since $k - C < k < k^2(\xi, j-C)$, we can apply Lemma~\ref{L:engulfing1}(b) to see that each element of $\mathcal{R}_L[R] = \{R' \in \mathcal{R}_L : R' \cap R \neq \emptyset\}$ is contained in $R_3(\xi, k-2C, j-2C)$. So if $R' = R_3(\xi', k-C, j-C)$, then $\pi_1(\xi') \in \mathcal{L}_{j} \cap I_{j-2C}$ and $\pi_2(\xi') \in \mathcal{L}_{k} \cap I_{k-2C}$ for fixed intervals $I_{j-2C} \in \mathcal{I}_{j-2C}$ and $I_{k-2C} \in \mathcal{I}_{k-2C}$. By construction, $(\pi_1(\xi'), \pi_2(\xi'))$ differs for each $\xi' \in \Xi$, while $|\mathcal{L}_j \cap I_{j-2C}|, |\mathcal{L}_k \cap I_{k-2C}| \lesssim 1$, so we conclude \[|\mathcal{R}_L[R]| \lesssim 1.\]We use (\ref{E:subsets}) to see that $\mathcal{R}_M$ is boundedly overlapping as well. 
    
    We claim that for any pair $\tau, \tau' \in \mathcal{T}_{j,k}$ with $\tau \sim \tau'$, we have some $R_M \in \mathcal{R}_M$ such that $\tau + \tau' \subset R_M$. By Lemma \ref{L:whitney props} (i), we see that $\tau = R_3(\xi_0, k', j)$ for some $k' \le k^2(\xi, j- C)$ and $\xi_0 \in \tau \cap \Xi$. It follows that $R_S = R_3(\xi_0, k-C/3, j-C/3) \in \mathcal{R}_S$. The dimensions of $R_S$ are wide enough relative to those of $\tau, \tau'$ to ensure that $\tau \cup \tau' \subset R_S$. By taking $C$ sufficiently large, we can use Lemma \ref{L:comparable rectangles} to conclude that \[\tau \cup \tau' \subset R_3(\xi, k-C/3, j-C/3) \subset R_2(\xi, k -2C/3+2, j - 2C/3+2).\]Since $R_2(\xi, j - 2C/3+2, k -2C/3+2)$ is (essentially) a Euclidean rectangle by construction, \begin{align}
        \tau + \tau' &\subset R_2(\xi + 2, k-2C/3 + 2, j - 2C/3) + R_2(\xi + 2, k-2C/3 + 2, j - 2C/3)\notag \\ &\subset \{2\xi : \xi \in R_2(\xi_0, k-2C/3, j -2C/3).\}\label{E:sum set control using rectangles}
    \end{align}We set $R_M(\tau + \tau') := R_2(\xi_0, k-2C/3, j -2C/3)$ and denote final set in (\ref{E:sum set control using rectangles}) by $2R_M(\tau + \tau')$. We will shortly use that $2R_M \cap 2R'_M \neq \emptyset$ if and only if $R_M \cap R'_M = \emptyset$. 

    Recall the definition $\Theta = \{\tau + \tau' : \tau, \tau' \in \mathcal{T}_{j,k}, \tau \sim \tau'\}$. Fix a choice of $\theta = \tau + \tau' \in \Theta$. We can cover \begin{equation}\label{E:cover almost orthog}\{\theta' \in \Theta : \theta' \cap \theta \neq \emptyset \} \subset \bigcup_{\substack{R \in \mathcal{R}_M \\ R \cap R_M(\theta) \neq \emptyset}} \{\theta' \in \Theta: R_M(\theta') = R\}.\end{equation}Since $\mathcal{R}_M$ is boundedly overlapping, the union in (\ref{E:cover almost orthog}) is always over $\lesssim 1$ elements, so (\ref{E:almost orthog goal}) follows from proving \[|\{\theta' \in \Theta: R_M(\theta') = R\}| \lesssim 1\]for any fixed $R \in \mathcal{R}_M$. This would follow from proving \begin{equation}\label{E:almost done with almost orthog}|\{\tau \in \mathcal{T}_{j,k}: \tau \subset R_L\}| \lesssim 1\end{equation}for any fixed $R_L \in \mathcal{R}_L$. For $\tau \subset R_L$, $I_j = \pi_1(\tau) \subset I_{j-C}$, $I_k = \pi_2(\tau) \subset I_{k-C}$ \footnote{Here we mildly violate our conventions: $I_k \in \mathcal{I}_k \cup \mathcal{I}_{k-2}$}. There are $\lesssim 1$ many distinct pairs $(I_j, I_k)$ or $(I_j, I_{k-2})$ with these properties, each corresponding to at most one choice of $\tau \subset R_L$, and hence (\ref{E:almost done with almost orthog}) is true, completing our proof of (\ref{E:almost orthog goal}).
  \end{proof}

  We use Lemma \ref{L:almost orthog} to establish the following corollary.

  \begin{corollary}\label{C:TVV ao}
    If $s < 2$, then for any function $\Omega \subset Q$ and for any $(j,k) \in \N^2$, we have

    \begin{equation}\label{E:TVV ao}\left|\left|\sum_{\substack{\tau, \tau' \in \mathcal{T}_{j,k} \\ \tau \sim \tau'}} E1_{\tau \cap \Omega'}E1_{\tau' \cap \Omega'}\right|\right|^s_s \lesssim \sum_{\tau, \tau' \in \mathcal{T}_{j,k}} ||E1_{\tau \cap \Omega'}E1_{\tau' \cap \Omega'}||_s^s \end{equation}
  \end{corollary}

  The proof of this corollary given (\ref{E:almost orthog goal}) is standard. We refer the reader to \cite{TVV98a} for details. Finally, we prove that we can apply bilinear extension for pairs $\tau \sim \tau'$.

  \begin{lemma}\label{L:better bil ext}
    Suppose $\tau, \tau' \in \mathcal{T}_{j,k}$ with $\tau \sim \tau'$. Then for any $f_{\tau}, f_{\tau'}$ respectively supported on $\tau, \tau'$, and any $r$ in a small neighborhood of $s'$, \begin{equation}\label{E:better bil ext}
      ||E1_{\tau \cap \Omega'}E1_{\tau' \cap \Omega'}||_s \lesssim  2^{(-j-k)\left(\frac2r + \frac2s - 2\right)}|\tau \cap \Omega'|^{1/r}|\tau' \cap \Omega'|^{1/r}.
    \end{equation}
  \end{lemma}

  \begin{proof}
    Without loss of generality, we assume that $k > j$. Recall we defined $\mathcal{L}_{j}$ to be the collection of dyadic rationals with denominator $2^{-j}$. Set \[\Xi = \{\pi_1^{-1}(\{t\}) \cap \pi_2^{-1}(\{s\}) : t \in \mathcal{L}_{j+C+1}, s \in \mathcal{L}_{k+C+1}\}\]so that $\mathcal{R}[\tau] = \{R_3(\xi, k-C, j-C): \xi \in \Xi \cap \tau\}$ covers $\tau$. Similarly, $\mathcal{R}[\tau']$ covers $\tau'$. Both sets have cardinality $\lesssim 1$. By Lemma \ref{L:whitney props}(i), each point $\xi\in \tau$ satisfy $k^2(\xi, j) \ge k^2(\xi, j-C) \ge k$ and by Lemma \ref{L:whitney props}(ii), each pair $\xi \in \tau, \eta \in \tau'$ are $(2^{-k}, 2^{-j})$ separated, for each $R \in \mathcal{R}[\tau]$ and $R' \in \mathcal{R}[\tau']$, $R, R'$ satisfy the assumption of Corollary \ref{C:rescaled bilinear rest}, so we conclude that \begin{equation}\label{E:bil ext for rects}||E1_{R \cap \Omega'}E1_{R' \cap \Omega'}||_s \lesssim  2^{(-j-k)\left(\frac2r + \frac2s - 2\right)}|R \cap \Omega'|^{1/r}|R' \cap \Omega'|^{1/r} \le 2^{(-j-k)\left(\frac2r + \frac2s - 2\right)}|\tau \cap \Omega'|^{1/r}|\tau' \cap \Omega'|^{1/r}.\end{equation}Summing (\ref{E:bil ext for rects}) over the $\lesssim 1$ pairs in $\mathcal{R}[\tau] \times \mathcal{R}[\tau']$, we see (\ref{E:better bil ext}).
  \end{proof}
  
  \subsection{The proof of Proposition \ref{P:RR}}

  We prove Proposition \ref{P:RR} following an argument of Vargas. The argument is similar to well-known bilinear-to-linear reductions for elliptic surfaces. Unfortunately, controlling the sum over the two separation parameters $j, k$ adds complexity compared to the single separation parameter used for elliptic surfaces. We use our assumption on the fiber lengths to gain extra control on the size of $\tau \cap \Omega$. With some book-keeping, this allows us to control the sum over both scales.

\begin{proof}
  It suffices to prove that \begin{equation}\label{E:rr real goal}||E1_{\Omega'}E1_{\Omega'}||_{L^s} \lesssim |\Omega(K)|^{2/s'}.\end{equation}We use Lemma \ref{L:whitney props}(ii) to see that \begin{equation}\label{E:PRR step 1}||E1_{\Omega'}E1_{\Omega'}||_{L^s} = \left|\left|\sum_{j,k \in \N^2} \sum_{\substack{\tau, \tau' \in \mathcal{T}_{j,k} \\ \tau \sim \tau'}} E1_{\tau \cap \Omega'}E1_{\tau' \cap \Omega'}\right|\right|_{L^{s}} \le \sum_{j,k \in \N^2} \left|\left|\sum_{\substack{\tau, \tau' \in \mathcal{T}_{j,k} \\ \tau \sim \tau'}} E1_{\tau \cap \Omega'}E1_{\tau' \cap \Omega'}\right|\right|_{L^s}.\end{equation}We bound each summands in the latter sum individually. Using Corollary \ref{C:TVV ao}, we have  \begin{equation}\label{E:PRR one scale}\left|\left|\sum_{\substack{\tau, \tau' \in \mathcal{T}_{j,k} \\ \tau \sim \tau'}} E1_{\tau \cap \Omega'}E1_{\tau' \cap \Omega'}\right|\right|^s_{L^s}\lesssim \sum_{\substack{\tau, \tau' \in \mathcal{T}_{j,k} \\ \tau \sim \tau'}} ||E1_{\tau \cap \Omega'}E1_{\tau' \cap \Omega'}||^s_{L^s}\end{equation}Using Lemma \ref{L:better bil ext} on each pair $\tau \sim \tau'$, we see that \begin{equation}\label{E:PRR one scale + bil ext}||\sum_{\substack{\tau, \tau' \in \mathcal{T}_{j,k} \\ \tau \sim \tau'}} E1_{\tau \cap \Omega'}E1_{\tau' \cap \Omega'}||^s_s \lesssim 2^{(j+k)(2 + \frac{2s}{r} -2s)}\sum_{\substack{\tau, \tau' \in \mathcal{T}_{j,k} \\ \tau \sim \tau'}} |\tau \cap \Omega(K)|^{s/r}|\tau' \cap \Omega(K)|^{s/r}.\end{equation}By Lemma \ref{L:whitney props}(ii), each $\tau \in \mathcal{T}_{j,k}$ has only $\sim 1$ many $\tau' \in \mathcal{T}_{j,k}$ with $\tau \sim \tau'$, so we can bound \[\sum_{\substack{\tau, \tau' \in \mathcal{T}_{j,k} \\ \tau \sim \tau'}} |\tau \cap \Omega(K)|^{s/r}|\tau' \cap \Omega(K)|^{s/r} \lesssim \sum_{\tau \in \mathcal{T}_{j,k}} |\tau \cap \Omega(K)|^{2s/r}.\]Therefore, \begin{equation}\label{E:rr almost there}
  ||\sum_{\substack{\tau, \tau' \in \mathcal{T}_{j,k} \\ \tau \sim \tau'}} E1_{\tau \cap \Omega'}E1_{\tau' \cap \Omega'}||^s_s \lesssim 2^{(j+k)(2 + \frac{2s}{r} -2s)} \sum_{\tau \in \mathcal{T}_{j,k}} |\tau \cap \Omega(K)|^{2s/r}
  \end{equation}

  Choose $J \in \N$ so that $|\pi_1(\Omega(K))| \sim 2^{-J}$. We claim that \begin{equation}\label{E:RR mass bd}
    |\tau \cap \Omega(K)| \lesssim \min\{2^{-j-k}, 2^{-J-k}, 2^{-j-K}, 2^{-J-K}\}.
  \end{equation}
  We know by construction that $|\pi_1(\tau)| \lesssim 2^{-j}$ and $|\pi_1(\Omega(K))| \sim 2^{-J}$, so $|\pi_1(\tau \cap \Omega(K))| \le \min\{2^{-j}, 2^{-J}\}$. We also know $|\pi_2(\tau)| \lesssim 2^{-k}$. Since $\pi_1$ and $\pi_2$ have quantitatively transverse kernels and their level sets intersect at a point, a change of variables gives \[|\tau \cap \Omega(K)| \lesssim |\pi_1(\tau)||\pi_2(\tau)| \lesssim 2^{-k}\min\{2^{-j}, 2^{-J}\},\]which gives the first two upper bounds in (\ref{E:RR mass bd}). For the second two upper bounds, we use the co-area formula to see that \[|\tau \cap \Omega(K)| \lesssim |\pi_1(\tau \cap \Omega(K))| \sup_{s \in \pi_1(\tau \cap \Omega(K))} \mathcal{H}^1(\pi_1^{-1}(s) \cap \left(\tau \cap \Omega(K)\right)).\]Since $e^{sN_1}(\xi)$ parametrizes $\pi_1^{-1}(\pi_1(\xi))$ and has unit speed, we conclude that \[\sup_{s \in \pi_1(\tau \cap \Omega(K))} \mathcal{H}^1(\pi_1^{-1}(s) \cap \left(\tau \cap \Omega(K)\right)) \lesssim \sup_{\xi \in \Omega(K)} \mathcal{H}^1(\{t : e^{tN_1}(\xi) \in \Omega(K)\}) \sim 2^{-K}.\]Hence \[|\tau \cap \Omega(K)|\lesssim |\pi_1(\tau \cap \Omega(K))| 2^{-K} \lesssim \min\{2^{-j}, 2^{-J}\} 2^{-K},\]giving the second pair of upper bounds in (\ref{E:RR mass bd}).

  Now, we further bound $\sum_{\tau \in \mathcal{T}_{j,k}} |\tau \cap \Omega(K)|^{2s/r}$ in different ways depending on the values $j$ and $k$. If $j \ge J, k \ge K$, we bound \begin{equation}\label{E:RR jk}\sum_{\tau \in \mathcal{T}_{j,k}}|\tau \cap \Omega(K)|^{2s/r} \lesssim 2^{(\frac{2s}{r}-1)(-j-k)}\sum_{\tau \in \mathcal{T}_{j,k}} |\tau \cap \Omega(K)| \lesssim 2^{(\frac{2s}{r}-1)(-j-k)}|\Omega(K)|.\end{equation}
  If $j \ge J, k < K$, we similarly bound \begin{equation}\label{E:RR jK}\sum_{\tau \in \mathcal{T}_{j,k}}|\tau \cap \Omega(K)|^{2s/r} \lesssim 2^{(\frac{2s}{r}-1)(-j-K)}\sum_{\tau \in \mathcal{T}_{j,k}} |\tau \cap \Omega(K)| \lesssim 2^{(\frac{2s}{r}-1)(-j-K)}|\Omega(K)|.\end{equation}
  Also, for $j < J, k \ge K$: \begin{equation}\label{E:RR Jk}\sum_{\tau \in \mathcal{T}_{j,k}}|\tau \cap \Omega(K)|^{2s/r} \lesssim 2^{(\frac{2s}{r}-1)(-J-k)}\sum_{\tau \in \mathcal{T}_{j,k}} |\tau \cap \Omega(K)| \lesssim 2^{(\frac{2s}{r}-1)(-J-k)}|\Omega(K)|.\end{equation}
  Finally, if $j < J, k < K$, we instead bound \begin{equation}\label{E:RR JK}\sum_{\tau \in \mathcal{T}_{j,k}}|\tau \cap \Omega(K)|^{2s/r} \lesssim \sum_{\tau \in \mathcal{T}_{j,k}} |\tau \cap \Omega(K)||\Omega(K)|^{2s/r-1} \lesssim |\Omega(K)|^{2s/r}.\end{equation}Combining each of these with (\ref{E:PRR step 1}), (\ref{E:PRR one scale + bil ext}) gives 
  \begin{align}
    ||E1_{\Omega'}E1_{\Omega'}||_{L^s} &\le \sum_{j \ge J, k \ge K} 2^{(j+k)(\frac{2}{s} + \frac{2}{r} -2)} 2^{(\frac{2}{r}-\frac1s)(-j-k)}|\Omega(K)|^{1/s} \label{E:rr jk}\\ &+ \sum_{j \ge J, k < K}2^{(j+k)(\frac{2}{s} + \frac{2}{r} -2)} 2^{(\frac{2}{r}-\frac1s)(-j-K)}|\Omega(K)|^{1/s} \label{E:rr jK}\\ &+\sum_{j < J, k \ge K} 2^{(j+k)(\frac{2}{s} + \frac{2}{r} -2)}2^{(\frac{2}{r}-\frac1s)(-J-k)}|\Omega(K)|^{1/s} \label{E:rr Jk}\\
    &+ \sum_{j < J, k < K} 2^{(j+k)(\frac{2}{s} + \frac{2}{r} -2)}|\Omega(K)|^{2/r}\label{E:rr JK}
  \end{align} 
  To conclude (\ref{E:rr real goal}) and hence (\ref{E:restricted restricted linear}), it remains to prove \begin{equation}\label{E:last step}\max\{\text{(\ref{E:rr jk})}, \text{(\ref{E:rr Jk})}, \text{(\ref{E:rr jK})}, \text{(\ref{E:rr JK})}\} \lesssim |\Omega|^{2/s'}.\end{equation}We know $s > \frac{3}{2}$, so $\frac{3}{s} -2  < 0$, which implies that \[(\ref{E:rr jk}) = \sum_{j \ge J, k \ge K} 2^{(j+k)(\frac{3}{s}-2)}|\Omega(K)|^{1/s} \lesssim |\Omega(K)|^{1/s}2^{(J+K)(\frac{3}{s} - 2)} \sim |\Omega(K)|^{2/s'}\]
  We can choose $r > s'$ to have $\frac{2}{s} + \frac{2}{r} - 2 > 0$. Then \begin{align*}(\ref{E:rr jK}) &= 2^{-K(\frac{2}{r} - \frac1s)}\sum_{j \ge J, k < K} 2^{k(\frac{2}{s} + \frac2r - 2)} 2^{j(\frac{3}{s}-2)}|\Omega(K)|^{1/s}\\ &\lesssim |\Omega(K)|^{1/s}2^{(J+K)(\frac{3}{s} - 2)} \sim |\Omega(K)|^{2/s'}\end{align*}By similar reasoning as for (\ref{E:rr jK}), we also have \[(\ref{E:rr Jk}) \lesssim |\Omega(K)|^{2/s'} \text{ and }(\ref{E:rr JK}) \lesssim |\Omega(K)|^{2/s'}.\]Hence, (\ref{E:last step}) holds, completing the proof. 
\end{proof}

  \section{Linear restriction for general sets}

  We now derive the full linear restriction estimate. 

  \begin{proposition}\label{P:Bil2RST}
    Suppose for some $s > \frac32$ Corollary \ref{C:rescaled bilinear rest} holds.

  Then for any $\Omega \subset [-1,1]^2$, we have the restricted strong-type bound \begin{equation}\label{E:RST goal}
    ||E1_{\Omega}||_{L^{2s}} \lesssim |\Omega|^{1/s'}
  \end{equation}
  \end{proposition}

  Since Corollary \ref{C:rescaled bilinear rest} holds for $s > \frac{11}{7}$, (\ref{E:RST goal}) holds with the same range of $s$. By interpolating the restricted strong-type bound (\ref{E:RST goal}) with the classical $L^2 \to L^4$ estimate for $\mathcal{E}$, we arrive at strong-type scaling-line estimates $L^{s'} \to L^{2s}$ for $s > \frac{11}{7}$. Labeling $2s = p$, we have Theorem \ref{T:1}. Theorem \ref{T:2} follows similarly: if (\ref{E:T2 assumption}) is true for any $p > p_0, 2q' < p$, then Corollary \ref{C:rescaled bilinear rest} holds for $s > \frac{p_0}{2}$. We can again interpolate to find strong-type $L^{s'} \to L^{2s}$ bounds for any $s > \frac{p_0}{2}$, and hence (\ref{E:T2 goal}) is true. 
  
  Our proof of Proposition \ref{P:Bil2RST} follows similar lines to the second named author's proof of scaling-line restriction estimates for the hyperbolic paraboloid. We begin by proving those sets which nearly extremize the inequality in Proposition \ref{P:RR} are unions of small numbers of tiles from $\mathcal{T}_{j,k}^0$.

  \subsection{Near-extremizers for Proposition \ref{P:RR}}

  We specifically aim to prove the following decomposition for sets of constant $N_1$-fiber length. Of course, \emph{mutatis mutandis} the same result holds for sets of $N_2$-fiber length, but going forward we only state the $N_1$-fiber length result.
  \begin{proposition}\label{P:S22}
    Under the assumptions of Proposition \ref{P:Bil2RST}, the following holds. 
    
    Suppose $\Omega$ has constant $N_1$-fiber length. Take $\varepsilon$ dyadic so that for any measurable $\Omega' \subset \Omega$, \begin{equation}\label{E:eps def}||E1_{\Omega'}||_{L^{2s}} \lesssim \varepsilon^2|\Omega|^{1/s'}.\end{equation}If $\Omega = \Omega(K)$ and $|\pi_1(\Omega)|=2^{-J}$, then there exists $\{\Omega_{\delta} : \delta \in (0, \varepsilon] \cap 2^{-\N}\}$ so that the following hold:
    \begin{enumerate}[label = (\arabic*)]
        \item $\Omega = \bigsqcup_{0 < \delta \le \varepsilon} \Omega_{\delta}$, up to a set of measure zero.
        \item For each $\delta \in (0, \varepsilon] \cap 2^{-\N}$, $\Omega_{\delta} \subset \bigcup_{\tau \in \scriptR_{\delta}} \tau$, where $\scriptR_{\delta} \subset \scriptT^0_{J,K}$ has cardinality $\delta^{-O(1)}$. 
        \item For each $\Omega' \subseteq \Omega_{\delta}$, \begin{equation}\label{E:22goal}||E1_{\Omega'}||_{L^{2s}} \lesssim \delta|\Omega|^{1/s'}.\end{equation}
    \end{enumerate}
  \end{proposition}

  \begin{proof}
    Let $C$ be a fixed large constant, whose value will be determined over the course of the proof. Let $S = \pi_1(\Omega)$. For any $\eta \le \varepsilon$ dyadic and $t \in S$, let $I_{\eta}(t)$ be the maximal interval $I$ containing $t$ for which $|I \cap S| \ge \eta^C |I|$. Bounding $|I_{\eta}(t)| \le \eta^{-C} |I_{\eta}(t) \cap S| \le \eta^{-C} |S| \lesssim \eta^{-C} 2^{-J}$, we ensure that \begin{equation}\label{E:Ieta bound}
        |I_{\eta}(t)| \lesssim \eta^{-C} 2^{-J}
    \end{equation}
    For any Lebesgue point $t \in S$, $I_{\eta}(t)$ exists for all $\eta$. Define $T_{\eta}$ to be those points where (\ref{E:Ieta bound}) is within $\eta^{2C}$ of being sharp: \[T_{\eta} = \{t \in S: |I_{\eta}(t)| \ge \eta^{C}2^{-J}\}.\]Finally, set \[S_{\eta} = \{t : \eta \text{ is the largest value for which }t \in T_{\eta}\}.\]
    
    The $S_{\eta}$ are disjoint and by construction cover all of $S$ outside a measure zero set. For $\eta < \varepsilon$, we have $S_{\eta} = T_{\eta} \setminus T_{2\eta}$. Define $\Omega_{\eta}^1 = \Omega \cap \pi_1^{-1}(S_{\eta})$, these are also disjoint sets and cover all but a measure zero subset of $\Omega$. 

    \begin{lemma}\label{L:S23}
        For each $0 < \eta \le \varepsilon$
        \begin{enumerate}[label = (\roman*)]
            \item $S_{\eta}$ is contained in a union of $\lesssim \eta^{-3C}$ many dyadic intervals $\mathcal{I}_{\eta} \subset \mathcal{I}_{J}$. 
            \item For each $\Omega' \subset \Omega_{\eta}^1$, \begin{equation}\label{E:intermediate RRB}
                ||E1_{\Omega'}||_{L^{2s}} \lesssim \eta^2 |\Omega|^{1/s'}.
            \end{equation}
        \end{enumerate}
    \end{lemma}
    \begin{remark}\label{R:S23}
      This lemma depends on the construction of the $\Omega^1_{\eta}$ and some other facts about $\Omega$ and the extension operator. In particular, we know 
      \begin{itemize}
        \item (\ref{E:eps def}) holds to ensure that we can take $\eta \le \varepsilon$;
        \item $|\pi_1(\Omega)| \lesssim 2^{-J}$ to give the right dimensions on the strips; and 
        \item $\Omega$ has constant $\pi_1$ fiber length.
      \end{itemize}

      We will later apply the construction before the lemma to decompose sets $\Omega^2_{\eta,\rho}$ into sets $\Omega^3_{\eta, \rho, \delta}$ and then and then apply the lemma with $\Omega^3_{\eta, \rho, \delta}$ taking the place of $\Omega^1_{\eta}$; $\pi_2$ taking the place of $\pi_1$; $\pi_2$-strips $\tilde{S}_{\delta}$ taking the place of $S_{\eta}$; and a parameter $\rho$ taking the place of $\varepsilon$ by satisfying \begin{equation}\label{E:intermediate intermediate RBB}
        ||E1_{\Omega'}||_{L^{2s}} \lesssim \rho^2 |\Omega|^{1/s'}.
      \end{equation}for any $\Omega' \subset \Omega^2_{\eta, \rho}$. Before doing so, we will confirm 

      \begin{itemize}
        \item (\ref{E:intermediate intermediate RBB}) so we can take $\delta \le \rho$;
        \item $|\pi_2(\Omega^2_{\eta, \rho})| \lesssim \rho^{-C} 2^{-K}$, so the strips have width $\rho^{-C} 2^{-K}$; and 
        \item $\Omega^2_{\eta,\rho}$ has constant $\pi_2$ fiber length.
      \end{itemize}

      In this setting, Lemma \ref{L:S23} tells using that for any $0 < \delta \le \rho$,

      \begin{enumerate}[label = (\roman{*})]
        \item $\tilde{S}_{\eta}$ is contained in a union of $\lesssim \rho^{-C}\delta^{-3C}$ many dyadic intervals $\tilde{\mathcal{I}}_{\delta} \subset \mathcal{I}_K$. 
        \item For each $\Omega' \subset \Omega^3_{\eta, \rho, \delta}$, \begin{equation}\label{E:intermediate RRB2}
          ||E1_{\Omega'}||_{L^{2s}} \lesssim \delta^2|\Omega|^{1/s'}
        \end{equation}
      \end{enumerate}
      We conclude (i) from conclusion (i) of Lemma \ref{L:S23} by further sub-dividing the $\delta^{-3C}$ many strips of width $\rho^{-C}2^{-K}$ into $\rho^{-C}\delta^{-3C}$ many strips of width $2^{-K}$.
    \end{remark}

    \begin{proof}
        Proof of (i): For each $t \in S_{\eta}$, we have some interval $I_{\eta}(t)$. These intervals are maximal in length and hence disjoint. For each $I_{\eta}(t)$, we know that $|I_{\eta}(t)| \ge \eta^C2^{-J}$ since $t \in T_{\eta}$. We also have $|I_{\eta}(t) \cap S| \ge \eta^C |I_{\eta}(t)|$ from how we defined $I$. Since $|S| \lesssim 2^{-J}$, there can only $\lesssim \eta^{-2C}$ many $I_{\eta}(t)$ covering $S_{\eta}$. Each has length at most $\eta^{-C}2^{-J}$ from (\ref{E:Ieta bound}), so each of the $\eta^{-2C}$ many intervals can be covered by $\sim \eta^{-C}$ many elements of $\mathcal{I}_J$. Together, we can cover $S_{\eta}$ with $\lesssim \eta^{-3C}$ many dyadic intervals in $\mathcal{I}_J$. We take $\mathcal{I}_{\eta}$ to be this set of intervals.

        Proof of (ii): If $\eta = \varepsilon$, this follows from (\ref{E:eps def}), so assume $\eta < \varepsilon$. In particular, $S_{\eta} = T_{\eta} \setminus T_{2\eta}$. We prove (\ref{E:intermediate RRB}) by sharpening our proof of Proposition \ref{P:RR} with $\Omega_{\eta}^1$ taking the place of $\Omega(K)$, which is acceptable since $\Omega_{\eta}^1 = \Omega_{\eta}^1(K)$ by construction. The only difference is a sharpening of (\ref{E:RR mass bd}): if $\tau \in \mathcal{T}_{j,k}$ and $|j - J| \le \frac{C}{4}\log(1/\eta)$, then \begin{equation}\label{E:RR mass bd2} |\tau \cap \Omega_{\eta}^1| \lesssim \eta^{3C/4}\min\{2^{-j-k}, 2^{-J-k}, 2^{-j-K}, 2^{-J-K}\}. \end{equation}To see this, write $I_j = \pi_1(\tau)$, a dyadic interval of length $2^{-j}$ and assume $|j - J| \le \frac{C}{4}\log(1/\eta)$. If \begin{equation}\label{E:iRBB1}|I_j \cap S| > (2\eta)^C2^{-j} = (2\eta)^C |I|,\end{equation}any point $t \in I_j \cap S$ must have $|I_{2\eta}(t)| \ge |I_j| \ge 2^{-j} \gtrsim (2\eta)^C 2^{-J}$, and so $t \in T_{2\eta}$. We conclude that $I_j \cap S_{\eta} = \emptyset$. Since $\Omega_{\eta}^1 \cap \tau \subset \pi_1^{-1}(I_j \cap S_{\eta})$, the same is true for $\tau \cap \Omega_{\eta}^1$, so (\ref{E:RR mass bd2}) is immediate. If (\ref{E:iRBB1}) is false, we use that $2^{-j} \le \eta^{-C/4}2^{-J}$ to conclude \[|I_j \cap S_{\eta}| \le |I_j \cap S| \le (2\eta)^C2^{-j} \lesssim \eta^{\frac{3C}{4}}\min\{2^{-j}, 2^{-J}\}.\]
        Using the coarea formula as in the proof of Proposition \ref{P:RR}, we bound \[|\Omega_{\eta}^1 \cap \tau| \lesssim |I_j \cap S|\min\{2^{-k}, 2^{-K}\} \lesssim \eta^{\frac{3C}{4}}\min\{2^{-j}, 2^{-J}\}\min\{2^{-k}, 2^{-K}\},\]giving (\ref{E:RR mass bd2}). Summing geometric series as in Proposition \ref{P:RR}, we see that \[||E1_{\Omega'}||_{L^{2s}} \lesssim \eta^{\Theta_{r,s}(C)}|\Omega_{\eta}^1|^{1/s'},\]Choosing $C$ sufficiently large, we have $\Theta_{r,s}(C) \ge 2$, giving (\ref{E:intermediate RRB}).
    \end{proof}
    So far, we have decomposed $\Omega$ into sets $\Omega_1^{\eta}$, which satisfy conclusions (1) and (3) of Proposition \ref{P:S22} with $\eta$ taking the place of $\delta$ and a weaker version of conclusion (2): instead of being covered by $2^{-J} \times 2^{-K}$ tiles, the sets are only covered by $2^{-J}$-strips. We would like to carry out the same decomposition with the projection $\pi_2$ taking the place of the projection $\pi_1$, but before doing so, we need to replicate our control on the lengths of the $\pi_1$ fibers for the $\pi_2$ fibers, while sharpening (\ref{E:eps def}). 

    For dyadic choices of $\rho$, $0 < \rho < \eta$, define \[U_{\rho} = \{\xi_2 \in \pi_2(\Omega^1_{\eta}) : \mathcal{H}^1(\pi_2^{-1}(\xi_2) \cap \Omega^1_{\eta}) \sim \rho^C 2^{-J}\}.\]Set \[U_{\eta} = \{\xi_2 \in \pi_2(\Omega^1_{\eta}) : \mathcal{H}^1(\pi_2^{-1}(\xi_2) \cap \Omega^1_{\eta}) \gtrsim \eta^C 2^{-J}\}.\]Define $\Omega^2_{\eta, \rho} = \pi_2^{-1}(U_{\rho}) \cap \Omega_{\eta}^1$. Choosing implicit constants appropriately, we ensure that $\bigsqcup_{\eta \ge \rho} \Omega^2_{\eta, \rho} = \Omega^1_{\eta}.$ These sets have $\pi_2$ fiber length roughly constant. We now prove a sharpened version of (\ref{E:eps def}) for $\Omega_{\eta, \rho}$, which will allow us to apply Lemma \ref{L:S23} with $\pi_2$ taking the place of $\pi_1$ and $\rho$ taking the place of $\varepsilon$. 

    \begin{lemma}\label{L:S24}
      For any $\Omega' \subset \Omega^2_{\eta, \rho}$, \begin{equation}\label{E:rho def}
           ||E1_{\Omega'}||_{L^{2s}} \lesssim \rho^2 |\Omega|^{1/s'}.
      \end{equation}
    That is, (\ref{E:intermediate intermediate RBB}) is true.
    \end{lemma}

    \begin{proof}
      If $\eta = \rho$, (\ref{E:intermediate intermediate RBB}) is simply (\ref{E:22goal}), so assume $\eta < \rho$ and in particular, $\Omega^2_{\eta,\rho}$ has constant $\pi_2$ fiber length $\sim \rho^C 2^{-J}$.

      As in the proof of Lemma \ref{L:S23}, we sharpen the proof of Proposition \ref{P:RR} by improving (\ref{E:RR mass bd}). Indeed, we see that if $\tau \in \mathcal{T}_{j,k}$, then \begin{equation}\label{E:RR mass bd3}|\Omega' \cap \tau| \lesssim \min\{2^{-J-K}, 2^{-j - K}, 2^{-j-k}, \rho^C 2^{-J-k}\}.\end{equation}We need more care to prove (\ref{E:intermediate intermediate RBB}). In Proposition \ref{P:RR} and in Lemma \ref{L:S23}, we subdivide $\N^2$ into four regions with boundaries at the lines $j = J$ and $k = K$. Here, we still subdivide into four regions, but they are not defined by a pair of lines. Instead, we start with the unbounded rectangle:

      \[R_4 = \{(j,k) : j \le J + \frac{C|\log \rho|}{2}, k \ge K - \frac{C|\log \rho|}{2}\}.\]On this rectangle, we have $\min\{2^{-J-K}, 2^{-j - K}, 2^{-j-k}, \rho^C 2^{-J-k}\} = \rho^C 2^{-J-k}$. Next, we have \[R_1 = \{(j, k) : k \le K - \frac{C|\log \rho|}{2}, j \le J\}, R_3 = \{(j,k) : j \ge J + \frac{C|\log \rho|}{2}, k \ge K\}.\]Finally, we have \[R_2 = \{(j,k) : k \le K, j \ge J + \frac{C|\log \rho|}{2}\} \cup \{(j,k) : j \ge J, k \le K - \frac{C|\log \rho|}{2}\}.\]We follow the same argument as in Proposition \ref{P:RR} with our bound (\ref{E:RR mass bd3}) taking the place of (\ref{E:RR mass bd}) to conclude that \begin{align*}
        ||E1_{\Omega'}E1_{\Omega'}||_{L^s} \lesssim &\sum_{(j,k) \in R_1} 2^{(j+k)(\frac2s + \frac2r - 2)}2^{-2\frac{J+K}{r}}\\ &+ \sum_{(j,k) \in R_1} 2^{-j(2 - \frac3s)}2^{k(\frac2s + \frac2r - 2)}2^{-\frac{J}{s}}2^{-\frac{2K}{r}}\\
        &+ \sum_{(j,k) \in R_2} 2^{-j(2 - \frac3s)}2^{k(\frac2s + \frac2r - 2)}2^{-\frac{J}{s}}2^{-\frac{2K}{r}}\\
        &+ \sum_{(j,k) \in R_3} 2^{-(j+k)(2 - \frac3s)}2^{-\frac{J+K}{s}}\\ &+ \rho^{C(\frac2r - \frac1s)} \sum_{(j,k) \in R_4} 2^{-k(2-\frac3s)}2^{j(\frac2s + \frac2r -2)}2^{-\frac{2J}{r}-\frac{K}{s}}
      \end{align*}
      When we sum each of the geometric series as in the proof of Proposition \ref{P:RR}, we have gains of $\rho^{\Theta_{r,s}(C)}$ for $R_1, R_2$, and $R_3$ thanks to the smaller domain of summation and the gain from $\rho^{C(\frac2r - \frac1s)}$ in the $R_4$ sum dominates the loss from the larger domain of summation in the $R_4$. See \cite{S19} for details. In total we have \[||E1_{\Omega'}E1_{\Omega'}||_{L^s} \lesssim \rho^{\Theta_{r,s}(C)}|\Omega|^{2/s}.\]Taking $C$ sufficiently large gives (\ref{E:rho def}).
    \end{proof}

    When $\rho < \eta$, we have approximately constant $\pi_2$ fiber length, and when $\rho = \eta$ we can decompose $\Omega^2_{\eta, \rho}$ into $\sim |\log \eta|$ disjoint subsets with constant $\pi_2$ fiber length, which we continue to refer to as $\Omega^2_{\eta, \rho}$. So we have constant fiber length and in either case, we have (\ref{E:rho def}) for $\Omega^2_{\eta, \rho}$ from Lemma \ref{L:S24}. Now for each choice of $\Omega^2_{\eta, \rho}$, we carry out the same decomposition as described at the start of the proof for the $\pi_2$ projection and replacing $\eta$ with $\delta$. We arrive at disjoint sets $\Omega^3_{\eta, \rho, \delta}$, where $\delta$ is a dyadic number $\le \rho$. Finally, the fact that $\Omega^3_{\eta,\rho,\delta}$ has constant fiber length $\sim \rho^C 2^{-J}$ and measure $2^{-J-K}$ implies that $|\pi_2(\Omega^3_{\eta, \rho, \delta})| \lesssim \rho^{-C} 2^{-K}$. So we have the  conclusions (i) and (ii) discussed in Remark \ref{R:S23} for $\Omega^3_{\eta, \rho, \delta}$.

    We can cover each $\Omega^3_{\eta,\rho,\delta}$ with a set $\mathcal{T}_{\eta, \rho,\delta} \subset \mathcal{T}^0_{J,K}$ with cardinality $\lesssim (\eta\rho\delta)^{-3C}$ coming from pairwise intersections of elements of $\pi_1^{-1}(I) \cap \pi_2^{-1}(\tilde{I})$, with $I \in \mathcal{I}_{\eta}$ and $\tilde{I} \in \tilde{\mathcal{I}}_{\delta}$. Note that the collection of tiles here is $\mathcal{T}^0_{J,K}$, as opposed to $\mathcal{T}_{J,K}$. Define $\Omega_{\delta} = \bigcup_{\eta \ge \rho \ge \delta} \Omega^3_{\eta,\rho,\delta}$. By construction, $\bigsqcup_{\delta \le \varepsilon} \Omega_{\delta} = \Omega$ up to a set of measure zero, giving the conclusion (1) Proposition \ref{P:S22}. Since $\eta, \rho$ are dyadic and at least $\delta$, we can cover $\Omega^3_{\eta,\rho,\delta}$ with $\mathcal{R}_{\delta} = \bigcup_{\eta \ge \rho \ge \delta} \mathcal{T}_{\eta, \rho,\delta}$ which comprises $\lesssim |\log \delta|^2\delta^{-9C}$ tiles from $\mathcal{T}_{\eta, \rho, \delta}$. This gives conclusion (2) of the proposition.

    Finally, using the triangle inequality and summing (\ref{E:intermediate RRB2}) over $|\log \delta|^2$ many choices of $\eta \ge \rho \ge \delta$ gives \[||E1_{\Omega'}||_{L^{2s}} \lesssim |\log \delta|^2 \delta^2|\Omega|^{1/s'} \lesssim \delta|\Omega|^{1/s'}\]for any $\Omega' \subset \Omega_{\delta}$, yielding conclusion (3).
  \end{proof}

  \subsection{Linear restriction estimates for general sets $\Omega$.}

  We assume that we actually have fixed one choice of $\Omega$ in the following pair of definitions in order to suppress the dependence of various parameters on $\Omega$.
  \begin{definition}\label{D:Keps}
    \begin{itemize}
      \item For each $K \in \N$, define $J(K)$ so that $|\Omega(K)| \sim 2^{-J(K) - K}$. In other words, $2^{-J(K)} \sim |\pi_1(\Omega(K))|$. 
      \item For dyadic choices of $\varepsilon$, let \[\mathcal{K}(\varepsilon) = \left\{K : \sup_{\Omega' \subset \Omega(K)} \frac{||1E_{\Omega'}||_{L^{2s}}}{|\Omega(K)|^{1/s'}} \sim \varepsilon^2 \right\}.\]Choosing the implicit constant appropriately, we ensure that \[\{\mathcal{K}(\varepsilon) \subset \N: \varepsilon \in 2^{-\N}\}\]partitions $\mathbb{N}$. 
      \item Let $\mathcal{R}_{\delta}(K) \subset \mathcal{T}_{j,k}$ be the set of $\lesssim \delta^{-C}$ many tiles which Proposition \ref{P:RR} (ii) finds for $\Omega(K)$. 
    \end{itemize}
  \end{definition}
  We expand the sum 
  \begin{proposition}\label{P:S31}
    Denote by $\Omega_{\delta}(K)$ the partition of $\Omega(K)$ given by Proposition \ref{P:S22}. For $0 < \delta \le \varepsilon$, we have \begin{equation}\label{E:S31 goal}
      \left|\left| \sum_{K \in \mathcal{K}(\varepsilon)} E1_{\Omega_{\delta}(K)}\right|\right|_{L^{2s}}^{2s} \lesssim |\log \delta|^{O(1)}\left(\sum_{K \in \mathcal{K}(\varepsilon)} ||E1_{\Omega_{\delta}(K)}||_{L^{2s}}^{2s} + \delta |\Omega|^{\frac{2s}{s'}}\right).
    \end{equation}
  \end{proposition}  

  \begin{lemma}\label{L:S32}
    Suppose $K < K'$ and $J := J(K) > J' := J(K')$. Suppose $\tau \in \mathcal{R}_{\delta}(K)$ and $\tau' \in \mathcal{R}_{\delta}(K')$. Then \begin{equation}\label{E:S32}
      ||E1_{\tau \cap \Omega_{\delta}(K)}E1_{\tau' \cap \Omega_{\delta}(K')}||_{L^s} \lesssim 2^{-\Theta_s(|K - K'|)}\max\{|\Omega(K)|^{\frac2{s'}},|\Omega(K')|^{\frac2{s'}}\}
    \end{equation}
  \end{lemma}

  \begin{proof}
    Without loss of generality, we may assume \[\max\{|\Omega(K)|, |\Omega(K')|\} = |\Omega(K)|.\]We have defined $\tau, \tau'$ so that $|\tau| \sim |\Omega(K)|$ and $|\tau'| \sim |\Omega(K')|$. Proposition \ref{P:RR} is immediate if $|K' - K| \lesssim 1$ or $|(J + K) - (J'+K')| \gtrsim |K'-K|$, so we can assume $K' - K \ge C_K$ for some admissible choice of $C_K$. By using the reverse triangle inequality on $|(J + K) - (J'+K')| \ll |K' - K|$, we see that $|J - J'| \sim |K - K'|$ and hence $|J - J'| \ge C_J$ for some admissible choice of $C_J$. We set $\tilde{K} = K + \frac{K' - K}{5}, \tilde{J} = J' + \frac{J - J'}{5}$. 

    By construction, we have $\tau = \pi_1^{-1}(I_J) \cap \pi_2^{-1}(I_K)$ and $\tau' = \pi_1^{-1}(I_{J'}) \cap \pi_2^{-1}(I_{K'})$. The worst case in proving (\ref{E:S32}) is when $I_{K'} \subset I_K$ and $I_{J} \subset I_{J'}$, so we assume this is the case. We set $\tilde{I}_{k} = N_{2^{-k}}(I_{K'})$ and $V_k = I_{K} \cap (\tilde{I}_{k-1} \setminus \tilde{I}_{k})$ for $k < \tilde{K}$ and $V_{\tilde{K}} = I_{K} \cap \tilde{I}_{\tilde{K}}$. In this way, we decompose $I_K = \bigsqcup_{k=K}^{\tilde{K}} V_{k}$ so when $k < \tilde{K}$, $V_k$ is separated from $I_{K'}$ by $\sim 2^{-k}$. Similarly, we decompose $I_{J'} = \bigsqcup_{j=J'}^{\tilde{J}} H_j$, where for $j < \tilde{J}$, each $H_j$ is separated from $I_{J'}$ by $\sim 2^{-j}$. We decompose $\tau = \bigsqcup_{k=K}^{\tilde{K}} \tau_k$, where $\tau_k = \pi_1^{-1}(I_J) \cap \pi_2^{-1}(V_k)$. Similarly, we decompose $\tau' = \bigsqcup_{j = J'}^{\tilde{J}} \tau'_j$, where $\tau'_j = \pi_1^{-1}(H_j) \cap \pi_2^{-1}(I_{K'})$. 

    We decompose the left-hand side of (\ref{E:S31 goal}) as follows \begin{align}
      ||E1_{\tau \cap \Omega_{\delta}(K)}&E1_{\tau' \cap \Omega_{\delta}(K')}||_{L^s}\notag \\ &\lesssim ||E1_{{\tau_{\tilde{K}}} \cap \Omega_{\delta}(K)}E1_{\tau' \cap \Omega_{\delta}(K')}||_{L^s}   \label{E:S32 easy1} \\ &+ ||E1_{{\tau\setminus \tau_{\tilde{K}}} \cap \Omega_{\delta}(K)}E1_{\tau'_{\tilde{J}} \cap \Omega_{\delta}(K')}||_{L^s}   \label{E:S32 easy2}\\ &+ \sum_{j = J'}^{\tilde{J}}\sum_{k = K}^{\tilde{K}} ||E1_{\tau_k \cap \Omega_{\delta}(K)}E1_{\tau'_j \cap \Omega_{\delta}(K')}||_{L^s}\label{E:S32 hard}
    \end{align}
    We claim that any $N_2$-fiber intersecting $\tau_{\tilde{K}}$ has fiber length $\lesssim 2^{-J}$. Take one such fiber $\gamma = e^{tN_2}(\xi_0)$ and let $U = \{t \in \R : e^{tN_2}(\xi_0) \in \tau_{K'-C}\}$. We need to prove $|U| \lesssim 2^{-J}$. Since $N_2(\xi)$ is uniformly transverse to $\text{ker}(d\pi_1(\xi)) = N_1(\xi)\R$, we know that $|\frac{d}{dt}\pi_1(e^{tN_2}(\xi_0))| \gtrsim 1$. For $\pi_1(\{e^{tN_2}(\xi_0) : t \in U\})$ to stay in the interval $I_J$ of length $2^{-J}$, we must have $|U| \lesssim 2^{-J}$. We conclude that $\tau_{K'-C}$ has $N_2$-fiber length at most $2^{-J}$. By essentially the same reasoning see that $\tau'_{\tilde{J}}$ has $N_1$-fiber length at most $2^{-K'}$. We can then use Cauchy-Schwarz and Proposition \ref{P:RR} to conclude that \[||E1_{\tau_{\tilde{K}} \cap \Omega_{\delta}(K)}E1_{{\tau'} \cap \Omega_{\delta}(K')}||_{L^s} \lesssim ||E1_{{\tilde{K}} \cap \Omega_{\delta}(K)}||_{L^{2s}}||E1_{{\tau'} \cap \Omega_{\delta}(K')}||_{L^{2s}} \lesssim |\tau_{\tilde{K}}|^{1/s'}|\tau'|^{1/s'} \lesssim 2^{-\Theta_s|K'-K|}|\Omega(K)|^{2/s'}.\]This is a good enough bound for (\ref{E:S32 easy1}). We can follow the same reasoning and the fact that $|J - J'| \sim |K - K'|$ to prove that \[||E1_{{\tau\setminus \tau_{\tilde{K}}} \cap \Omega_{\delta}(K)}E1_{\tau'_{\tilde{J}} \cap \Omega_{\delta}(K')}||^s_{L^s} \lesssim 2^{-\Theta_s|J - J'|}|\Omega(K)|^{2/s'} \le 2^{-\Theta_s|K - K'|}|\Omega(K)|^{2/s'}.\]Hence, we have the right bound for (\ref{E:S32 easy2}) as well. It it remains to prove that
    \begin{equation}\label{E:S32 hard reasoning}
     \sum_{j = J'}^{\tilde{J}}\sum_{k = K}^{\tilde{K}} ||E1_{\tau_k \cap \Omega_{\delta}(K)}E1_{\tau'_j \cap \Omega_{\delta}(K')}||_{L^s} \lesssim 2^{-\Theta_s|K'-K|}|\Omega(K)|^{2/s'}
    \end{equation}

    Our bilinear extension estimate extends to hold on pairs of tiles $\tau_k, \tilde{\tau}_j$.

    \begin{lemma}\label{L:bil ext holds on weird tiles}
      Suppose $k < \tilde{K}$ and $j < \tilde{J}$. Then for $r$ in a neighborhood of $s'$, we have
      \begin{equation}\label{E:bil ext holds on weird tiles}
          \sum_{j = J'}^{\tilde{J}}\sum_{k = K}^{\tilde{K}} ||E1_{\tau_k \cap \Omega_{\delta}(K)}E1_{\tau'_j \cap \Omega_{\delta}(K')}||_{L^s} \lesssim 2^{(j+k)\left(\frac2s + \frac2r -2\right)}|\tau_k|^{1/r}|\tau'_j|^{1/r}.
        \end{equation}
    \end{lemma}

    The proof of Lemma \ref{L:bil ext holds on weird tiles} is simple but requires some book-keeping, so we defer it to the end of the proof of Lemma \ref{L:S32}. Since $k < \tilde{K}, j < \tilde{J}, |\tau_k| < |\tau|,$ and $|\tau'_j| < |\tau'|$ (\ref{E:bil ext holds on weird tiles}) implies 
    
    \begin{equation}\label{E:chill bounds}
      ||E1_{\tau_k \cap \Omega_{\delta}(K)}E1_{\tau'_j \cap \Omega_{\delta}(K')}||_{L^s} \lesssim 2^{2(\tilde{J} + \tilde{K})\left(\frac{1}{r} - \frac{1}{s'}\right)}|\tau|^{1/r}|\tau'|^{1/r}.
    \end{equation}We substitute (\ref{E:chill bounds}) into (\ref{E:S32 hard reasoning}) to bound \begin{align*}\sum_{j = J'}^{\tilde{J}}\sum_{k = K}^{\tilde{K}} &||E1_{\tau_k \cap \Omega_{\delta}(K)}E1_{\tau'_j \cap \Omega_{\delta}(K')}||_{L^s} \lesssim |J' - \tilde{J}||K - \tilde{K}|2^{(2\tilde{J} + 2\tilde{K})\left(\frac{1}{r} - \frac{1}{s'}\right)}(|\tau||\tau'|)^{1/r}  \\ &\lesssim |K - K'|^2 2^{(2\tilde{J} + 2\tilde{K} - J - J' - K - K')\left(\frac{1}{r} - \frac{1}{s'}\right)}(|\tau||\tau'|)^{1/s'}\\ & \le |K - K'|^2 2^{(2(\tilde{J} - J) + 2(\tilde{K} - K') - (J' - J) - (K - K'))\left(\frac{1}{r} - \frac{1}{s'}\right)}|\Omega(K)|^{2/s'} \\ &\le |K - K'|^22^{\frac{3}{5}\left(\frac{1}{r} - \frac{1}{s'}\right)(-|J' - J| - |K'-K|)} |\Omega(K)|^{2/s'} \tag{$\ast$}\\ &= |K - K'|^2 2^{-\Theta_s(|K-K'|)}|\Omega(K)|^{2/s'} \lesssim 2^{-\Theta_s(|K-K'|)}|\Omega(K)|^{2/s'}.\end{align*}In the inequality labeled ($\ast$), we use that $\tilde{J} - J = \frac{4}{5}(J' - J)$ and $J'< J$ so $2(\tilde{J} - J) - (J' - J) = -\frac{3}{5}|J' - J|$. Similarly for $K, \tilde{K}$, and $K'$. Finally, we have the right bound for (\ref{E:S32 hard}). All that remains is a proof of Lemma \ref{L:bil ext holds on weird tiles}. 

    \begin{proof}[Proof of Lemma \ref{L:bil ext holds on weird tiles}]
      Both $\tau_k, \tau'_j$ comprise at most two connected components. For example, $\tau_k = \left(\pi_1^{-1}(I_J) \cap \pi_2^{-1}(I^+_k)\right) \cup  \left(\pi_1^{-1}(I_J) \cap \pi_2^{-1}(I^-_k)\right)$, where $I^+_k, I_k^-$ are (possibly empty) intervals of length $\lesssim 2^{-k}$ and separated from $\pi_2(\tau'_j)$ by $\sim 2^{-k}$. We have a similar decomposition for $\tau'_j$. If we can prove the desired bound for each pair of connected components, then we have the estimate for $\tau_k, \tau'_j$ as well by summing over the $\le 4$ combinations of connected components. In a small abuse of notation, we simply assume $\tau_k, \tau'_j$ are each a single component and hence of the form $\tau_k = \pi_1^{-1}(I_J) \cap \pi_2^{-1}(I^*_k), \tau_j = \pi_1^{-1}(I^*_j) \cap \pi_2^{-1}(I_{K'})$. Here, $I_J$ and $I^*_j$ are $\sim 2^{-j}$ separated, $I^*_j$ is an interval with length $\lesssim 2^{-j}$, $I_{K'}$ and $I^*_k$ are $\sim 2^{-k}$ separated and $I_{K'}$ is an interval with length $\lesssim 2^{-k}$.
      
      Without loss of generality, assume $k > j$. We claim that any $(\xi, \eta) \in \tau_k \times \tau'_j$ and $k' \le k$ such that $(\xi, \eta)$ is $(2^{-k'}, 2^{-j})$-separated. Since $\pi_1$ is $\sim 1$ Lipschitz, our assumption that $\pi_1(\tau_k) = I_J$ and $\pi_1(\tau'_j) = I^*_j$ are $\sim 2^{-j}$ separated implies that $|\xi - \eta| \sim 2^{-j}$, justifying the second entry of the separation tuple. As for the first, suppose for contradiction that $\sup_{t \in [-1,2]} |N_2(\xi + t(\eta - \xi)) \wedge (\xi - \eta)| \ll 2^{-k}$. Then \[|\pi_2(\eta) - \pi_2(\xi)| \lesssim \int_0^1 |\langle d\pi_2(\xi + t(\eta - \xi)), (\eta-\xi)\rangle| ~dt = \int_0^1 |N_2(\xi + t(\eta-\xi)) \wedge (\eta-\xi)| ~dt \ll 2^{-k}.\]Since we assumed $|\pi_2(\eta) - \pi_2(\xi)| \gtrsim 2^{-k}$, this is indeed a contradiction, so we see that $\sup_{t \in [-1,2]} |N_2(\xi + t(\eta - \xi)) \wedge (\xi - \eta)| \sim 2^{-k'} \gtrsim 2^{-k}$. We have the desired separation for the points $(\xi, \eta)$. 

      For any $\xi$, if there exists $\eta$ so that $(\xi, \eta)$ is $(2^{-k'}, 2^{-j})$-separated, then we claim that $k^2(\xi, j) \ge k'$. We define $k'$ so that the direction of $N_2(\xi + t(\eta - \xi))$ and $\eta - \xi$ is less than $2^{j-k'}$ for $t \in [0,1]$. A consequence is that $\xi + t(\eta - \xi)$ stays within $2^{-k'}$ of $\xi + 2^{-j}t N_2(\xi)$ for $t \in [0,1]$, which implies that $N_2(\xi + 2^{-j}t N_2(\xi))$ also stays within $2^{-k'}$ of $\eta - \xi$ for all $t$ and in particular $t = 0$. It follows that, \[\sup_{t \in [0,1]}|N_2(\xi + 2^{-j}t N_2(\xi)) \wedge N_2(\xi)| \lesssim 2^{-k'}\]and hence $k^2(\xi, k) \ge k'$.

      We can partition $\tau_k$ into $\lesssim 1$ many sets in $\mathcal{T}_{j,k}$ by partitioning $I_J, I^+_k$ with $\lesssim 1$ many elements of $\mathcal{I}_{j+C}, \mathcal{I}_{k+C}$. We call this collection of sets $\mathcal{R}[\tau_k]$. We similarly define $\mathcal{R}[\tau'_j]$. For each pair $(R_0, R'_0) \in  \mathcal{R}[\tau_k] \times \mathcal{R}[\tau'_j]$, let $\xi, \eta$ respectively be the center points of $R_0, R'_0$ (that is, $R_0 = R_3(\xi, k+C, j+C)$ and $R'_0 = R_3(\eta, k+C, j+C)$). We know that $(\xi, \eta)$ is $(2^{-k'}, 2^{-j})$-separated for some $k' \le k(\xi, j)$. The sets $R = R_3(\xi, k'+C, j-C), R' = R_3(\eta, k'+C, j+C)$ satisfy the assumptions of Corollary \ref{C:rescaled bilinear rest} with separation $(2^{-k'}, 2^{-j})$, so we can bound

      \begin{equation}\label{E: S32 bil ext appl 1}
        ||Ef_R Eg_{R'}||_{L^s} \lesssim 2^{(j+k')\left(\frac2r + \frac2s -2\right)}||f_R||_{L^r}||g_{R'}||_{L^r} \le 2^{(j+k)\left(\frac2r + \frac2s -2\right)}||f_R||_{L^r}||g_{R'}||_{L^r}
      \end{equation}
      for any $f_R$ supported on $R$ and any $g_{R'}$ supported on $R'$. Since $R_0 \subset R, R'_0 \subset R'$, we can take $f_R = 1_{R_0 \cap \tau_k \cap \Omega_{\delta}(K)}, g_{R'} = 1_{R'_0 \cap \tau'_j \cap \Omega_{\delta}(K')}$ in (\ref{E: S32 bil ext appl 1}) and conclude \begin{align}
        ||E1_{R_0 \cap \tau_k \cap \Omega_{\delta}(K)} E1_{R'_0 \cap \tau'_j \cap \Omega_{\delta}(K')}||_{L^s} &\lesssim 2^{(j+k)\left(\frac2r + \frac2s -2\right)}|R_0 \cap \tau_k \cap \Omega_{\delta}(K)|^{1/r}|R'_0 \cap \tau'_j \cap \Omega_{\delta}(K')|^{1/r} \notag \\ &\le 2^{(j+k)\left(\frac2r + \frac2s -2\right)}|\tau_k|^{1/r}|\tau'_j|^{1/r}.\label{E: S32 bil ext appl 2}
      \end{align}
      Finally, we sum over the $\lesssim 1$ many pairs $(R_0, R'_0) \in \mathcal{R}[\tau_k] \times \mathcal{R}[\tau'_j]$ to conclude 

      \begin{align}
        ||E1_{\tau_k \cap \Omega_{\delta}(K)} E1_{\tau'_j \cap \Omega_{\delta}(K')}||_{L^s} &\le \sum_{(R_0, R'_0) \in \mathcal{R}[\tau_k] \times \mathcal{R}[\tau'_j]} ||E1_{R_0 \cap \tau_k \cap \Omega_{\delta}(K)} E1_{R'_0 \cap \tau'_j \cap \Omega_{\delta}(K')}||_{L^s}\\ &\lesssim  |\mathcal{R}[\tau_k] \times \mathcal{R}[\tau'_j]| 2^{(j+k)\left(\frac2r + \frac2s -2\right)}|\tau_k|^{1/r}|\tau'_j|^{1/r}\\ &\lesssim 2^{(j+k)\left(\frac2r + \frac2s -2\right)}|\tau_k|^{1/r}|\tau'_j|^{1/r}.
      \end{align}
      Hence, we have (\ref{E:bil ext holds on weird tiles}). 
    \end{proof}
  \end{proof}
   With a good off-diagonal estimate from Lemma \ref{L:S32}, the proof of Proposition \ref{P:S31} is straightforward. 
  \begin{proof}[Proof of Proposition \ref{P:S31}]
    Let $A$ be an admissible constant much larger than $C$. It suffices to prove \begin{equation}\label{E:S31new goal}
      \left|\left| \sum_{K \in \mathcal{K}} E1_{\Omega_{\delta}(K)}\right|\right|_{L^{2s}}^{2s} \lesssim \sum_{K \in \mathcal{K}} ||E1_{\Omega_{\delta}(K)}||_{L^{2s}}^{2s} + \delta^2 |\Omega|^{\frac{2s}{s'}}.
    \end{equation}for any $A |\log \delta|$ separated set $\mathcal{K} \subset \mathcal{K}(\varepsilon)$. We can also assume $\{J(K) : K \in \mathcal{K}\}$ is $A |\log \delta|$ separated when proving (\ref{E:S31new goal}), as removing the sparsity requirements on $J$ and $K$ will only make (\ref{E:S31new goal}) worse by a polylogarithmic factor acceptable in (\ref{E:S31 goal}).

    We write \[\left(\sum_{K \in \mathcal{K}} E1_{\Omega_{\delta}(K)}\right)^{2s} = \left(\sum_{K_1, K_2, K_3, K_4 \in \mathcal{K}} \prod_{i=1}^4 E1_{\Omega_{\delta}(K_i)}\right)^{s/2}.\]We separate out the diagonal and off-diagonal terms \[\left(\sum_{K_1, K_2, K_3, K_4} \prod_{i=1}^4 E1_{\Omega_{\delta}(K_i)}\right)^{s/2} \lesssim \left(\sum_{K \in \mathcal{K}}|E1_{\Omega_{\delta}(K)}|^4\right)^{\frac{s}{2}} + \left(\sum_{\substack{K_1, K_2, K_3, K_4\\ \text{not all the same}}} \prod_{i=1}^4 E1_{\Omega_{\delta}(K_i)}\right)^{s/2}.\]
    Since $s < 2$, we bound by convexity both the diagonal and off-diagonal terms: \begin{align*}
      \left(\sum_{K \in \mathcal{K}}|E1_{\Omega_{\delta}(K)}|^4\right)^{\frac{s}{2}} &\lesssim \sum_{K \in \mathcal{K}}|E1_{\Omega_{\delta}(K)}|^{2s}\\
      \left(\sum_{\substack{K_1, K_2, K_3, K_4\\ \text{not all the same}}} \prod_{i=1}^4 E1_{\Omega_{\delta}(K_i)}\right)^{s/2} &\lesssim \sum_{\substack{K_1, K_2, K_3, K_4\\ \text{not all the same}}} \prod_{i=1}^4 |E1_{\Omega_{\delta}(K_i)}|^{s/2} 
    \end{align*}
    Integrating, we conclude that \begin{equation*}\left|\left| \sum_{K \in \mathcal{K}} E1_{\Omega_{\delta}(K)}\right|\right|_{L^{2s}}^{2s} \lesssim \sum_{K \in \mathcal{K}}||E1_{\Omega_{\delta}(K)}||_{L^{2s}}^{2s} + \sum_{\substack{K_1, K_2, K_3, K_4\\ \text{not all the same}}} \left|\left|\prod_{i=1}^4 E1_{\Omega_{\delta}(K_i)}\right|\right|_{L^{s/2}}^{s/2}.\end{equation*}

    It remains to prove \begin{equation}\label{E:S31 intermediate 2}\sum_{\substack{K_1, K_2, K_3, K_4\\ \text{not all the same}}} \left|\left|\prod_{i=1}^4 E1_{\Omega_{\delta}(K_i)}\right|\right|_{L^{s/2}}^{s/2} \lesssim \delta^2 |\Omega|^{\frac{2s}{s'}}.\end{equation}We may relabel the indices so that \begin{equation}\label{E:relabling the indices ig}
        K_1 \neq K_2 \text{ and } \Omega(K_1)| \ge \max\{|\Omega(K_2)|, |\Omega(K_3)|, |\Omega(K_4)|\
    \end{equation}for each choice of $K_1, K_2, K_3, K_4$ contributing to (\ref{E:S31 intermediate 2}). We use H\"older with Proposition \ref{P:RR} to bound \begin{equation}\label{E:S31 intermediate}\sum_{\substack{K_1, K_2, K_3, K_4\\ \text{not all the same}}} \left|\left|\prod_{i=1}^4 E1_{\Omega_{\delta}(K_i)}\right|\right|_{L^{s/2}}^{s/2} \lesssim \sum_{K_1 \in \mathcal{K}} \sum_{K_2 \neq K_1 \in \mathcal{K}}\left|\left|E1_{\Omega_{\delta}(K_1)}E1_{\Omega_{\delta}(K_2)}\right|\right|_{L^{s}}^{s/2}|\Omega(K_1)|^{\frac{s}{s'}}.\end{equation}
    By Lemma \ref{L:S23} and our assumption $|\Omega(K_1)| \ge |\Omega(K_2)|$, we can bound \begin{equation}\label{E:bil intermediate bound}\left|\left|E1_{\Omega_{\delta}(K_1) \cap \tau}E1_{\Omega_{\delta}(K_2) \cap \tau'}\right|\right|_{L^{s}}^{s/2} \lesssim 2^{-O_s(|K_1 - K_2|)} |\Omega(K_1)|^{\frac{s}{s'}}\end{equation}whenever $K_1 < K_2$ and $J_1 > J_2$ or $K_1 > K_2$ and $J_1 < J_2$ and $\tau \in \mathcal{R}_{\delta}(K_1), \tau' \in \mathcal{R}_{\delta}(K_2)$. By Proposition \ref{P:S22}, $|\mathcal{R}_{\delta}(K_1)|, |\mathcal{R}_{\delta}(K_2)| \lesssim \delta^{-O(1)}$ and the tiles from $\mathcal{R}_{\delta}(K_1), \mathcal{R}_{\delta}(K_2)$ respectively cover $\Omega_{\delta}(K_1), \Omega_{\delta}(K_2)$. We sum (\ref{E:bil intermediate bound}) over choices of $\tau \in  \mathcal{R}_{\delta}(K_1), \tau' \in \mathcal{R}_{\delta}(K_2)$ to conclude \begin{equation}\label{E:bil final bound}\left|\left|E1_{\Omega_{\delta}(K_1)}E1_{\Omega_{\delta}(K_2)}\right|\right|_{L^{s}}^{s/2} \lesssim \delta^{-O(1)} 2^{-O_s(|K_1 - K_2|)} |\Omega(K_1)|^{\frac{s}{s'}}\end{equation}
    If $J_1 \le J_2$ and $K_1< K_2$, we can use Cauchy-Schwarz and Proposition \ref{P:RR} to bound \[\left|\left|E1_{\Omega_{\delta}(K_1)}E1_{\Omega_{\delta}(K_2)}\right|\right|_{L^{s}} \lesssim |\Omega(K)|^{\frac{1}{s'}}|\Omega(K')|^{\frac1{s'}} \lesssim 2^{-\Theta_s|K - K'|}|\Omega|(K)|^{\frac2{s'}}\]and conclude (\ref{E:bil final bound}), indeed without the $\delta^{-O(1)}$ error. The only remaining cases are $K_1 = K_2$ or $J_1> J_2, K_1 > K_2$, but both contradict (\ref{E:relabling the indices ig}). Plugging (\ref{E:bil final bound}) into (\ref{E:S31 intermediate}) and recalling the $K_i$ are $A |\log \delta|$ separated, we have \begin{align*}\sum_{\substack{K_1, K_2, K_3, K_4\\ \text{not all the same}}} \left|\left|\prod_{i=1}^4 E1_{\Omega_{\delta}(K_i)}\right|\right|_{L^{s/2}}^{s/2} &\lesssim \delta^{-O(1)}\sum_{K_1 \in \mathcal{K}} \sum_{K_2 \neq K_1} 2^{-O_s(|K_1 - K_2|)}|\Omega(K_1)|^{\frac{2s}{s'}}\\
    & \lesssim \sum_{K_1 \in \mathcal{K}} \delta^{O_s(A)-O(1)}|\Omega(K_1)|^{\frac{2s}{s'}} \\&\lesssim \delta^{O_s(A)-O(1)} \sum_{K_1 \in \mathcal{K}} |\Omega(K_1)|^{\frac{2s}{s'}}.\end{align*}
    Since $2s > s'$, \[\sum_{\substack{K_1, K_2, K_3, K_4\\ \text{not all the same}}}\left|\left|\prod_{i=1}^4 E1_{\Omega_{\delta}(K_i)}\right|\right|_{L^{s/2}}^{s/2} \lesssim \delta^{O_s(A)-O(1)} \sum_{K_1 \in \mathcal{K}} |\Omega(K_1)|^{\frac{2s}{s'}} \lesssim \delta^{O_s(A)-O(1)} |\Omega|^{\frac{2s}{s'}}.\]Making a large enough choice of $A$ gives (\ref{E:S31 intermediate 2}).
  \end{proof}

  We are finally ready to prove our scaling line estimate. We decompose $\Omega$ into subsets with different fiber lengths and further decompose the fiber lengths depending on the sharpness of Proposition \ref{P:RR} when applied to $\Omega(K)$. 

  \begin{proof}[Proof of Proposition \ref{P:Bil2RST}]
    We decompose \[\Omega = \bigsqcup_{0 < \varepsilon \lesssim 1} \bigsqcup_{0 < \delta \le \varepsilon} \left(\bigcup_{K \in \mathcal{K}(\varepsilon)} \Omega_{\delta}(K)\right),\]where $\mathcal{K}(\varepsilon)$ is as defined in Definition \ref{D:Keps}. Using the triangle inequality and (\ref{P:S31}), we can bound \begin{align}\left|\left|E1_{\Omega}\right|\right|_{L^{2s}} &\lesssim  \sum_{0 < \varepsilon \lesssim 1} \sum_{0 < \delta \le \varepsilon}\left|\left|\sum_{K \in \mathcal{K}(\varepsilon)} E1_{\Omega_{\delta}(K)}\right|\right|_{L^{2s}}\\ &\lesssim \sum_{0 < \varepsilon \lesssim 1} \sum_{0 < \delta \le \varepsilon}|\log \delta|^{O(1)}\left(\left(\sum_{K \in \mathcal{K}(\varepsilon)} ||E1_{\Omega_{\delta}(K)}||^{2s}_{L^{2s}}\right)^{\frac{1}{2s}} + \delta^{\frac{1}{2s}} |\Omega|^{\frac{1}{s'}}\right)\label{E:mixed bound}.\end{align}We use conclusion (3) of Proposition \ref{P:S22} to bound $||E1_{\Omega_{\delta}(K)}||_{L^{2s}} \lesssim \delta|\Omega(K)|^{\frac{1}{s'}}$. Since $\frac{2s}{s'} > 1$, we have \[\sum_{K \in \mathcal{K}(\varepsilon)} ||E1_{\Omega_{\delta}(K)}||^{2s}_{L^{2s}} \lesssim \delta^{2s} \sum_{K \in \mathcal{K}(\varepsilon)} |\Omega(K)|^{\frac{2s}{s'}} \lesssim \delta^{2s}|\Omega|^{\frac{2s}{s'}}.\]Plugging this into (\ref{E:mixed bound}) gives \[||E1_{\Omega}||_{L^{2s}} \lesssim \sum_{0 < \varepsilon \lesssim 1} \sum_{0 < \delta \le \varepsilon} |\log \delta|^{O(1)} \delta^{\frac{1}{2s}}|\Omega|^{\frac{1}{s'}} \lesssim |\Omega|^{\frac{1}{s'}}.\]
  \end{proof}

\printbibliography
\end{document}